\documentclass{article}

\usepackage[letterpaper, margin=1in]{geometry}
\usepackage{amsmath,amssymb}
\usepackage{mathrsfs}
\usepackage{graphicx}
\usepackage[amsmath,thmmarks,framed]{ntheorem}
\usepackage{authblk}
\usepackage{xcolor}
\definecolor{darkgreen}{rgb}{0.0,0.5,0.0}
\usepackage[colorlinks,linkcolor=red,citecolor=darkgreen,urlcolor=blue]{hyperref}
\hypersetup{
  pdftitle={Semigroup Hessian estimates for underdamped Langevin dynamics},
  pdfauthor={Xuda Ye}
}
\usepackage[numbers,square,comma,sort&compress]{natbib}

\newtheorem{theorem}{Theorem}
\newtheorem{lemma}{Lemma}
\newtheorem{corollary}[theorem]{Corollary}

\newtheorem{assumption}{Assumption}
\newtheorem{definition}{Definition}
\newenvironment{proof}[1][Proof]{\par\medskip\noindent\textit{#1.}\ }{\hfill$\square$\par\medskip}

\newcommand{\Le}{\leqslant}
\newcommand{\Ge}{\geqslant}
\newcommand{\D}{\mathrm{d}}
\newcommand{\norm}[1]{\lVert #1 \rVert}

\title{Semigroup Hessian Estimates for Underdamped Langevin Dynamics}
\author[1]{Xuda Ye}
\affil[1]{Department of Mathematics, Purdue University, West Lafayette, IN 47907, USA}
\author[2]{Zhihan Zhang}
\affil[2]{Department of Computer Science and Operations Research, Université de Montréal, Montréal, QC H3T 1N8, Canada}
\date{}

\numberwithin{equation}{section}

\begin{document}
\maketitle

\begin{abstract}
For the semigroup $e^{t\mathcal L}$ of the underdamped Langevin dynamics, we bound the Hessian $\nabla^2(e^{t\mathcal L}f)$ when $\nabla f$ is uniformly bounded. When the potential is strongly convex outside a ball and has bounded second and third derivatives, and the damping rate is above an explicit threshold, the bounds on $\nabla^2(e^{t\mathcal L}f)$ decay exponentially in time with constants free of the dimension. The proof combines the Wasserstein contraction of the dynamics with a Bismut--Elworthy--Li formula for the degenerate noise. As an application, we prove the uniform-in-time error of the stochastic gradient UBU integrator, with an explicit coefficient of the numerical bias in the step size, proportional to the covariance of the stochastic gradient noise.


\end{abstract}

\section{Introduction}
\label{sec: intro}

Let $U$ be a potential on $\mathbb R^d$, let $\gamma > 0$ be the damping rate and consider the underdamped Langevin dynamics
\begin{equation}
    \D X_t = V_t \,\D t , \qquad
    \D V_t = -\bigl( \nabla U(X_t) + \gamma V_t \bigr) \D t + \sqrt{2\gamma}\,\D B_t ,
    \label{eq: uld}
\end{equation}
with $(B_t)_{t \Ge 0}$ a standard Brownian motion in $\mathbb R^d$. The dynamics leaves invariant the Boltzmann distribution of the state $z = (x,v) \in \mathbb R^{2d}$,
\begin{equation}
    \pi(\D z) \propto \exp\Bigl( -U(x) - \frac12 |v|^2 \Bigr) \D z .
    \label{eq: boltzmann}
\end{equation}
Write $Z_t = (X_t,V_t)$ for the solution, $\pi(f) := \int f \,\D\pi$ for the mean of a function $f$ on $\mathbb R^{2d}$ under $\pi$, $\nabla^2_{xx}f$, $\nabla^2_{xv}f$, $\nabla^2_{vv}f\in\mathbb R^{d\times d}$ for the blocks of the Hessian of $f$, and $\nabla^3_{xvv}f$, $\nabla^3_{vvv}f\in\mathbb R^{d\times d\times d}$ for the derivatives of $\nabla^2_{vv}f$ in the position $x$ and in the velocity $v$. The generator of the dynamics is
\begin{equation}
    \mathcal L f = v \cdot \nabla_x f - \bigl( \nabla U(x) + \gamma v \bigr) \cdot \nabla_v f + \gamma \Delta_v f ,
    \label{eq: generator}
\end{equation}
and $(e^{t\mathcal L})_{t\Ge0}$ is the associated semigroup. By definition, $(e^{t\mathcal L} f)(z) = \mathbb E^z\bigl[ f(Z_t) \bigr]$, where $\mathbb E^z$ is the expectation under the law of the solution started at $z$.

This paper bounds the Hessian of $e^{t\mathcal L}f$ with $\nabla f$ alone, and these bounds decay exponentially in time. For a potential with bounded second and third derivatives that is strongly convex outside a ball, a damping rate $\gamma$ above an explicit threshold, and every $f \in C^1(\mathbb R^{2d})$ with $\norm{\nabla f}_\infty <\infty$, Theorem~\ref{thm: hessian} and Corollary~\ref{cor: uniform} prove the Hessian bound
\begin{equation*}
\begin{aligned}
    \norm{\nabla^2_{vv} (e^{t\mathcal L} f)}_\infty + \norm{\nabla^2_{xv} (e^{t\mathcal L} f)}_\infty & \Le C \min(1,t)^{-\frac12} e^{-\lambda t} \norm{\nabla f}_\infty ,
    \qquad\forall t>0, \\
    \norm{\nabla^2_{xx} (e^{t\mathcal L} f)}_\infty & \Le C \min(1,t)^{-\frac32} e^{-\lambda t} \norm{\nabla f}_\infty , \qquad\forall t>0,
\end{aligned}
\end{equation*}
and Theorem~\ref{thm: third} and Corollary~\ref{cor: third uniform} produce the third derivative bound
\begin{equation*}
    \norm{\nabla^3_{xvv} (e^{t\mathcal L} f)}_\infty + \norm{\nabla^3_{vvv} (e^{t\mathcal L} f)}_\infty \Le C \min(1,t)^{-1} e^{-\lambda t} \norm{\nabla f}_\infty , \qquad\forall t>0,
\end{equation*}
where $\norm{\,\cdot\,}_\infty$ is the supremum of the tensor norm in the Euclidean distance, which is the absolute value for a scalar, the Euclidean norm for a vector and the matrix 2-norm for a matrix. The constant $C$ and the rate $\lambda > 0$ do not depend on the dimension $d$. The exponents of $t$ are sharp.

The error analysis of numerical integrators of the dynamics relies on the Poisson equation and on bounds on the derivatives of its solution. Define the function
\begin{equation}
    \phi(z) := \int_0^\infty \bigl( (e^{t\mathcal L}f)(z) - \pi(f) \bigr) \D t ,
    \label{eq: poisson solution}
\end{equation}
which, by the Hessian bound, is a strong solution to the Poisson equation
\begin{equation*}
    \mathcal L\phi = v \cdot \nabla_x \phi - \bigl( \nabla U(x) + \gamma v \bigr) \cdot \nabla_v \phi + \gamma \Delta_v \phi = \pi(f) - f .
\end{equation*}

As an application of these semigroup Hessian estimates, we prove in Theorem~\ref{thm: sgubu} the uniform-in-time error of the stochastic gradient UBU (SG-UBU) integrator \citep{sanz2021wasserstein, lu2025mse}. Corollary~\ref{cor: bias limit} shows that the numerical invariant distribution with the step size $h$, denoted by $\pi_h$, has an explicit bias estimate
\begin{equation*}
    \Bigl| \pi_h(f) - \pi(f) - \frac h2\, \pi\bigl( \nabla^2_{vv}\phi : \mathcal C \bigr) \Bigr| \Le C \norm{\nabla f}_\infty d^2 h^2 \Bigl( 1 + \log\frac1h \Bigr) ,
\end{equation*}
where $\mathcal C(x)$ is the covariance matrix of the stochastic gradient noise, $A : B := \sum_{i,j=1}^d A_{ij} B_{ij}$ for matrices $A, B \in \mathbb R^{d\times d}$, and $C$ also depends on the regularity of the stochastic gradients. As a consequence, the numerical bias has an explicit coefficient,
\begin{equation*}
    \lim_{h \to 0} \frac{\pi_h(f) - \pi(f)}{h} = \frac12\, \pi\bigl( \nabla^2_{vv}\phi : \mathcal C \bigr) ,
\end{equation*}
which depends linearly on the covariance matrix $\mathcal C(x)$. Corollary~\ref{cor: bias limit} sharpens \citep[Theorem~9]{lu2025mse} in the sense that it proves the limit without global convexity of the potential, and only $\norm{\nabla f}_\infty$ enters the bound. Corollary~\ref{lem: w1 limit} also gives the explicit limit of $\mathcal W_1(\pi_h,\pi)/h$ as $h \to 0$, with the Wasserstein-1 distance on $\mathbb R^{2d}$.

The Hessian estimates rest on two ingredients. The decay in time comes from the Wasserstein contraction \citep{schuh2024global}, which gives $\norm{\nabla (e^{t\mathcal L}f)}_\infty \Le C e^{-\lambda t}\norm{\nabla f}_\infty$ when $\gamma$ is large. The short-time singularities come from a derivative formula of Bismut--Elworthy--Li (BEL) type \citep{bismut1984large, elworthy1994formulae} for the degenerate noise. BEL formulas and integration by parts formulas hold for classes of degenerate diffusions, including the kinetic case \citep{priola2006formulae, zhang2010stochastic, guillin2012degenerate, wang2013derivative, wang2014integration}. For bounded test functions, the formula of \citet{guillin2012degenerate} gives
\begin{equation*}
    \norm{\nabla_x (e^{t\mathcal L} f)}_\infty \Le C t^{-\frac32} \norm{f}_\infty , \qquad
    \norm{\nabla_v (e^{t\mathcal L} f)}_\infty \Le C t^{-\frac12} \norm{f}_\infty ,\qquad \forall t \in (0,1],
\end{equation*}
while this paper applies it to test functions with bounded gradient. For uniformly elliptic diffusions, the decay in a Wasserstein distance and a BEL formula bound the Hessian of the semigroup by the Lipschitz constant of the test function \citep{gorham2019measuring}, and the same two ingredients estimate the second and third derivatives of the semigroup of the overdamped Langevin dynamics \citep{ye2026mse}.

The paper is organized as follows. Section~\ref{sec: related} reviews related works. Section~\ref{sec: pre} fixes the notation and the assumptions. Section~\ref{sec: semigroup} states the Hessian and third derivative estimates and solves the Poisson equation, and Section~\ref{sec: sgubu} states the uniform-in-time error of SG-UBU and the coefficient of its numerical bias. Section~\ref{sec: conclusion} concludes. Appendix~\ref{app: proofs} proves the results of Section~\ref{sec: semigroup}, and Appendix~\ref{app: sgubu} those of Section~\ref{sec: sgubu}.

\section{Related works}
\label{sec: related}

This section reviews three lines of work: (i) entropy, hypoelliptic, Lyapunov and hypocoercivity methods for the convergence of the underdamped Langevin dynamics to equilibrium; (ii) coupling methods for the contraction of the dynamics and of its discretizations, and for the error of these discretizations; (iii) bounds on the derivatives of $e^{t\mathcal L}f$ and of the Poisson solution.

Entropy methods \citep{desvillettes2001trend} and hypoelliptic methods \citep{herau2004isotropic, helffer2005hypoelliptic} establish the trend to equilibrium of the kinetic Fokker--Planck equation. The spectral method of \citet{herau2004isotropic} extends to a wider class of hypoelliptic operators \citep{eckmann2003spectral}. Lyapunov functions give exponential convergence of the law of the dynamics to $\pi$ \citep{wu2001large, mattingly2002ergodicity}. The theory of hypocoercivity \citep{villani2009hypocoercivity} gives exponential convergence to equilibrium in a modified $H^1$ norm, and the $L^2$ hypocoercivity \citep{dolbeault2015hypocoercivity} gives it in a weighted $L^2$ norm. Generalized $\Gamma$ calculus \citep{baudoin2017bakry, monmarche2019generalized} casts hypocoercivity in the Bakry--\'Emery framework \citep{bakry1985diffusions}. Through the variational approach \citep{albritton2024variational}, \citet{cao2023explicit} obtain explicit rates of convergence in $L^2$ for the dynamics. Later works extend this approach to piecewise deterministic samplers \citep{lu2022explicit}, kinetic Fokker--Planck equations with fat-tailed local equilibria or weakly confining potentials \citep{brigati2023time, brigati2025construct, brigati2024explicit}, and non-reversible lifts of reversible diffusions \citep{eberle2026nonreversible, brigati2026hypocoercivity, eberle2025convergence}.

Coupling methods give contraction in Wasserstein distances. A synchronous coupling gives contraction for potentials close to a quadratic \citep{bolley2010trend} and for strongly convex potentials \citep{cheng2018underdamped, dalalyan2020sampling}. Reflection couplings give contraction for overdamped Langevin diffusions with locally nonconvex potentials \citep{eberle2016reflection} and, combined with synchronous couplings, for the underdamped Langevin dynamics \citep{eberle2019couplings}. When the potential is strongly convex outside a ball and the damping rate is above a threshold, \citet{schuh2024global} proves contraction with a rate free of the dimension. Couplings of these kinds also give contraction of discretizations of the dynamics \citep{monmarche2021high, leimkuhler2024contraction, schuh2024convergence, chak2026reflection}. \citet{schuh2024convergence} also bound the Wasserstein-1 distance between $\pi$ and the invariant distributions of the Euler--Maruyama integrator and of a splitting integrator. For dissipative potentials, \citet{chau2022stochastic} combine the contraction of \citet{eberle2019couplings} with error bounds on time intervals of unit length to bound, uniformly in time, the Wasserstein-1 error of stochastic gradient Hamiltonian Monte Carlo. Harnack inequalities bound the relative entropy between the laws of two solutions started at different states \citep{guillin2012degenerate}. When the damping rate is above a threshold, this relative entropy decays exponentially in time, for strongly convex potentials \citep{altschuler2025shifted} and for potentials that are strongly convex outside a ball \citep{habring2026reverse}. \citet{altschuler2025shifted} also bound the relative entropy between the laws of the numerical solution of an integrator and of the solution of the dynamics by the local error, with no contraction of the integrator, uniformly in time for strongly convex potentials and a damping rate proportional to the square root of the Lipschitz constant of $\nabla U$.

Bounds on the derivatives of $e^{t\mathcal L}f$ that are integrable in time control the derivatives of the solution $\phi$ of the Poisson equation. The Poisson equation underlies the error analysis of time averages and invariant distributions of numerical integrators \citep{mattingly2010convergence, kopec2015langevin, leimkuhler2016computation}. For uniformly elliptic diffusions, \citet{crisan2021uniform} use the exponential decay of the derivatives of the semigroup up to the fourth order to obtain a uniform-in-time weak error of first order for the Euler--Maruyama integrator. Discrete Poisson equations bound the mean-square error of time averages of integrators of the overdamped Langevin dynamics \citep{ye2026mse, ye2026sglm} and of the underdamped Langevin dynamics \citep{lu2025mse}. Hypocoercive resolvent estimates bound the Poisson solution in $L^2$ \citep{bernard2022hypocoercivity}. For diffusions with a constant diffusion coefficient and a strongly dissipative drift, a BEL formula bounds the Hessian of the Poisson solution by the Lipschitz constant of the test function \citep{fang2019multivariate}. \citet{talay2002stochastic} proves that $e^{t\mathcal L}f - \pi(f)$ and all its derivatives decay exponentially for smooth test functions. \citet{herau2007short} and \citet{villani2009hypocoercivity} bound the first derivatives of the solution of the kinetic Fokker--Planck equation in weighted $L^2$ norms, with the short-time factors $t^{-\frac12}$ in the velocity and $t^{-\frac32}$ in the position. Generalized $\Gamma$ calculus \citep{baudoin2017bakry} gives pointwise gradient bounds that decay exponentially when the damping rate is $1$ and the eigenvalues of $\nabla^2 U$ are uniformly bounded in a closed subinterval of $(0,1)$. For a class of degenerate operators that includes $\mathcal L$ when $\nabla^2 U$, $\nabla^3 U$ and $\nabla^4 U$ are bounded, \citet{priola2009global} bounds the derivatives of $e^{t\mathcal L}f$ up to the third order in the supremum norm by H\"older norms of a bounded test function $f$, with constants that depend on the dimension and do not decay in time. These derivative bounds for the dynamics assume a smooth or bounded test function, hold in weighted $L^2$ norms, or bound first derivatives only.

\section{Preliminaries}
\label{sec: pre}

\subsection{Notations}
\label{subsec: notation}

Throughout, $\langle \cdot, \cdot \rangle$ and $a \cdot b$ both denote the Euclidean inner product, and $|\cdot|$ the Euclidean norm of a vector and the absolute value of a scalar. The $k \times k$ identity matrix is $I_k$, and $A \succeq B$ means that $A - B$ is positive semidefinite. A tensor $T$ of order $k$ contracts with vectors $u_1, \dots, u_j$, $j \Le k$, in its last $j$ indices,
\begin{equation*}
    \bigl( T[u_1,\dots,u_j] \bigr)_{i_1 \dots i_{k-j}} = \sum_{l_1,\dots,l_j} T_{i_1 \dots i_{k-j} l_1 \dots l_j}\, (u_1)_{l_1} \cdots (u_j)_{l_j} .
\end{equation*}
For example, $A[u] = Au$ for a matrix $A$ and a vector $u$. For a tensor $T$ of order $k \Ge 1$ we write
\begin{equation}
    \norm{T} := \sup\bigl\{ |T[u_1,\dots,u_k]| : |u_1| \Le 1, \dots, |u_k| \Le 1 \bigr\} .
    \label{eq: tensor norm}
\end{equation}
In particular, $\norm{T}$ is
the Euclidean norm when $T$ is a vector and the matrix 2-norm when $T$ is a matrix.

A state $z$ consists of a position $x$ and a velocity $v$ in $\mathbb R^d$, and a direction $w$ consists of a position block $w^x$ and a velocity block $w^v$ in $\mathbb R^d$. The state and the direction are column vectors in $\mathbb R^{2d}$,
\begin{equation*}
    z = \begin{pmatrix} x \\ v \end{pmatrix}\in \mathbb R^{2d} , \qquad
    Z_t = \begin{pmatrix} X_t \\ V_t \end{pmatrix}\in \mathbb R^{2d} , \qquad
    w = \begin{pmatrix} w^x \\ w^v \end{pmatrix} \in \mathbb R^{2d},
\end{equation*}
and the position and velocity blocks of any vector in $\mathbb R^{2d}$ carry the superscripts $x$ and $v$ in the same way. Partial derivatives in the position and in the velocity are written $\nabla_x$ and $\nabla_v$, and $\Delta_v = \sum_{i=1}^{d}\partial_{v_i}^2$ is the Laplacian in the velocity.

A subscripted $\partial$ is a derivative in the initial state of the dynamics when its subscript is a direction, one of the blocks $x$ and $v$ of the initial state, or the initial state $z$ itself, as in $\partial_x X_t\in \mathbb R^{d\times d}$ and $\partial_z Z_t\in \mathbb R^{2d\times 2d}$. A derivative in the initial state carries the differentiated coordinate in its column, so that
\begin{equation*}
    \bigl( \partial_x X_t \bigr)_{ij} = \partial_{x_j}(X_t)_i , \qquad
    \partial_z Z_t = \begin{pmatrix} \partial_x X_t & \partial_v X_t \\ \partial_x V_t & \partial_v V_t \end{pmatrix},
\end{equation*}
and $ \partial_w Z_t = \bigl( \partial_z Z_t \bigr) w\in \mathbb R^{2d}$ is the derivative of $Z_t$ in the direction $w$. When its subscript is a coordinate, a subscripted $\partial$ is a derivative in that coordinate, of the argument or of the initial state, as in $\partial_{x_i}$ and $\partial_{x_j}(X_t)_i$. For a function of the time and the state, the operators $\nabla$, $\nabla^2$, $\nabla_x$, $\nabla_v$, $\Delta_v$ and $\mathcal L$ act in the state.

The blocks of the Hessian of a function $f$ on $\mathbb R^{2d}$ have the entries, with $i$ the row and $j$ the column,
\begin{equation*}
    \bigl( \nabla^2_{xx} f \bigr)_{ij} = \partial_{x_i}\partial_{x_j} f , \quad
    \bigl( \nabla^2_{xv} f \bigr)_{ij} = \partial_{x_i}\partial_{v_j} f , \quad
    \bigl( \nabla^2_{vx} f \bigr)_{ij} = \partial_{v_i}\partial_{x_j} f = \bigl( \nabla^2_{xv} f \bigr)_{ji} , \quad
    \bigl( \nabla^2_{vv} f \bigr)_{ij} = \partial_{v_i}\partial_{v_j} f .
\end{equation*}
For the sup norms we write
\begin{equation*}
    \norm{\nabla^j f}_\infty := \sup_{z\in\mathbb R^{2d}} \norm{\nabla^j f(z)} , \qquad
    \norm{\nabla_x f}_\infty := \sup_{z\in\mathbb R^{2d}} |\nabla_x f(z)| , \qquad
    \norm{\nabla^2_{xv} f}_\infty := \sup_{z\in\mathbb R^{2d}} \norm{\nabla^2_{xv} f(z)} ,
\end{equation*}
and likewise for $\nabla_v f$, $\nabla^2_{xx}f$ and $\nabla^2_{vv}f$.
For a differentiable $f$, $\nabla f(z)[w] = \nabla f(z) \cdot w$ is the derivative of $f$ at $z$ in $w$. For a twice differentiable $f$, $\nabla^2 f(z)[w,w'] = \langle \nabla^2 f(z)\, w, w' \rangle$ is its second derivative in $w$ and $w'$.

The Hessian of $e^{t\mathcal L}f$ contracted in the single direction $w$ is the vector
\begin{equation*}
    \nabla^2 (e^{t\mathcal L} f)(z)[w]
    = \begin{pmatrix} \nabla^2_{xx} (e^{t\mathcal L} f)(z) & \nabla^2_{xv} (e^{t\mathcal L} f)(z) \\ \nabla^2_{vx} (e^{t\mathcal L} f)(z) & \nabla^2_{vv} (e^{t\mathcal L} f)(z) \end{pmatrix}
    \begin{pmatrix} w^x \\ w^v \end{pmatrix} \in \mathbb R^{2d} .
\end{equation*}

For probability distributions $\nu$ and $\nu'$ with finite first moments on the same Euclidean space, the Wasserstein-1 distance is
\begin{equation*}
    \mathcal W_1(\nu,\nu') := \inf \bigl\{ \mathbb E|\xi - \xi'| : \mathrm{Law}(\xi) = \nu ,\ \mathrm{Law}(\xi') = \nu' \bigr\} ,
\end{equation*}
where the infimum runs over the pairs of random variables $(\xi,\xi')$ whose marginals are $\nu$ and $\nu'$. By Kantorovich--Rubinstein duality \citep{villani2021topics},
\begin{equation*}
    \mathcal W_1(\nu,\nu') = \sup \Bigl\{ \int f \,\D\nu - \int f \,\D\nu' : f \text{ has Lipschitz constant at most } 1 \Bigr\} .
\end{equation*}

\subsection{Assumptions}
\label{subsec: assumptions}

We first impose the regularity condition on the potential $U(x)$.

\begin{assumption}[regularity]
\label{ass: U}
    Let $M > 0$. The potential $U \in C^3(\mathbb R^d)$ satisfies
    \begin{equation}
        |\nabla U(0)| \Le M \sqrt d, \qquad \norm{\nabla^2 U(x)} \Le M , \qquad \norm{\nabla^3 U(x)} \Le M , \qquad \forall x \in \mathbb R^d .
        \label{eq: U bounds}
    \end{equation}
\end{assumption}

In \eqref{eq: U bounds}, the bound on $|\nabla U(0)|$ carries the factor $\sqrt d$, the root mean square of the norm of a vector distributed as $\mathcal N(0, I_d)$. The bounds on $\norm{\nabla^2 U(x)}$ and $\norm{\nabla^3 U(x)}$ carry no power of $d$, because these norms are suprema over unit directions. By \eqref{eq: U bounds} the gradient grows at most linearly,
\begin{equation}
    |\nabla U(x)| \Le M \bigl( \sqrt d + |x| \bigr) , \qquad \forall x \in \mathbb R^d .
    \label{eq: gradient growth}
\end{equation}

When long-time properties are concerned, we also require $U(x)$ to be strongly convex outside a bounded ball and the damping rate to be large.

\begin{assumption}[convexity]
\label{ass: convexity}
    Let $m > 0$ and $R \Ge 0$. The potential $U \in C^2(\mathbb R^d)$ satisfies
    \begin{equation}
        \nabla^2 U(x) \succeq m I_d , \qquad \forall\, |x| \Ge R .
        \label{eq: U convexity}
    \end{equation}
\end{assumption}

\begin{assumption}[damping]
\label{ass: damping}
    The damping rate $\gamma$ satisfies the damping condition
    \begin{equation}
        \gamma > \frac{2M+m}{\sqrt m} .
        \label{eq: friction condition}
    \end{equation}
\end{assumption}

The damping condition \eqref{eq: friction condition} ensures the condition on the damping rate in the Wasserstein contraction of \citet{schuh2024global}. From here on, the letter $C$ denotes a positive constant depending only on $m$, $M$, $R$ and $\gamma$, and it changes from line to line. In a statement under Assumption~\ref{ass: U} alone, $C$ depends only on $M$ and $\gamma$. In particular, no constant $C$ depends on the dimension $d$.

\section{Semigroup estimates for underdamped Langevin dynamics}
\label{sec: semigroup}

This section states the Hessian and third derivative estimates for $t \in (0,1]$ under Assumption~\ref{ass: U} and extends them to every $t > 0$ under Assumptions~\ref{ass: U}--\ref{ass: damping}. The estimates for every $t > 0$ then give a strong solution to the Poisson equation. Under Assumptions~\ref{ass: U} and~\ref{ass: convexity} alone, a Poincaré inequality gives the exponential decay of $e^{t\mathcal L}$ in $L^2(\pi)$ and a solution to the Poisson equation in $L^2(\pi)$.

\subsection{Variation processes}
\label{subsec: variations}

By \eqref{eq: U bounds}, the underdamped Langevin dynamics \eqref{eq: uld} has a unique strong solution with moments of every order, and $\bigl( X_t, V_t - \sqrt{2\gamma}B_t \bigr)$ solves an ordinary differential equation path by path. Let $Z_t = (X_t,V_t)$ be the solution with the initial state $z = (x,v)$. Then the map $z \mapsto Z_t$ is twice continuously differentiable, and the first and second derivatives of $Z_t$ in $z$ solve the equations obtained by differentiating \eqref{eq: uld} in $z$. The first variation $J_t := \partial_z Z_t \in \mathbb R^{2d\times2d}$ solves
\begin{equation}
    \frac{\D}{\D t} J_t = \begin{pmatrix} 0 & I_d \\ -\nabla^2U(X_t) & -\gamma I_d \end{pmatrix} J_t , \qquad J_0 = I_{2d} ,
    \label{eq: first variation}
\end{equation}
and the second variation $K_t[w] := \partial_w J_t \in \mathbb R^{2d\times2d}$, the derivative of $J_t$ in the initial state $z$ along the direction $w \in \mathbb R^{2d}$, solves
\begin{equation}
    \frac{\D}{\D t} K_t[w] = \begin{pmatrix} 0 & I_d \\ -\nabla^2U(X_t) & -\gamma I_d \end{pmatrix} K_t[w]  - \begin{pmatrix} 0 & 0 \\ \nabla^3U(X_t)\bigl[ (J_t w)^x \bigr] & 0 \end{pmatrix} J_t , \qquad K_0[w] = 0 ,
\label{eq: second variation}
\end{equation}
where $(J_tw)^x \in \mathbb R^d$ denotes the position block of $J_tw\in\mathbb R^{2d}$. The second variation $K_t$ is written with a direction $w$ mainly for notational convenience.

Lemma~\ref{lem: variations} bounds the variation processes $J_t$ and $K_t$ for $t\in[0,1]$, and these bounds yield the short-time estimate of $\nabla^2(e^{t\mathcal L}f)$ in Lemma~\ref{lem: derivatives}.

\begin{lemma}
\label{lem: variations}
    Under Assumption~\ref{ass: U}, for every initial state $z \in \mathbb R^{2d}$ and $w \in \mathbb R^{2d}$,
    \begin{equation}
        \norm{J_t} \Le C , \qquad \norm{K_t[w]} \Le C t |w|, \qquad \forall t\in[0,1].
        \label{eq: variation bounds}
    \end{equation}
    \emph{(Proof in Appendix~\ref{app: variations}.)}
\end{lemma}

\begin{lemma}
\label{lem: derivatives}
    Under Assumption~\ref{ass: U}, for every $f \in C^1(\mathbb R^{2d})$ with $\norm{\nabla f}_\infty <\infty$ and $t \Ge 0$, $e^{t\mathcal L}f \in C^1(\mathbb R^{2d})$, and
    \begin{equation}
        \norm{\nabla (e^{t\mathcal L} f)}_\infty \Le C \norm{\nabla f}_\infty , \qquad \forall t\in[0,1].
        \label{eq: gradient bound}
    \end{equation}
    If moreover $f \in C^2(\mathbb R^{2d})$ with $\norm{\nabla^2 f}_\infty < \infty$, then $e^{t\mathcal L}f \in C^2(\mathbb R^{2d})$ for every $t \in [0,1]$, and
    \begin{equation}
        \norm{\nabla^2 (e^{t\mathcal L} f)}_\infty \Le C \bigl( \norm{\nabla^2 f}_\infty + t \norm{\nabla f}_\infty \bigr) , \qquad \forall t\in[0,1].
        \label{eq: second derivative pass}
    \end{equation}
    \emph{(Proof in Appendix~\ref{app: variations}.)}
\end{lemma}

\subsection{BEL formula for degenerate noise}
\label{subsec: bel}

The BEL formula \citep{elworthy1994formulae} expresses the derivative of the semigroup without differentiating the test function. It integrates the inverse of the diffusion coefficient applied to $J_s w$ against $\D B_s$, so it requires the diffusion to be uniformly elliptic, which \eqref{eq: uld} is not. The linear operator $\mathcal A_Z$ of Definition~\ref{def: A} maps a direction $w$ to the rate $\mathcal A_Z w$ of a shift of the Brownian path. Under this shift, the solution started from a state displaced along $w$ rejoins the original solution at the terminal time, to first order in the displacement. The Girsanov density of this shift gives the BEL formula of Lemma~\ref{lem: bismut}, which \citet{guillin2012degenerate} prove for bounded test functions.

\begin{definition}
\label{def: A}
    Given $0 < t \Le 1$ and $z \in \mathbb R^{2d}$, let $Z_s = (X_s,V_s)$ be the solution to \eqref{eq: uld} with the initial state $z=(x,v)$ for $s\in [0,t]$. The linear operator $\mathcal A_Z$ acts on $w = (w^x,w^v) \in \mathbb R^{2d}$ by
    \begin{equation}
        (\mathcal A_Z w) (s) := \frac{1}{\sqrt{2\gamma}} \Bigl( \nabla^2 U(X_s)\, w^\circleddash(s) + \gamma\, \dot w^\circleddash(s) + \ddot w^\circleddash(s) \Bigr) \in \mathbb R^d, \qquad s \in [0,t] ,
        \label{eq: eta}
    \end{equation}
    where $w^\circleddash(s)$ is the cubic interpolant and each dot above $w^\circleddash(s)$ is a derivative in $s$,
    \begin{equation}
        w^\circleddash(s) := \Bigl( 1 - \frac st \Bigr)^2 \Bigl[ \Bigl( 1 + \frac{2s}t \Bigr) w^x + s\, w^v \Bigr] ,
        \quad w^\circleddash(0) = w^x , \quad \dot w^\circleddash(0) = w^v , \quad w^\circleddash(t) = \dot w^\circleddash(t) = 0 .
        \label{eq: hermite}
    \end{equation}
\end{definition}

For $s \in [0,t]$, the cubic interpolant \eqref{eq: hermite} satisfies
\begin{equation}
    |w^\circleddash(s)| \Le |w^x| + t |w^v| , \qquad
    |\dot w^\circleddash(s)| \Le \frac{2}{t} |w^x| + |w^v| , \qquad
    |\ddot w^\circleddash(s)| \Le \frac{6}{t^2} |w^x| + \frac{4}{t} |w^v| .
    \label{eq: hermite bounds}
\end{equation}
By \eqref{eq: U bounds} and \eqref{eq: hermite bounds}, the process $\mathcal A_Z w$ is pathwise bounded, and it is adapted since $(X_s)_{s \in [0,t]}$ is adapted.

\begin{lemma}[BEL formula]
\label{lem: bismut}
    Under Assumption~\ref{ass: U}, let $0 < t \Le 1$, $z \in \mathbb R^{2d}$ and $w = (w^x,w^v) \in \mathbb R^{2d}$. Then
    \begin{equation}
        \biggl( \int_0^t |(\mathcal A_Z w)(s)|^2 \,\D s \biggr)^{\frac12} \Le C \Bigl( t^{-\frac32} |w^x| + t^{-\frac12} |w^v| \Bigr) .
        \label{eq: eta size}
    \end{equation}
    For every $f \in C^1(\mathbb R^{2d})$ with $\norm{\nabla f}_\infty <\infty$,
    \begin{equation}
        \nabla (e^{t\mathcal L} f)(z)[w] = -\, \mathbb E^z \biggl[ f(Z_t) \int_0^t (\mathcal A_Z w)(s) \cdot \D B_s \biggr] .
        \label{eq: bismut}
    \end{equation}
    \emph{(Proof in Appendix~\ref{app: bismut}.)}
\end{lemma}

The stochastic integral of \eqref{eq: bismut} has mean zero, so a constant may be subtracted from the test function,
\begin{equation}
    \nabla (e^{t\mathcal L} f)(z)[w] = -\, \mathbb E^z \biggl[ \bigl( f(Z_t) - c \bigr) \int_0^t (\mathcal A_Z w)(s) \cdot \D B_s \biggr] , \qquad c \in \mathbb R.
    \label{eq: bismut centred}
\end{equation}
Also, the right side of \eqref{eq: bismut} tends to $\nabla f(z)[w]$ as $t \to 0$, as a formal argument shows. By the martingale representation theorem, there is an adapted process $H_s \in \mathbb R^d$, $s \in [0,t]$, such that
\begin{equation*}
    f(Z_t) = (e^{t\mathcal L}f)(z) + \int_0^t H_s \cdot \D B_s .
\end{equation*}
The noise enters \eqref{eq: uld} through the velocity with the factor $\sqrt{2\gamma}$, so $H_s$ is $\sqrt{2\gamma}$ times the velocity block of the derivative of $f(Z_t)$ in $Z_s$. For small $t$, $Z_t \approx z$, and \eqref{eq: first variation} over $[s,t]$ gives
\begin{equation*}
    H_s \approx \sqrt{2\gamma}\, \Bigl( (t-s)\, \nabla_x f(z) + \bigl( 1 - \gamma(t-s) \bigr) \nabla_v f(z) \Bigr) .
\end{equation*}
By the It\^o isometry, the right side of \eqref{eq: bismut} equals $-\mathbb E^z \int_0^t H_s \cdot (\mathcal A_Z w)(s) \,\D s$. Integration by parts in $s$, with the boundary values of $w^\circleddash$ in \eqref{eq: hermite}, shows that this expression tends to $\nabla_x f(z) \cdot w^x + \nabla_v f(z) \cdot w^v = \nabla f(z)[w]$.

\subsection{Uniform-in-time Hessian estimates}
\label{subsec: estimates}

The proof of Theorem~\ref{thm: hessian} needs a bound with no factor of the dimension on the variance of $\psi(Z_t)$ under $\mathbb E^z$, for $\psi$ with $\norm{\nabla\psi}_\infty <\infty$. The bound $|\psi(Z_t) - \psi(z)| \Le \norm{\nabla\psi}_\infty |Z_t - z|$ would bring in $\mathbb E^z|Z_t - z|^2$, which carries the factor $d$ of $\mathbb E|B_t|^2 = dt$.

\begin{lemma}
\label{lem: fluctuation}
    Under Assumption~\ref{ass: U}, for every $\psi \in C^1(\mathbb R^{2d})$ with $\norm{\nabla\psi}_\infty <\infty$ and $z \in \mathbb R^{2d}$,
    \begin{equation}
        \mathbb E^z \bigl| \psi(Z_t) - \mathbb E^z \psi(Z_t) \bigr|^2 \Le C\, t\, \norm{\nabla\psi}_\infty^2 , \qquad \forall t\in[0,1].
        \label{eq: fluctuation}
    \end{equation}
    \emph{(Proof in Appendix~\ref{app: fluctuation}.)}
\end{lemma}

\begin{theorem}[Hessian]
\label{thm: hessian}
    Under Assumption~\ref{ass: U}, for every $f \in C^1(\mathbb R^{2d})$ with $\norm{\nabla f}_\infty <\infty$, $e^{t\mathcal L}f \in C^2(\mathbb R^{2d})$ for every $t \in (0,1]$, and for every $z \in \mathbb R^{2d}$ and $w = (w^x,w^v) \in \mathbb R^{2d}$,
    \begin{equation}
        \bigl| \nabla^2 (e^{t\mathcal L} f)(z)[w] \bigr|
        \Le C \norm{\nabla f}_\infty \Bigl( t^{-\frac32} |w^x| + t^{-\frac12} |w^v| \Bigr) ,\qquad \forall t\in(0,1].
        \label{eq: hessian anisotropic}
    \end{equation}
    \emph{(Proof in Appendix~\ref{app: hessian}.)}
\end{theorem}

For $f\in C^1(\mathbb R^{2d})$ with  $\norm{\nabla f}_\infty <\infty$,
integrating \eqref{eq: hessian anisotropic} from $z' = (x',v')$ to $z = (x,v)$ gives
\begin{equation}
    \bigl| \nabla (e^{t\mathcal L} f)(z) - \nabla (e^{t\mathcal L} f)(z') \bigr|
    \Le C \norm{\nabla f}_\infty \Bigl( t^{-\frac32} |x - x'| + t^{-\frac12} |v - v'| \Bigr) , \qquad \forall t \in (0,1] ,
    \label{eq: hessian lipschitz}
\end{equation}
and thus \eqref{eq: hessian lipschitz} also holds for every Lipschitz function $f$ by mollification, with the Lipschitz constant of $f$ in place of $\norm{\nabla f}_\infty$.

For $t > 1$, Corollary~\ref{cor: uniform} combines Theorem~\ref{thm: hessian} with the exponential decay of $\nabla (e^{t\mathcal L} f)$ in Lemma~\ref{lem: contraction}.

\begin{lemma}
\label{lem: contraction}
    Under Assumptions~\ref{ass: U}--\ref{ass: damping}, there is $\lambda > 0$, determined by $m$, $M$, $R$ and $\gamma$, such that for every $f \in C^1(\mathbb R^{2d})$ with $\norm{\nabla f}_\infty <\infty$,
    \begin{equation}
        \norm{\nabla (e^{t\mathcal L} f)}_\infty \Le C\, e^{-\lambda t}\, \norm{\nabla f}_\infty , \qquad \forall t \Ge 0 .
        \label{eq: contraction}
    \end{equation}
\end{lemma}

Lemma~\ref{lem: contraction} follows from the $\mathcal W_1$ contraction of \citep[Theorem~5]{schuh2024global} at Dirac initial laws and Kantorovich--Rubinstein duality. Appendix~\ref{app: contraction} verifies the hypotheses of that theorem and shows that the rate $\lambda$ and the prefactor $C$ of \eqref{eq: contraction} depend on the radius $R$ and not on the dimension.

\begin{corollary}[Hessian]
\label{cor: uniform}
    Under Assumptions~\ref{ass: U}--\ref{ass: damping}, with $\lambda$ as in Lemma~\ref{lem: contraction}, for every $f \in C^1(\mathbb R^{2d})$ with $\norm{\nabla f}_\infty <\infty$ and $t > 0$, $e^{t\mathcal L}f \in C^2(\mathbb R^{2d})$ and
    \begin{equation}
        \norm{\nabla^2_{vv} (e^{t\mathcal L} f)}_\infty + \norm{\nabla^2_{xv} (e^{t\mathcal L} f)}_\infty \Le C \min(1,t)^{-\frac12} e^{-\lambda t} \norm{\nabla f}_\infty ,
        \label{eq: uniform}
    \end{equation}
    and
    \begin{equation}
        \norm{\nabla^2_{xx} (e^{t\mathcal L} f)}_\infty \Le C \min(1,t)^{-\frac32} e^{-\lambda t} \norm{\nabla f}_\infty .
        \label{eq: uniform xx}
    \end{equation}
    \emph{(Proof in Appendix~\ref{app: uniform}.)}
\end{corollary}

The exponents of $t$ in Theorem~\ref{thm: hessian} and Corollary~\ref{cor: uniform} cannot be improved, already for a quadratic potential in one dimension.

\begin{lemma}
\label{lem: sharp}
    Let $d = 1$, $U(x) = \frac12 x^2$ and $\gamma > 0$, which satisfies Assumptions~\ref{ass: U} and~\ref{ass: convexity} with $M = m = 1$ and $R = 0$. For the test functions $f_x(x,v) := |x|$ and $f_v(x,v) := |v|$, both with Lipschitz constant $1$, the functions $e^{t\mathcal L}f_x$ and $e^{t\mathcal L}f_v$ are smooth at every $t > 0$ and, as $t \to 0$,
    \begin{equation}
        \nabla^2_{xx} (e^{t\mathcal L} f_x)(0) \sim \sqrt{\frac{3}{\pi\gamma}}\, t^{-\frac32} , \quad
        \nabla^2_{xv} (e^{t\mathcal L} f_x)(0) \sim \sqrt{\frac{3}{\pi\gamma}}\, t^{-\frac12} , \quad
        \nabla^2_{vv} (e^{t\mathcal L} f_v)(0) \sim \sqrt{\frac{1}{\pi\gamma}}\, t^{-\frac12} .
        \label{eq: sharp}
    \end{equation}
    \emph{(Proof in Appendix~\ref{app: sharp}.)}
\end{lemma}

In \eqref{eq: hessian lipschitz}, the coefficient of $|x - x'|$  bounds $|\nabla^2_{xx} (e^{t\mathcal L} f_x)(0)|$, and the coefficient of $|v - v'|$ bounds $|\nabla^2_{xv} (e^{t\mathcal L} f_x)(0)|$ and $|\nabla^2_{vv} (e^{t\mathcal L} f_v)(0)|$. Hence no exponent of $t$ in \eqref{eq: hessian anisotropic} can be improved. If $f\in C^2(\mathbb R^{2d})$ with $\norm{\nabla f}_\infty + \norm{\nabla^2 f}_\infty < \infty$, then \eqref{eq: second derivative pass} implies that $\nabla^2(e^{t\mathcal L}f)$ is bounded as $t\to 0$, so choosing a nonsmooth $f$ in Lemma~\ref{lem: sharp} is necessary.

\begin{figure}[!htbp]
	\centering
	\includegraphics[width=0.75\textwidth]{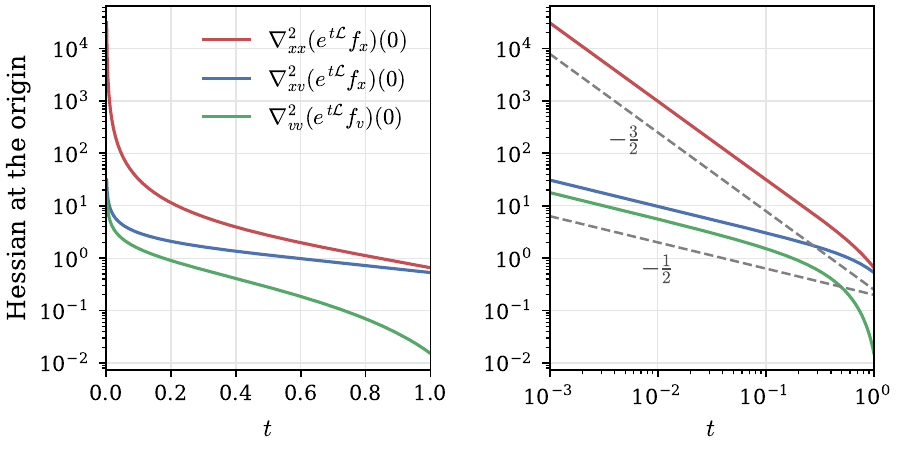}
	\caption{The Hessians $\nabla^2_{xx} (e^{t\mathcal L} f_x)(0)$, $\nabla^2_{xv} (e^{t\mathcal L} f_x)(0)$ and $\nabla^2_{vv} (e^{t\mathcal L} f_v)(0)$ of Lemma~\ref{lem: sharp} against $t$ for $\gamma = 1$, with $t$ on a linear axis (left) and in log-log scale (right). Dashed lines carry the slopes marked.}
	\label{fig: sharp}
\end{figure}

Figure~\ref{fig: sharp} plots the Hessians in \eqref{eq: sharp} against $t$ for $\gamma = 1$, evaluated from the Gaussian distribution of $Z_t$. In log-log scale, the curves have the slopes $-\frac32$, $-\frac12$ and $-\frac12$ of \eqref{eq: sharp} for small $t$.

\subsection{Estimate of third derivatives}
\label{subsec: third}

This subsection bounds the third derivative of $e^{t\mathcal L}f$ by $\norm{\nabla f}_\infty$ alone. For simplicity, we bound the derivative of the velocity Hessian $\nabla^2_{vv}(e^{t\mathcal L}f)$ only. For $\psi \in C^2(\mathbb R^{2d})$ with $\nabla^2_{vv}\psi$ continuously differentiable on $\mathbb R^{2d}$, the derivatives of $\nabla^2_{vv}\psi$ in $x$ and $v$ are the tensors $\nabla^3_{xvv}\psi$ and $\nabla^3_{vvv}\psi$ in $\mathbb R^{d\times d\times d}$ with the entries
\begin{equation*}
    \bigl( \nabla^3_{xvv} \psi \bigr)_{ijk} = \partial_{x_i}\partial_{v_j}\partial_{v_k} \psi , \qquad
    \bigl( \nabla^3_{vvv} \psi \bigr)_{ijk} = \partial_{v_i}\partial_{v_j}\partial_{v_k} \psi ,
\end{equation*}
with the sup norms $\norm{\nabla^3_{xvv} \psi}_\infty := \sup_{z\in\mathbb R^{2d}} \norm{\nabla^3_{xvv} \psi(z)}$ and $\norm{\nabla^3_{vvv} \psi}_\infty := \sup_{z\in\mathbb R^{2d}} \norm{\nabla^3_{vvv} \psi(z)}$.

\begin{theorem}[third derivative]
\label{thm: third}
    Under Assumption~\ref{ass: U}, for every $f \in C^1(\mathbb R^{2d})$ with $\norm{\nabla f}_\infty <\infty$ and $t \in (0,1]$, $\nabla^2_{vv}(e^{t\mathcal L}f)$ is continuously differentiable on $\mathbb R^{2d}$ and
    \begin{equation}
        \norm{\nabla^3_{xvv} (e^{t\mathcal L} f)}_\infty + \norm{\nabla^3_{vvv} (e^{t\mathcal L} f)}_\infty \Le C\, t^{-1} \norm{\nabla f}_\infty .
        \label{eq: third}
    \end{equation}
    \emph{(Proof in Appendix~\ref{app: third}.)}
\end{theorem}

The bound \eqref{eq: third} needs no fourth derivative of the potential. By \eqref{eq: second variation}, the second variation $K_t[w]$ is an integral over $s \in [0,t]$ whose integrand carries the factor $\nabla^3U(X_s)\bigl[ (J_sw)^x \bigr]$, and differentiating this factor in the initial state $z$ produces $\nabla^4U(X_s)$. The proof of Theorem~\ref{thm: third} instead moves the derivative in the initial state to the Brownian path on $[0,s]$, by the shift of the Brownian path in the proof of Lemma~\ref{lem: bismut}, which leaves $X_s$ unchanged.

Corollary~\ref{cor: third uniform} extends Theorem~\ref{thm: third} to every $t > 0$ under Assumptions~\ref{ass: U}--\ref{ass: damping}. For $t \Ge 1$, write $e^{t\mathcal L}f = e^{\mathcal L}(e^{(t-1)\mathcal L}f)$. Theorem~\ref{thm: third} at the time $1$ applied to $e^{(t-1)\mathcal L}f \in C^1(\mathbb R^{2d})$ and \eqref{eq: contraction} at the time $t-1$ give
\begin{equation*}
    \norm{\nabla^3_{xvv} (e^{t\mathcal L} f)}_\infty + \norm{\nabla^3_{vvv} (e^{t\mathcal L} f)}_\infty \Le C \norm{\nabla (e^{(t-1)\mathcal L} f)}_\infty \Le C e^{-\lambda t} \norm{\nabla f}_\infty .
\end{equation*}
With Theorem~\ref{thm: third} for $t \Le 1$, this bound proves Corollary~\ref{cor: third uniform}.

\begin{corollary}[third derivative]
\label{cor: third uniform}
    Under Assumptions~\ref{ass: U}--\ref{ass: damping}, with $\lambda$ as in Lemma~\ref{lem: contraction}, for every $f \in C^1(\mathbb R^{2d})$ with $\norm{\nabla f}_\infty <\infty$ and $t > 0$, $\nabla^2_{vv}(e^{t\mathcal L}f)$ is continuously differentiable on $\mathbb R^{2d}$ and
    \begin{equation}
        \norm{\nabla^3_{xvv} (e^{t\mathcal L} f)}_\infty + \norm{\nabla^3_{vvv} (e^{t\mathcal L} f)}_\infty \Le C \min(1,t)^{-1} e^{-\lambda t} \norm{\nabla f}_\infty .
        \label{eq: third uniform}
    \end{equation}
\end{corollary}

Like the bound \eqref{eq: hessian lipschitz}, Theorem~\ref{thm: third} and Corollary~\ref{cor: third uniform} also apply to every Lipschitz function $f$, with the Lipschitz constant of $f$ in place of $\norm{\nabla f}_\infty$. For such $f$, $e^{\frac t2\mathcal L}f \in C^1(\mathbb R^{2d})$ by mollification, and the two results apply to $e^{\frac t2\mathcal L}f$ at the time $\frac t2$, since $e^{t\mathcal L}f = e^{\frac t2\mathcal L}(e^{\frac t2\mathcal L}f)$.

The rate $t^{-1}$ in \eqref{eq: third} is sharp in the sup norm, already for a quadratic potential in one dimension.

\begin{lemma}
\label{lem: third sharp}
    Let $d = 1$, $U(x) = \frac12 x^2$ and $\gamma > 0$, which satisfies Assumptions~\ref{ass: U} and~\ref{ass: convexity} with $M = m = 1$ and $R = 0$. For the test functions $f_x(x,v) := |x|$ and $f_v(x,v) := |v|$, as $t \to 0$,
    \begin{equation}
        \norm{\nabla^3_{xvv} (e^{t\mathcal L} f_x)}_\infty \sim \frac{3}{\sqrt{2\pi e}\,\gamma}\, t^{-1} , \qquad
        \norm{\nabla^3_{vvv} (e^{t\mathcal L} f_v)}_\infty \sim \frac{1}{\sqrt{2\pi e}\,\gamma}\, t^{-1} .
        \label{eq: third sharp}
    \end{equation}
    \emph{(Proof in Appendix~\ref{app: third sharp}.)}
\end{lemma}

\begin{figure}[!htbp]
	\centering
	\includegraphics[width=0.75\textwidth]{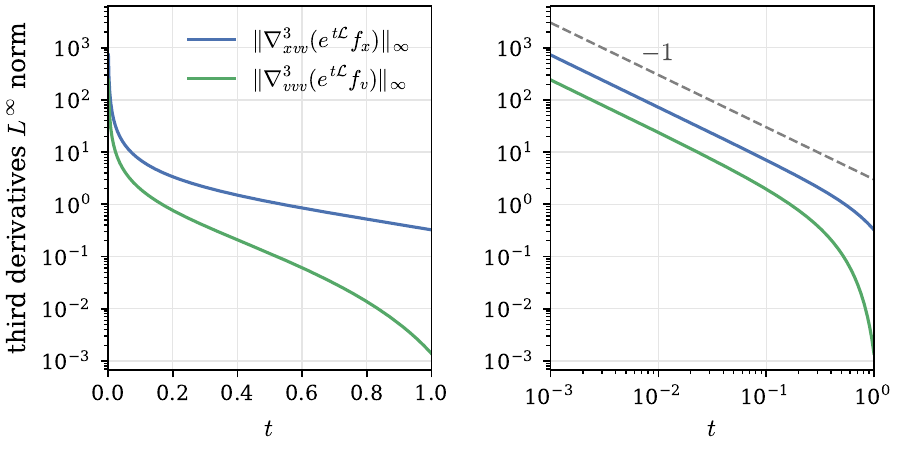}
	\caption{The sup norms $\norm{\nabla^3_{xvv} (e^{t\mathcal L} f_x)}_\infty$ and $\norm{\nabla^3_{vvv} (e^{t\mathcal L} f_v)}_\infty$ of Lemma~\ref{lem: third sharp} against $t$ for $\gamma = 1$, with $t$ on a linear axis (left) and in log-log scale (right). The dashed line carries the slope marked.}
	\label{fig: third sharp}
\end{figure}

Unlike the Hessians in \eqref{eq: sharp}, $\nabla^3_{xvv}(e^{t\mathcal L}f_x)$ and $\nabla^3_{vvv}(e^{t\mathcal L}f_v)$ vanish at the state $0$, so \eqref{eq: third sharp} is stated in the sup norm.
Figure~\ref{fig: third sharp} plots the sup norms in \eqref{eq: third sharp} against $t$ for $\gamma = 1$, evaluated from the Gaussian distribution of $Z_t$. In log-log scale, the curves have the slope $-1$ of \eqref{eq: third sharp} for small $t$.

\subsection{Poisson equation and Poincaré inequality}
\label{subsec: poisson}

Under Assumptions~\ref{ass: U}--\ref{ass: damping}, for $f \in C^1(\mathbb R^{2d})$ with $\norm{\nabla f}_\infty <\infty$, let $\phi$ be the function \eqref{eq: poisson solution},
\begin{equation*}
    \phi(z) = \int_0^\infty \bigl( (e^{t\mathcal L}f)(z) - \pi(f) \bigr) \D t .
\end{equation*}
The Hessian decay in Corollary~\ref{cor: uniform} ensures that $\phi$ is a strong solution to the Poisson equation $\mathcal L\phi = \pi(f) - f$.

\begin{lemma}[Poisson]
\label{lem: poisson}
    Under Assumptions~\ref{ass: U}--\ref{ass: damping}, for every $f \in C^1(\mathbb R^{2d})$ with $\norm{\nabla f}_\infty <\infty$, the integral in \eqref{eq: poisson solution} converges at every $z \in \mathbb R^{2d}$, and $\phi \in C^1(\mathbb R^{2d})$, $\nabla_v\phi \in C^1(\mathbb R^{2d})$ with
    \begin{equation}
        \norm{\nabla\phi}_\infty + \norm{\nabla^2_{vv}\phi}_\infty + \norm{\nabla^2_{xv}\phi}_\infty \Le C \norm{\nabla f}_\infty .
        \label{eq: poisson bounds}
    \end{equation}
    Moreover, $\phi$ solves the Poisson equation in the strong sense,
    \begin{equation}
        \mathcal L\phi = v \cdot \nabla_x \phi - \bigl( \nabla U(x) + \gamma v \bigr) \cdot \nabla_v \phi + \gamma \Delta_v \phi = \pi(f) - f .
        \label{eq: poisson equation}
    \end{equation}
    \emph{(Proof in Appendix~\ref{app: poisson}.)}
\end{lemma}

The function $\phi$ does not necessarily belong to $C^2(\mathbb R^{2d})$, because $\nabla^2_{xx}(e^{t\mathcal L}f)$ is not necessarily integrable at $t = 0$, as shown in Lemma~\ref{lem: sharp}. 
Since $\mathcal L$ has diffusion only in the velocity variable $v$,
It\^o's formula \citep[Chapter~IV, Remarks of Theorem~3.3]{revuz1999continuous} can be applied when $\phi \in C^1(\mathbb R^{2d})$ and $\nabla_v\phi \in C^1(\mathbb R^{2d})$, and it gives, for every solution $(Z_t)_{t \Ge 0}$ of \eqref{eq: uld} and $t \Ge 0$,
\begin{equation}
    \phi(Z_t) - \phi(Z_0) = \int_0^t (\mathcal L\phi)(Z_s) \,\D s + \sqrt{2\gamma} \int_0^t \nabla_v\phi(Z_s) \cdot \D B_s .
    \label{eq: ito formula}
\end{equation}
Let $(Z_t)_{t \Ge 0}$ solve \eqref{eq: uld} with $Z_0\sim\pi$ and independent of $(B_t)_{t \Ge 0}$. Then \eqref{eq: ito formula} and \eqref{eq: poisson equation} give, for $t \Ge 0$,
\begin{equation}
    \int_0^t \bigl( f(Z_s) - \pi(f) \bigr) \D s = \phi(Z_0) - \phi(Z_t) + \sqrt{2\gamma} \int_0^t \nabla_v\phi(Z_s) \cdot \D B_s ,
    \label{eq: poisson decomposition}
\end{equation}
which expresses the integral of $f(Z_s) - \pi(f)$ using the boundary term of $\phi$ and the stochastic integral of $\nabla_v\phi$, and corresponds to \citep[Equation~(4.5)]{mattingly2010convergence}.

Alternately, the Poisson solution $\phi$ can be characterized via the exponential decay of $e^{t\mathcal L}$ in $L^2(\pi)$, which does not require the damping condition \eqref{eq: friction condition}. This decay follows from a Poincaré inequality in the position. The next two lemmas state the Poincaré inequality and this decay, under Assumptions~\ref{ass: U} and~\ref{ass: convexity} alone.

\begin{lemma}
\label{lem: poincare}
    Under Assumptions~\ref{ass: U} and~\ref{ass: convexity}, there is a constant $C_P$, depending only on $m$, $M$ and $R$, such that every $\psi \in C^1(\mathbb R^{2d})$ that depends on the position $x$ only, with $\norm{\nabla \psi}_\infty <\infty$, satisfies
    \begin{equation}
        \pi\bigl( (\psi - \pi(\psi))^2 \bigr) \Le C_P\, \pi\bigl( |\nabla \psi|^2 \bigr) .
        \label{eq: poincare}
    \end{equation}
    \emph{(Proof in Appendix~\ref{app: poincare}.)}
\end{lemma}

Lemma~\ref{lem: poincare} requires $C_P$ to be independent of the dimension $d$, hence the Poincaré inequality in \citep{ye2026sglm} cannot be directly applied. Using \citep[Proposition~4.3.1]{bakry2014analysis}, $\pi \propto \exp(-U(x)-\frac12|v|^2)$ itself also satisfies the Poincaré inequality with the constant $\max\{C_P,1\}$: for every $\psi \in C^1(\mathbb R^{2d})$ with $\norm{\nabla \psi}_\infty <\infty$,
\begin{equation}
    \pi\bigl( (\psi - \pi(\psi))^2 \bigr) \Le \max\{C_P,1\}\, \pi\bigl( |\nabla \psi|^2 \bigr) .
    \label{eq: poincare full}
\end{equation}

\begin{lemma}
\label{lem: l2 decay}
    Under Assumptions~\ref{ass: U} and~\ref{ass: convexity}, let $\psi \in L^2(\pi)$ satisfy $\pi(\psi) = 0$. Then
    \begin{equation}
        \pi\bigl( (e^{t\mathcal L}\psi)^2 \bigr) \Le C e^{-t/C}\, \pi(\psi^2) , \qquad \forall t \Ge 0 ,
        \label{eq: l2 decay}
    \end{equation}
    where $(e^{t\mathcal L}\psi)(z) = \mathbb E^z\bigl[ \psi(Z_t) \bigr]$ is defined for $\pi$-almost every $z$.
    \emph{(Proof in Appendix~\ref{app: l2}.)}
\end{lemma}

Lemma~\ref{lem: l2 decay} follows from \citep[Theorem~1]{cao2023explicit} and the Poincaré inequality \eqref{eq: poincare}. The rate and the prefactor of that theorem, like $C_P$, do not depend on the dimension.

Under Assumptions~\ref{ass: U} and~\ref{ass: convexity}, for $f \in C^1(\mathbb R^{2d})$ with $\norm{\nabla f}_\infty <\infty$, the Poincaré inequality \eqref{eq: poincare full} gives $\pi\bigl( (f - \pi(f))^2 \bigr) \Le C \norm{\nabla f}_\infty^2$, and then the exponential decay \eqref{eq: l2 decay} with $\psi = f - \pi(f)$ gives
\begin{equation*}
    \pi\bigl( ( e^{t\mathcal L}f - \pi(f) )^2 \bigr)^{\frac12} \Le C e^{-t/C} \norm{\nabla f}_\infty , \qquad \forall t \Ge 0 .
\end{equation*}
Hence the integral \eqref{eq: poisson solution} converges in $L^2(\pi)$, and its limit $\phi$ satisfies $\pi(\phi^2)^{\frac12} \Le C \norm{\nabla f}_\infty$. By the semigroup property, $e^{s\mathcal L}\phi - \phi = -\int_0^s \bigl( e^{t\mathcal L}f - \pi(f) \bigr) \D t$ for every $s > 0$. Since $e^{t\mathcal L}f \to f$ in $L^2(\pi)$ as $t \to 0$, $\phi$ solves the Poisson equation in the sense that
\begin{equation*}
    \lim_{s \to 0} \frac{e^{s\mathcal L}\phi - \phi}{s} = \pi(f) - f \qquad \text{in } L^2(\pi) .
\end{equation*}

\section{Uniform-in-time error in stochastic gradient sampling}
\label{sec: sgubu}

In stochastic gradient sampling, an unbiased stochastic gradient $\hat g(x;\xi) \in \mathbb R^d$ with a random variable $\xi$ replaces the gradient $\nabla U(x)$ of \eqref{eq: uld}. For the SG-UBU integrator \citep{sanz2021wasserstein, paulin2024correction, schuh2024convergence, lu2025mse}, this section bounds, uniformly in time, the error between the numerical solution and the solution of \eqref{eq: uld}. Also, we quantify the numerical bias with an explicit coefficient, given by the solution to the Poisson equation \eqref{eq: poisson equation}.

\subsection{SG-UBU integrator}

The UBU integrator splits the underdamped Langevin dynamics \eqref{eq: uld} into
\begin{equation*}
\begin{aligned}
    & \mathscr U : \quad \D X_t = V_t \,\D t , && \D V_t = -\gamma V_t \,\D t + \sqrt{2\gamma}\,\D B_t , \\
    & \mathscr B : \quad \D X_t = 0 , && \D V_t = -\nabla U(X_t) \,\D t .
\end{aligned}
\end{equation*}
For a given Brownian motion $(B_t)_{t \Ge 0}$, let $\mathscr U_t$ and $\mathscr B_t$ be the solution maps of $\mathscr U$ and $\mathscr B$ on $\mathbb R^{2d}$ at the time $t \Ge 0$. The map $\mathscr U_t$ is explicitly expressed as
\begin{equation}
    \mathscr U_t : (x,v) \mapsto \biggl( x + \frac{1-e^{-\gamma t}}{\gamma}\, v + \sqrt{2\gamma} \int_0^t \frac{1-e^{-\gamma(t-s)}}{\gamma} \,\D B_s ,\; e^{-\gamma t} v + \sqrt{2\gamma} \int_0^t e^{-\gamma(t-s)} \,\D B_s \biggr) ,
    \label{eq: U flow}
\end{equation}
and $\mathscr B_t : (x,v) \mapsto (x, v - t\nabla U(x))$. Replacing the gradient $\nabla U(x)$ in $\mathscr B_t$ by the stochastic gradient $\hat g(x;\xi)$ gives the unbiased approximation
\begin{equation}
    \mathscr B_t(\xi) : (x,v) \mapsto \bigl( x ,\, v - t\, \hat g(x;\xi) \bigr) .
    \label{eq: B flow}
\end{equation}

One step of SG-UBU with step size $h > 0$ is $\mathscr U_{h/2} \circ \mathscr B_h(\xi) \circ \mathscr U_{h/2}$. The numerical solution $\mathbf Z_n = (\mathbf X_n, \mathbf V_n)$ is obtained by the iterative scheme
\begin{equation}
    \mathbf Z_{n+1} = \bigl( \mathscr U_{h/2} \circ \mathscr B_h(\xi_n) \circ \mathscr U_{h/2} \bigr) \mathbf Z_n , \qquad
    n = 0,1,2,\dots ,
    \tag{SG-UBU}
    \label{eq: sgubu}
\end{equation}
where $(\xi_n)_{n \Ge 0}$ are independent and have the law of $\xi$. In the step from $\mathbf Z_n$, the map $\mathscr U_{h/2}$ applied first is driven by the increments of $(B_t)_{t \Ge 0}$ on $[nh, (n+\frac12)h]$, and the map $\mathscr U_{h/2}$ applied last by the increments on $[(n+\frac12)h, (n+1)h]$. The initial state $\mathbf Z_0$, the sequence $(\xi_n)_{n \Ge 0}$ and $(B_t)_{t \Ge 0}$ are independent.

\subsection{Error analysis results}

For simplicity, the initial state is chosen as $\mathbf X_0 = 0$ and $\mathbf V_0 \sim \mathcal N(0,I_d)$. The solution to the underdamped Langevin dynamics \eqref{eq: uld} is denoted by $Z_t = (X_t,V_t)$, and its initial state is $X_0 = 0$ and $V_0 = \mathbf V_0$.

Denote the stochastic gradient noise by
\begin{equation}
    \hat\eta(x;\xi) := \hat g(x;\xi) - \nabla U(x) \in \mathbb R^d .
    \label{eq: sg noise}
\end{equation}
We impose the following assumption on the stochastic gradient noise.

\begin{assumption}[noise]
\label{ass: sg}
    For some $\Sigma\Ge0$, the stochastic gradient noise \eqref{eq: sg noise} satisfies $\mathbb E_\xi\bigl[ \hat\eta(x;\xi) \bigr] = 0$, the moment bound
    \begin{equation}
        \mathbb E_\xi\bigl[ |\hat\eta(x;\xi)|^4 \bigr] \Le \Sigma^2 d^2 , \qquad \forall x \in \mathbb R^d ,
        \label{eq: sg variance}
    \end{equation}
    and the Lipschitz bound
    \begin{equation}
        \mathbb E_\xi\bigl[ |\hat\eta(x;\xi) - \hat\eta(x';\xi)|^2 \bigr] \Le \Sigma |x - x'|^2 , \qquad \forall x, x' \in \mathbb R^d ,
        \label{eq: sg lipschitz}
    \end{equation}
    where $\mathbb E_{\xi}$ denotes the expectation with respect to $\xi$.
\end{assumption}

The Lipschitz bound \eqref{eq: sg lipschitz} is used only in the proofs of Theorems~\ref{thm: bias} and~\ref{thm: coefficient limit}. We first state the uniform-in-time moment bound of SG-UBU, a similar form of which is proved in \citep[Theorem~3]{lu2025mse} for a globally strongly convex potential.

\begin{lemma}
\label{lem: sgubu moments}
    Under Assumptions~\ref{ass: U}--\ref{ass: sg}, let the step size satisfy
    \begin{equation}
    	h \Le h_0 := \min\biggl\{ \frac1{2\gamma} \Bigl( 1 - \frac{4M^2}{\gamma^2 m} \Bigr) ,\; \frac{\gamma m}{32M^2 + 576\,\Sigma} \biggr\} .
        \label{eq: step size}
    \end{equation}
    Then the SG-UBU solution satisfies
    \begin{equation}
        \sup_{n \Ge 0} \mathbb E |\mathbf Z_n|^4 \Le C d^2 .
        \label{eq: sgubu moments}
    \end{equation}
    \emph{(Proof in Appendix~\ref{app: sgubu moments}.)}
\end{lemma}

To state the error analysis results, we introduce a few new notations. For $n \Ge 0$, define
\begin{equation}
\begin{aligned}
    & \mathbf Y_n := \mathbf X_n + \frac{1-e^{-\frac{\gamma h}2}}{\gamma}\, \mathbf V_n + \sqrt{2\gamma} \int_0^{\frac h2} \frac{1-e^{-\gamma(\frac h2-s)}}{\gamma} \,\D B_{nh+s} \in \mathbb R^d , \\
    & \bar{\mathbf Z}_{n+1} := \bigl( \mathscr U_{h/2} \circ \mathscr B_h \circ \mathscr U_{h/2} \bigr) \mathbf Z_n \in \mathbb R^{2d} ,
\end{aligned}
    \label{eq: full gradient step}
\end{equation}
where the maps $\mathscr U_{h/2}$ applied first and last are driven by the increments of $(B_t)_{t \Ge 0}$ on $[nh, (n+\frac12)h]$ and $[(n+\frac12)h, (n+1)h]$. By \eqref{eq: U flow}, $\mathbf Y_n$ is the position of $\mathscr U_{h/2} \mathbf Z_n$, at which \eqref{eq: sgubu} evaluates the stochastic gradient. The state $\bar{\mathbf Z}_{n+1}$ is the result of the step from $\mathbf Z_n$ with $\nabla U(\mathbf Y_n)$ in place of $\hat g(\mathbf Y_n;\xi_n)$.

By \eqref{eq: B flow}, the states $\mathscr B_h(\xi_n) \mathscr U_{h/2} \mathbf Z_n$ and $\mathscr B_h \mathscr U_{h/2} \mathbf Z_n$ differ only in the velocity, by $-h\, \hat\eta(\mathbf Y_n;\xi_n)$. Since the map $\mathscr U_{h/2}$ in \eqref{eq: U flow} is affine in $(x,v)$ for a given Brownian motion, $\mathbf Z_{n+1} - \bar{\mathbf Z}_{n+1}$ is linear in $\hat\eta(\mathbf Y_n;\xi_n)$,
\begin{equation}
    \mathbf Z_{n+1} - \bar{\mathbf Z}_{n+1} = -h \biggl( \frac{1-e^{-\frac{\gamma h}2}}{\gamma}\, \hat\eta(\mathbf Y_n;\xi_n) ,\; e^{-\frac{\gamma h}2}\, \hat\eta(\mathbf Y_n;\xi_n) \biggr) .
    \label{eq: sg error}
\end{equation}
Since $\xi_n$ is independent of $\mathbf Z_n$ and of $(B_t)_{t \Ge 0}$, the difference \eqref{eq: sg error} has mean zero given $\mathbf Z_n$ and $(B_t)_{t \Ge 0}$.

\begin{theorem}[uniform-in-time error]
\label{thm: sgubu}
    Under Assumptions~\ref{ass: U}--\ref{ass: sg}, let $h \Le \min\{ h_0, \frac1d \}$ with $h_0$ in \eqref{eq: step size}. Then for every $f \in C^1(\mathbb R^{2d})$ with $\norm{\nabla f}_\infty < \infty$,
    \begin{equation}
        \sup_{n\Ge0} \,\Bigl| \mathbb E [f(\mathbf Z_n)] - \mathbb E [f(Z_{nh})] \Bigr| \Le C \norm{\nabla f}_\infty \bigl( \Sigma d h + d h^2 \bigr) .
        \label{eq: sgubu error}
    \end{equation}
    \emph{(Proof in Appendix~\ref{app: sgubu error}.)}
\end{theorem}

By Kantorovich--Rubinstein duality, \eqref{eq: sgubu error} gives, in the Wasserstein-1 distance on $\mathbb R^{2d}$,
\begin{equation}
    \sup_{n \Ge 0} \mathcal W_1\bigl( \mathrm{Law}(\mathbf Z_n), \mathrm{Law}(Z_{nh}) \bigr) \Le C \bigl( \Sigma d h + d h^2 \bigr) .
    \label{eq: sgubu w1}
\end{equation}
The term $\Sigma d h$ of \eqref{eq: sgubu w1} agrees with \citep[Theorem~2.15 and Corollary~2.16]{paulin2026theoretical}.

\subsection{Coefficient of the numerical bias}
\label{subsec: bias}

A probability distribution $\pi_h$ on $\mathbb R^{2d}$ is an invariant distribution of SG-UBU if the step \eqref{eq: sgubu} from an initial state $\mathbf Z_0 \sim \pi_h$ gives $\mathbf Z_1 \sim \pi_h$. We first state the existence and the uniqueness of the invariant distribution, with a fourth moment bounded uniformly in $h$.

\begin{lemma}
\label{lem: invariant}
    Under Assumptions~\ref{ass: U}--\ref{ass: sg}, let $h \Le h_0$ with $h_0$ in \eqref{eq: step size}. Then  \eqref{eq: sgubu} has a unique invariant distribution $\pi_h$, which satisfies
    \begin{equation}
        \int_{\mathbb R^{2d}} |z|^4 \,\pi_h(\D z) \Le C d^2 .
        \label{eq: invariant moments}
    \end{equation}
    \emph{(Proof in Appendix~\ref{app: invariant}.)}
\end{lemma}

For a test function $f$, the numerical bias of SG-UBU is $\pi_h(f) - \pi(f)$, where $\pi_h(f) := \int_{\mathbb R^{2d}} f \,\D\pi_h$. We show that the numerical bias has an explicit coefficient in the step size $h$.

\begin{theorem}[numerical bias]
\label{thm: bias}
    Under Assumptions~\ref{ass: U}--\ref{ass: sg}, let $h \Le \min\{ h_0, \frac1d \}$ with $h_0$ in \eqref{eq: step size}. Then for every $f \in C^1(\mathbb R^{2d})$ with $\norm{\nabla f}_\infty < \infty$, the invariant distribution $\pi_h$ satisfies
    \begin{equation}
        \Bigl| \pi_h(f) - \pi(f) - h\, \mathcal A_h(f) \Bigr| \Le C \norm{\nabla f}_\infty (1+\Sigma)^2\, d^2 h^2 ,
        \label{eq: bias}
    \end{equation}
    where $\mathcal A_h(f)$ is the coefficient functional,
    \begin{equation}
        \mathcal A_h(f) :=
        h \sum_{j=1}^\infty \mathbb E^{\pi_h}_\xi \biggl[ \int_0^1 (1-\theta)\, \nabla^2_{vv}\bigl( e^{jh\mathcal L} f \bigr)\bigl( \bar{\mathbf Z}_1 + \theta ( \mathbf Z_1 - \bar{\mathbf Z}_1 ) \bigr)\bigl[ \hat\eta(\mathbf Y_0;\xi), \hat\eta(\mathbf Y_0;\xi) \bigr] \D\theta \biggr] ,
        \label{eq: bias coefficient}
    \end{equation}
    and $\mathbb E^{\pi_h}_\xi$ denotes the expectation over the initial state $\mathbf Z_0 \sim \pi_h$, the random variable $\xi$ of the step from $\mathbf Z_0$ and the Brownian motion $(B_t)_{t \Ge 0}$. Moreover,
    \begin{equation}
        |\mathcal A_h(f)| \Le C \norm{\nabla f}_\infty \Sigma d .
        \label{eq: bias coefficient bound}
    \end{equation}
    \emph{(Proof in Appendix~\ref{app: bias}.)}
\end{theorem}

Denote by $\mathcal C(x) := \mathbb E_{\xi}\bigl[ \hat\eta(x;\xi) \otimes \hat\eta(x;\xi) \bigr] \in \mathbb R^{d\times d}$ the covariance matrix of the stochastic gradient noise. As $h \to 0$, $\mathbf Z_1$ and $\bar{\mathbf Z}_1$ become the initial state $\mathbf Z_0 = (\mathbf X_0, \mathbf V_0) \sim \pi$ and $\mathbf Y_0$ becomes the position $\mathbf X_0$, hence the functional $\mathcal A_h(f)$ formally converges,
\begin{equation}
	\lim_{h\to0} \mathcal A_h(f) = \frac12\int_0^\infty \mathbb E_\xi^{\pi}
	\Bigl[
	\nabla^2_{vv}(e^{t\mathcal L}f)(\mathbf Z_0) \bigl[ \hat\eta(\mathbf X_0;\xi), \hat\eta(\mathbf X_0;\xi) \bigr] 
	\Bigr] \D t = \frac12\, \pi\bigl( \nabla^2_{vv}\phi : \mathcal C \bigr) ,
	\label{eq: formal limit}
\end{equation}
where $\mathbb E_\xi^{\pi}$ denotes the expectation over $\mathbf Z_0 \sim \pi$ and $\xi$, and $\phi$ is the solution \eqref{eq: poisson solution} to the Poisson equation \eqref{eq: poisson equation}. Theorem~\ref{thm: coefficient limit} quantifies the convergence in \eqref{eq: formal limit} in the step size $h$.

\begin{theorem}[coefficient]
\label{thm: coefficient limit}
    Under Assumptions~\ref{ass: U}--\ref{ass: sg}, let $h \Le \min\{ h_0, \frac1d \}$ with $h_0$ in \eqref{eq: step size}. Then for every $f \in C^1(\mathbb R^{2d})$ with $\norm{\nabla f}_\infty < \infty$,
    \begin{equation}
        \Bigl| \mathcal A_h(f) - \frac12\, \pi\bigl( \nabla^2_{vv}\phi : \mathcal C \bigr) \Bigr| \Le C \norm{\nabla f}_\infty\, \Sigma (1+\Sigma)\, d^2 h \Bigl( 1 + \log\frac1h \Bigr) .
        \label{eq: coefficient limit}
    \end{equation}
    \emph{(Proof in Appendix~\ref{app: coefficient limit}.)}
\end{theorem}

Combining Theorems~\ref{thm: bias} and~\ref{thm: coefficient limit} gives the numerical bias of SG-UBU to first order in $h$, with an explicit coefficient, given by the velocity Hessian $\nabla^2_{vv}\phi$ and the covariance matrix $\mathcal C(x)$.

\begin{corollary}
\label{cor: bias limit}
    Under Assumptions~\ref{ass: U}--\ref{ass: sg}, let $h \Le \min\{ h_0, \frac1d \}$ with $h_0$ in \eqref{eq: step size}. Then for every $f \in C^1(\mathbb R^{2d})$ with $\norm{\nabla f}_\infty < \infty$,
    \begin{equation}
        \Bigl| \pi_h(f) - \pi(f) - \frac h2\, \pi\bigl( \nabla^2_{vv}\phi : \mathcal C \bigr) \Bigr| \Le C \norm{\nabla f}_\infty (1+\Sigma)^2\, d^2 h^2 \Bigl( 1 + \log\frac1h \Bigr) .
        \label{eq: bias expansion}
    \end{equation}
    As a consequence,
    \begin{equation}
        \lim_{h \to 0} \frac{\pi_h(f) - \pi(f)}{h} = \frac12\, \pi\bigl( \nabla^2_{vv}\phi : \mathcal C \bigr) .
        \label{eq: bias coefficient limit}
    \end{equation}
\end{corollary}

Corollary~\ref{cor: bias limit} sharpens \citep[Theorem~9]{lu2025mse} in three ways: the potential need not be globally convex, the limit \eqref{eq: bias coefficient limit} is proved rigorously, and the test function needs only a bounded gradient.

\subsection{Linear response in the numerical bias}
\label{subsec: response}

This subsection represents the coefficient of the numerical bias in Corollary~\ref{cor: bias limit} by a function on $\mathbb R^{2d}$. The representation uses the exponential decay \eqref{eq: l2 decay} of $e^{t\mathcal L}$ in $L^2(\pi)$.

The coefficient $\frac12\, \pi\bigl( \nabla^2_{vv}\phi : \mathcal C \bigr)$ is linear in the test function $f$, and by the Riesz representation theorem, a function $\mathfrak A$ on $\mathbb R^{2d}$ that does not depend on $f$ represents it,
\begin{equation}
    \frac12\, \pi\bigl( \nabla^2_{vv}\phi : \mathcal C \bigr) = \pi(\mathfrak A f) = \int_{\mathbb R^{2d}} \mathfrak A(x,v)\, f(x,v) \,\pi(\D x\, \D v) .
    \label{eq: response}
\end{equation}
To construct $\mathfrak A$, note that under $\pi$ the velocity has the law $\mathcal N(0,I_d)$ and is independent of the position, and that $( v \otimes v - I_d ) : \mathcal C(x)$ has mean zero under $\pi$. Hence integration by parts in $v$ and \eqref{eq: poisson solution} give
\begin{equation}
    \pi\bigl( \nabla^2_{vv}\phi : \mathcal C \bigr) = \pi\Bigl( \phi\, \bigl( v \otimes v - I_d \bigr) : \mathcal C(x) \Bigr) = \int_0^\infty \pi\Bigl( (e^{t\mathcal L}f)\, \bigl( v \otimes v - I_d \bigr) : \mathcal C(x) \Bigr) \D t .
    \label{eq: coefficient integral}
\end{equation}
The adjoint $\mathcal L^*$ of $\mathcal L$ in $L^2(\pi)$ is\footnote{The operator $\mathcal L^*$ is not the Fokker--Planck operator, which is the adjoint of $\mathcal L$ in $L^2(\mathbb R^{2d})$. For overdamped Langevin dynamics, the generator $\mathcal L$ is self-adjoint in $L^2(\pi)$.}
\begin{equation*}
    \mathcal L^* f = -v \cdot \nabla_x f + \bigl( \nabla U(x) - \gamma v \bigr) \cdot \nabla_v f + \gamma \Delta_v f ,
\end{equation*}
which is $\mathcal L$ with $v$ replaced by $-v$. Since $( v \otimes v - I_d ) : \mathcal C(x)$ is even in $v$, for every $t \Ge 0$,
\begin{equation*}
\begin{aligned}
    & \pi\Bigl( (e^{t\mathcal L}f)\, \bigl( v \otimes v - I_d \bigr) : \mathcal C(x) \Bigr) = \pi\Bigl( f\, e^{t\mathcal L^*} \bigl[ ( v \otimes v - I_d ) : \mathcal C(x) \bigr] \Bigr) \\
    & \qquad = \int_{\mathbb R^{2d}} f(x,v)\, \Bigl( e^{t\mathcal L} \bigl[ ( v \otimes v - I_d ) : \mathcal C(x) \bigr] \Bigr)(x,-v) \,\pi(\D x\, \D v) .
\end{aligned}
\end{equation*}
Hence \eqref{eq: coefficient integral} gives \eqref{eq: response} with the accurate expression of $\mathfrak A(x,v)$:
\begin{equation}
    \mathfrak A(x,v) := \frac12 \int_0^\infty \Bigl( e^{t\mathcal L} \bigl[ ( v \otimes v - I_d ) : \mathcal C(x) \bigr] \Bigr)(x,-v) \,\D t .
    \label{eq: response function}
\end{equation}
Since $\pi$ is even in $v$, \eqref{eq: sg variance} and \eqref{eq: l2 decay} give, for the integrand of \eqref{eq: response function} and every $t \Ge 0$,
\begin{equation*}
\begin{aligned}
    \int_{\mathbb R^{2d}} \biggl( \Bigl( e^{t\mathcal L} \bigl[ ( v \otimes v - I_d ) : \mathcal C(x) \bigr] \Bigr)(x,-v) \biggr)^2 \,\pi(\D x\, \D v) = \pi\biggl( \Bigl( e^{t\mathcal L} \bigl[ ( v \otimes v - I_d ) : \mathcal C(x) \bigr] \Bigr)^2 \biggr)
    \Le C e^{-t/C}\, \Sigma^2 d^2 .
\end{aligned}
\end{equation*}
Hence the integral \eqref{eq: response function} converges in $L^2(\pi)$, and $\mathfrak A \in L^2(\pi)$ with $\pi(\mathfrak A^2) \Le C \Sigma^2 d^2$.

If $\mathcal C(x) = \sigma^2 I_d$ is a constant matrix, $\mathfrak A(x,v)$ is explicit. Since $\mathcal L \bigl( U(x) + \frac12 |v|^2 \bigr) = \gamma \bigl( d - |v|^2 \bigr)$, \eqref{eq: response function} gives
\begin{equation}
    \mathfrak A(x,v) = \frac{\sigma^2}{2\gamma} \Bigl( U(x) + \frac12 |v|^2 - \pi(U) - \frac d2 \Bigr) .
    \label{eq: response constant}
\end{equation}
Hence $\mathfrak A(x,v)$ scales linearly with $\gamma^{-1}$.

With \eqref{eq: response}, the expansion \eqref{eq: bias expansion} compares $\pi_h(f)$ with $\pi\bigl( (1 + h\, \mathfrak A) f \bigr)$. In this sense, $\mathfrak A(x,v)$ is the linear response of $\pi$ under the SG-UBU discretization with the step size $h$.
The expansion \eqref{eq: bias expansion} and $\mathfrak A$ also give an explicit limit of the Wasserstein-1 distance between $\pi_h$ and $\pi$ on $\mathbb R^{2d}$.

\begin{corollary}
\label{lem: w1 limit}
Under Assumptions~\ref{ass: U}--\ref{ass: sg},
\begin{equation}
	\lim_{h \to 0} \frac{\mathcal W_1(\pi_h,\pi)}{h} = \norm{\mathfrak A \pi}_{\mathrm{KR}} ,
	\label{eq: w1 limit}
\end{equation}
where $\mathfrak A \pi$ is the signed measure $\mathfrak A(z) \pi(\D z)$ on $\mathbb R^{2d}$, and $\norm{\cdot}_{\mathrm{KR}}$ is the Kantorovich--Rubinstein norm of a signed measure $\nu$ on $\mathbb R^{2d}$ with $\nu(\mathbb R^{2d}) = 0$ and $\int_{\mathbb R^{2d}} |z|\, |\nu|(\D z) < \infty$,
\begin{equation*}
	\norm{\nu}_{\mathrm{KR}} := \sup \biggl\{ \int_{\mathbb R^{2d}} f \,\D\nu : f \text{ has Lipschitz constant at most } 1 \biggr\} .
\end{equation*}
\emph{(Proof in Appendix~\ref{app: w1 limit}.)}
\end{corollary}

The explicit Wasserstein-1 limit is also proved for the stochastic gradient Leimkuhler--Matthews (SGLM) integrator \citep{ye2026sglm}. As a comparison, \citep{ye2026sglm} requires bounded derivatives of $U$ up to $\nabla^7 U$, while \eqref{eq: w1 limit} only needs bounded derivatives up to $\nabla^3 U$.

\begin{figure}[!htbp]
	\centering
	\includegraphics[width=0.48\textwidth]{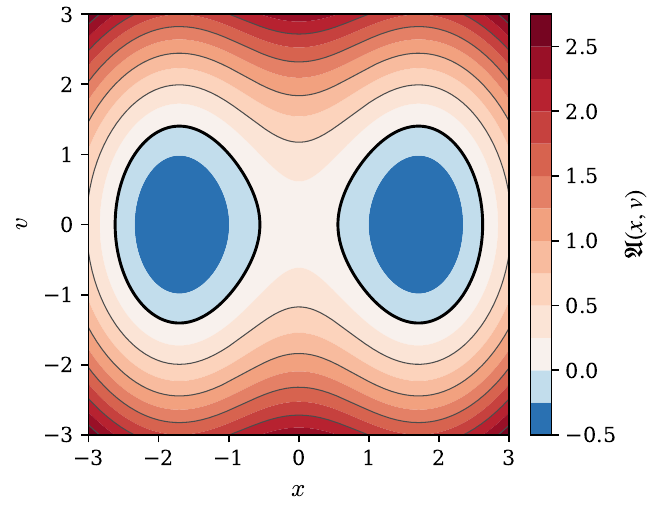}
	\includegraphics[width=0.43\textwidth]{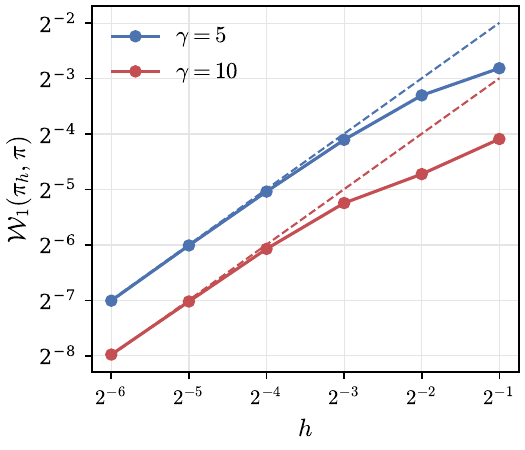}
	\caption{The function $\mathfrak A(x,v)$ for $\gamma = 10$, with the level $0$ drawn as the thick line (left), and $\mathcal W_1(\pi_h,\pi)$ against $h$ in log-log scale for $\gamma = 5$ and $\gamma = 10$ (right). The dashed lines are $\norm{\mathfrak A \pi}_{\mathrm{KR}}\, h$.}
	\label{fig: w1}
\end{figure}

We test \eqref{eq: w1 limit} in dimension $d = 1$ with the potential $U(x) = \frac12 x^2 + 4 e^{-0.4x^2}$ and the stochastic gradient $\hat g(x;\xi) = \nabla U(x) + \sqrt{10}\, \xi$ with $\xi \sim \mathcal N(0,1)$, so that $\mathcal C(x) = 10$. The damping rates are $\gamma = 5$ and $\gamma = 10$, where $\gamma = 10$ satisfies the damping condition \eqref{eq: friction condition} while $\gamma = 5$ does not. Figure~\ref{fig: w1} plots $\mathfrak A(x,v)$, and $\mathcal W_1(\pi_h,\pi)$ against $h$. The function $\mathfrak A(x,v)$ is plotted only for $\gamma = 10$ because of \eqref{eq: response constant}. The distance $\mathcal W_1(\pi_h,\pi)$ is estimated from samples of SG-UBU and compared with $\norm{\mathfrak A \pi}_{\mathrm{KR}}\, h$, where $\norm{\mathfrak A \pi}_{\mathrm{KR}} \approx 0.502$ for $\gamma = 5$ and $\norm{\mathfrak A \pi}_{\mathrm{KR}} \approx 0.251$ for $\gamma = 10$ are computed numerically. For both damping rates, $\mathcal W_1(\pi_h,\pi)/h$ approaches $\norm{\mathfrak A \pi}_{\mathrm{KR}}$ as $h$ decreases.

\section{Conclusion}
\label{sec: conclusion}

This paper bounds the Hessian of $e^{t\mathcal L}f$ by $\norm{\nabla f}_\infty$ at every $t > 0$, for potentials with bounded second and third derivatives that are strongly convex only outside a ball, and a damping rate above an explicit threshold, with constants free of the dimension. The derivatives $\nabla^3_{xvv}$ and $\nabla^3_{vvv}$ of the velocity Hessian of $e^{t\mathcal L}f$ are also bounded by $\norm{\nabla f}_\infty$, with the short-time factor $t^{-1}$ and no fourth derivative of the potential. The short-time exponents of both bounds cannot be improved. The bounds on the blocks $\nabla^2_{vv}$ and $\nabla^2_{xv}$ are integrable in time, so the function $\phi$ in \eqref{eq: poisson solution} is a strong solution to the Poisson equation. As an application, Theorem~\ref{thm: sgubu} proves the uniform-in-time error of SG-UBU, Corollary~\ref{cor: bias limit} gives the explicit coefficient of the numerical bias in the step size $h$, and Corollary~\ref{lem: w1 limit} gives the explicit limit of $\mathcal W_1(\pi_h,\pi)/h$ as $h \to 0$.

We leave the following questions for future work: to prove the results without the damping condition \eqref{eq: friction condition}, and to make the dependence of the constants on the damping rate $\gamma$ explicit.

\subsection*{Declaration of AI use}
The authors have used Claude Fable 5.2 to write the code for generating the figures and to verify the correctness of the proofs.

\bibliographystyle{plainnat}
\bibliography{references}

@article{elworthy1994formulae,
  title={Formulae for the derivatives of heat semigroups},
  author={Elworthy, K. David and Li, Xue-Mei},
  journal={Journal of Functional Analysis},
  volume={125},
  number={1},
  pages={252--286},
  year={1994},
  doi={10.1006/jfan.1994.1124}
}

@book{bismut1984large,
  title={Large Deviations and the {Malliavin} Calculus},
  author={Bismut, Jean-Michel},
  series={Progress in Mathematics},
  volume={45},
  publisher={Birkh{\"a}user},
  address={Boston},
  year={1984}
}

@article{gorham2019measuring,
  title={Measuring sample quality with diffusions},
  author={Gorham, Jackson and Duncan, Andrew B. and Vollmer, Sebastian J. and Mackey, Lester},
  journal={The Annals of Applied Probability},
  volume={29},
  number={5},
  pages={2884--2928},
  year={2019},
  doi={10.1214/19-AAP1467},
  note={arXiv:1611.06972}
}

@misc{ye2026mse,
  title={Mean square error analysis of stochastic {Runge--Kutta} integrators},
  author={Ye, Xuda},
  howpublished={arXiv:2609.11528},
  year={2026}
}

@misc{ye2026sglm,
  title={Stochastic gradient {Leimkuhler--Matthews} integrator with second order discretization error},
  author={Ye, Xuda and Wang, Yuliang and Lin, Guang and Kong, Cindy Xiangrui},
  howpublished={arXiv:2610.24781},
  year={2026}
}

@article{schuh2024global,
  title={Global contractivity for {Langevin} dynamics with distribution-dependent forces and uniform in time propagation of chaos},
  author={Schuh, Katharina},
  journal={Annales de l'Institut Henri Poincar{\'e}, Probabilit{\'e}s et Statistiques},
  volume={60},
  number={2},
  pages={753--789},
  year={2024},
  doi={10.1214/22-AIHP1337},
  note={arXiv:2206.03082}
}

@article{desvillettes2001trend,
  title={On the trend to global equilibrium in spatially inhomogeneous entropy-dissipating systems: the linear {Fokker--Planck} equation},
  author={Desvillettes, Laurent and Villani, C{\'e}dric},
  journal={Communications on Pure and Applied Mathematics},
  volume={54},
  number={1},
  pages={1--42},
  year={2001},
  doi={10.1002/1097-0312(200101)54:1<1::AID-CPA1>3.0.CO;2-Q}
}

@article{eckmann2003spectral,
  title={Spectral properties of hypoelliptic operators},
  author={Eckmann, Jean-Pierre and Hairer, Martin},
  journal={Communications in Mathematical Physics},
  volume={235},
  number={2},
  pages={233--253},
  year={2003},
  doi={10.1007/s00220-003-0805-9},
  note={arXiv:math-ph/0207046}
}

@article{herau2004isotropic,
  title={Isotropic hypoellipticity and trend to equilibrium for the {Fokker--Planck} equation with a high-degree potential},
  author={H{\'e}rau, Fr{\'e}d{\'e}ric and Nier, Francis},
  journal={Archive for Rational Mechanics and Analysis},
  volume={171},
  number={2},
  pages={151--218},
  year={2004},
  doi={10.1007/s00205-003-0276-3}
}

@book{helffer2005hypoelliptic,
  title={Hypoelliptic Estimates and Spectral Theory for {Fokker--Planck} Operators and {Witten} {Laplacians}},
  author={Helffer, Bernard and Nier, Francis},
  series={Lecture Notes in Mathematics},
  volume={1862},
  publisher={Springer},
  address={Berlin, Heidelberg},
  year={2005},
  doi={10.1007/b104762}
}

@article{wu2001large,
  title={Large and moderate deviations and exponential convergence for stochastic damping {Hamiltonian} systems},
  author={Wu, Liming},
  journal={Stochastic Processes and their Applications},
  volume={91},
  number={2},
  pages={205--238},
  year={2001},
  doi={10.1016/S0304-4149(00)00061-2}
}

@article{mattingly2002ergodicity,
  title={Ergodicity for {SDEs} and approximations: locally {Lipschitz} vector fields and degenerate noise},
  author={Mattingly, Jonathan C. and Stuart, Andrew M. and Higham, Desmond J.},
  journal={Stochastic Processes and their Applications},
  volume={101},
  number={2},
  pages={185--232},
  year={2002},
  doi={10.1016/S0304-4149(02)00150-3}
}

@article{villani2009hypocoercivity,
  title={Hypocoercivity},
  author={Villani, C{\'e}dric},
  journal={Memoirs of the American Mathematical Society},
  volume={202},
  number={950},
  year={2009},
  doi={10.1090/S0065-9266-09-00567-5},
  note={arXiv:math/0609050}
}

@article{dolbeault2015hypocoercivity,
  title={Hypocoercivity for linear kinetic equations conserving mass},
  author={Dolbeault, Jean and Mouhot, Cl{\'e}ment and Schmeiser, Christian},
  journal={Transactions of the American Mathematical Society},
  volume={367},
  number={6},
  pages={3807--3828},
  year={2015},
  doi={10.1090/S0002-9947-2015-06012-7},
  note={arXiv:1005.1495}
}

@article{baudoin2017bakry,
  title={{Bakry--{\'E}mery} meet {Villani}},
  author={Baudoin, Fabrice},
  journal={Journal of Functional Analysis},
  volume={273},
  number={7},
  pages={2275--2291},
  year={2017},
  doi={10.1016/j.jfa.2017.06.021},
  note={arXiv:1308.4938}
}

@article{priola2006formulae,
  title={Formulae for the derivatives of degenerate diffusion semigroups},
  author={Priola, Enrico},
  journal={Journal of Evolution Equations},
  volume={6},
  number={4},
  pages={577--600},
  year={2006},
  doi={10.1007/s00028-006-0273-8}
}

@article{priola2009global,
  title={Global {Schauder} estimates for a class of degenerate {Kolmogorov} equations},
  author={Priola, Enrico},
  journal={Studia Mathematica},
  volume={194},
  number={2},
  pages={117--153},
  year={2009},
  doi={10.4064/sm194-2-2},
  note={arXiv:0705.2810}
}

@article{fang2019multivariate,
  title={Multivariate approximations in {Wasserstein} distance by {Stein's} method and {Bismut's} formula},
  author={Fang, Xiao and Shao, Qi-Man and Xu, Lihu},
  journal={Probability Theory and Related Fields},
  volume={174},
  number={3--4},
  pages={945--979},
  year={2019},
  doi={10.1007/s00440-018-0874-5},
  note={arXiv:1801.07815}
}

@incollection{bakry1985diffusions,
  title={Diffusions hypercontractives},
  author={Bakry, Dominique and {\'E}mery, Michel},
  booktitle={S{\'e}minaire de Probabilit{\'e}s XIX 1983/84},
  series={Lecture Notes in Mathematics},
  pages={177--206},
  publisher={Springer},
  address={Berlin, Heidelberg},
  year={1985},
  doi={10.1007/BFb0075847}
}

@book{bakry2014analysis,
  title={Analysis and Geometry of {Markov} Diffusion Operators},
  author={Bakry, Dominique and Gentil, Ivan and Ledoux, Michel},
  series={Grundlehren der mathematischen Wissenschaften},
  volume={348},
  publisher={Springer},
  address={Cham},
  year={2014},
  doi={10.1007/978-3-319-00227-9}
}

@article{monmarche2019generalized,
  title={Generalized {$\Gamma$} calculus and application to interacting particles on a graph},
  author={Monmarch{\'e}, Pierre},
  journal={Potential Analysis},
  volume={50},
  number={3},
  pages={439--466},
  year={2019},
  doi={10.1007/s11118-018-9689-3},
  note={arXiv:1510.05936}
}

@article{albritton2024variational,
  title={Variational methods for the kinetic {Fokker--Planck} equation},
  author={Albritton, Dallas and Armstrong, Scott and Mourrat, Jean-Christophe and Novack, Matthew},
  journal={Analysis \& PDE},
  volume={17},
  number={6},
  pages={1953--2010},
  year={2024},
  doi={10.2140/apde.2024.17.1953},
  note={arXiv:1902.04037}
}

@article{cao2023explicit,
  title={On explicit {$L^2$}-convergence rate estimate for underdamped {Langevin} dynamics},
  author={Cao, Yu and Lu, Jianfeng and Wang, Lihan},
  journal={Archive for Rational Mechanics and Analysis},
  volume={247},
  number={5},
  pages={90},
  year={2023},
  doi={10.1007/s00205-023-01922-4},
  note={arXiv:1908.04746}
}

@article{lu2022explicit,
  title={On explicit {$L^2$}-convergence rate estimate for piecewise deterministic {Markov} processes in {MCMC} algorithms},
  author={Lu, Jianfeng and Wang, Lihan},
  journal={The Annals of Applied Probability},
  volume={32},
  number={2},
  pages={1333--1361},
  year={2022},
  doi={10.1214/21-AAP1710},
  note={arXiv:2007.14927}
}

@article{brigati2023time,
  title={Time averages for kinetic {Fokker--Planck} equations},
  author={Brigati, Giovanni},
  journal={Kinetic and Related Models},
  volume={16},
  number={4},
  pages={524--539},
  year={2023},
  doi={10.3934/krm.2022037},
  note={arXiv:2106.12801}
}

@article{brigati2025construct,
  title={How to construct explicit decay rates for kinetic {Fokker--Planck} equations?},
  author={Brigati, Giovanni and Stoltz, Gabriel},
  journal={SIAM Journal on Mathematical Analysis},
  volume={57},
  number={4},
  pages={3587--3622},
  year={2025},
  doi={10.1137/24M1700351},
  note={arXiv:2302.14506}
}

@misc{brigati2024explicit,
  title={Explicit convergence rates of underdamped {Langevin} dynamics under weighted and weak {Poincar{\'e}--Lions} inequalities},
  author={Brigati, Giovanni and Stoltz, Gabriel and Wang, Andi Q. and Wang, Lihan},
  howpublished={arXiv:2407.16033},
  year={2024}
}

@article{eberle2026nonreversible,
  title={Non-reversible lifts of reversible diffusion processes and relaxation times},
  author={Eberle, Andreas and L{\"o}rler, Francis},
  journal={Probability Theory and Related Fields},
  volume={194},
  number={1--2},
  pages={173--203},
  year={2026},
  doi={10.1007/s00440-024-01308-x},
  note={arXiv:2402.05041}
}

@article{brigati2026hypocoercivity,
  title={Hypocoercivity meets lifts},
  author={Brigati, Giovanni and L{\"o}rler, Francis and Wang, Lihan},
  journal={Kinetic and Related Models},
  volume={20},
  pages={34--55},
  year={2026},
  doi={10.3934/krm.2025020},
  note={arXiv:2412.10890}
}

@misc{eberle2025convergence,
  title={Convergence of non-reversible {Markov} processes via lifting and flow {Poincar{\'e}} inequality},
  author={Eberle, Andreas and Guillin, Arnaud and Hahn, Leo and L{\"o}rler, Francis and Michel, Manon},
  howpublished={arXiv:2503.04238},
  year={2025}
}

@article{bolley2010trend,
  title={Trend to equilibrium and particle approximation for a weakly selfconsistent {Vlasov--Fokker--Planck} equation},
  author={Bolley, Fran{\c{c}}ois and Guillin, Arnaud and Malrieu, Florent},
  journal={ESAIM: Mathematical Modelling and Numerical Analysis},
  volume={44},
  number={5},
  pages={867--884},
  year={2010},
  doi={10.1051/m2an/2010045},
  note={arXiv:0906.1417}
}

@inproceedings{cheng2018underdamped,
  title={Underdamped {Langevin} {MCMC}: a non-asymptotic analysis},
  author={Cheng, Xiang and Chatterji, Niladri S. and Bartlett, Peter L. and Jordan, Michael I.},
  booktitle={Proceedings of the 31st Conference On Learning Theory},
  series={Proceedings of Machine Learning Research},
  volume={75},
  pages={300--323},
  year={2018},
  note={arXiv:1707.03663}
}

@article{dalalyan2020sampling,
  title={On sampling from a log-concave density using kinetic {Langevin} diffusions},
  author={Dalalyan, Arnak S. and Riou-Durand, Lionel},
  journal={Bernoulli},
  volume={26},
  number={3},
  pages={1956--1988},
  year={2020},
  doi={10.3150/19-BEJ1178},
  note={arXiv:1807.09382}
}

@article{eberle2016reflection,
  title={Reflection couplings and contraction rates for diffusions},
  author={Eberle, Andreas},
  journal={Probability Theory and Related Fields},
  volume={166},
  number={3--4},
  pages={851--886},
  year={2016},
  doi={10.1007/s00440-015-0673-1},
  note={arXiv:1305.1233}
}

@article{eberle2019couplings,
  title={Couplings and quantitative contraction rates for {Langevin} dynamics},
  author={Eberle, Andreas and Guillin, Arnaud and Zimmer, Raphael},
  journal={The Annals of Probability},
  volume={47},
  number={4},
  pages={1982--2010},
  year={2019},
  doi={10.1214/18-AOP1299},
  note={arXiv:1703.01617}
}

@article{monmarche2021high,
  title={High-dimensional {MCMC} with a standard splitting scheme for the underdamped {Langevin} diffusion},
  author={Monmarch{\'e}, Pierre},
  journal={Electronic Journal of Statistics},
  volume={15},
  number={2},
  pages={4117--4166},
  year={2021},
  doi={10.1214/21-EJS1888},
  note={arXiv:2007.05455}
}

@article{leimkuhler2024contraction,
  title={Contraction and convergence rates for discretized kinetic {Langevin} dynamics},
  author={Leimkuhler, Benedict and Paulin, Daniel and Whalley, Peter A.},
  journal={SIAM Journal on Numerical Analysis},
  volume={62},
  number={3},
  pages={1226--1258},
  year={2024},
  doi={10.1137/23M1556289},
  note={arXiv:2302.10684}
}

@article{sanz2021wasserstein,
  title={Wasserstein distance estimates for the distributions of numerical approximations to ergodic stochastic differential equations},
  author={Sanz-Serna, Jes{\'u}s Mar{\'\i}a and Zygalakis, Konstantinos C.},
  journal={Journal of Machine Learning Research},
  volume={22},
  number={242},
  pages={1--37},
  year={2021},
  note={arXiv:2104.12384}
}

@article{paulin2024correction,
  title={Correction to ``{Wasserstein} distance estimates for the distributions of numerical approximations to ergodic stochastic differential equations''},
  author={Paulin, Daniel and Whalley, Peter A.},
  journal={Journal of Machine Learning Research},
  volume={25},
  number={376},
  pages={1--9},
  year={2024},
  note={arXiv:2402.08711}
}

@misc{paulin2026theoretical,
  title={Theoretical guarantees for stochastic gradient sampling methods via {Gaussian} convolution inequalities},
  author={Paulin, Daniel and Whalley, Peter A.},
  howpublished={arXiv:2604.24632},
  year={2026}
}

@book{villani2021topics,
  title={Topics in Optimal Transportation},
  author={Villani, C{\'e}dric},
  series={Graduate Studies in Mathematics},
  volume={58},
  publisher={American Mathematical Society},
  address={Providence, RI},
  year={2003},
  doi={10.1090/gsm/058}
}

@book{revuz1999continuous,
  title={Continuous Martingales and {Brownian} Motion},
  author={Revuz, Daniel and Yor, Marc},
  edition={3},
  series={Grundlehren der mathematischen Wissenschaften},
  volume={293},
  publisher={Springer},
  address={Berlin},
  year={1999},
  doi={10.1007/978-3-662-06400-9}
}

@misc{schuh2024convergence,
  title={Convergence of kinetic {Langevin} samplers for non-convex potentials},
  author={Schuh, Katharina and Whalley, Peter A.},
  howpublished={arXiv:2405.09992},
  year={2024}
}

@article{chau2022stochastic,
  title={Stochastic gradient {Hamiltonian} {Monte} {Carlo} for non-convex learning},
  author={Chau, Huy N. and R{\'a}sonyi, Mikl{\'o}s},
  journal={Stochastic Processes and their Applications},
  volume={149},
  pages={341--368},
  year={2022},
  doi={10.1016/j.spa.2022.04.001},
  note={arXiv:1903.10328}
}

@inproceedings{altschuler2025shifted,
  title={Shifted composition {IV}: toward ballistic acceleration for log-concave sampling},
  author={Altschuler, Jason M. and Chewi, Sinho and Zhang, Matthew S.},
  booktitle={Proceedings of the 58th Annual ACM Symposium on Theory of Computing},
  pages={1739--1750},
  year={2026},
  doi={10.1145/3798129.3800881},
  note={arXiv:2506.23062}
}

@misc{habring2026reverse,
  title={Reverse transportation and {Harnack} inequalities for underdamped {Langevin} in the non-convex setting},
  author={Habring, Andreas and Nagler, Laurenz},
  howpublished={arXiv:2611.02020},
  year={2026}
}

@article{crisan2021uniform,
  title={Uniform in time estimates for the weak error of the {Euler} method for {SDEs} and a pathwise approach to derivative estimates for diffusion semigroups},
  author={Crisan, Dan and Dobson, Paul and Ottobre, Michela},
  journal={Transactions of the American Mathematical Society},
  volume={374},
  number={5},
  pages={3289--3330},
  year={2021},
  doi={10.1090/tran/8301},
  note={arXiv:1905.03524}
}

@article{chak2026reflection,
  title={Reflection coupling for unadjusted generalized {Hamiltonian} {Monte} {Carlo} in the nonconvex stochastic gradient case},
  author={Chak, Martin and Monmarch{\'e}, Pierre},
  journal={IMA Journal of Numerical Analysis},
  volume={46},
  number={4},
  pages={2216--2266},
  year={2026},
  doi={10.1093/imanum/draf045},
  note={arXiv:2310.18774}
}

@article{mattingly2010convergence,
  title={Convergence of numerical time-averaging and stationary measures via {Poisson} equations},
  author={Mattingly, Jonathan C. and Stuart, Andrew M. and Tretyakov, Michael V.},
  journal={SIAM Journal on Numerical Analysis},
  volume={48},
  number={2},
  pages={552--577},
  year={2010},
  doi={10.1137/090770527},
  note={arXiv:0908.4450}
}

@article{kopec2015langevin,
  title={Weak backward error analysis for {Langevin} process},
  author={Kopec, Marie},
  journal={BIT Numerical Mathematics},
  volume={55},
  number={4},
  pages={1057--1103},
  year={2015},
  doi={10.1007/s10543-015-0546-0},
  note={arXiv:1310.2599}
}

@article{leimkuhler2016computation,
  title={The computation of averages from equilibrium and nonequilibrium {Langevin} molecular dynamics},
  author={Leimkuhler, Benedict and Matthews, Charles and Stoltz, Gabriel},
  journal={IMA Journal of Numerical Analysis},
  volume={36},
  number={1},
  pages={13--79},
  year={2016},
  doi={10.1093/imanum/dru056},
  note={arXiv:1308.5814}
}

@misc{lu2025mse,
  title={Mean square error analysis of stochastic gradient and variance-reduced sampling algorithms},
  author={Lu, Jianfeng and Ye, Xuda and Zhou, Zhennan},
  howpublished={arXiv:2511.04413},
  year={2025}
}

@article{bernard2022hypocoercivity,
  title={Hypocoercivity with {Schur} complements},
  author={Bernard, {\'E}tienne and Fathi, Max and Levitt, Antoine and Stoltz, Gabriel},
  journal={Annales Henri Lebesgue},
  volume={5},
  pages={523--557},
  year={2022},
  doi={10.5802/ahl.129},
  note={arXiv:2003.00726}
}

@article{talay2002stochastic,
  title={Stochastic {Hamiltonian} systems: exponential convergence to the invariant measure, and discretization by the implicit {Euler} scheme},
  author={Talay, Denis},
  journal={Markov Processes and Related Fields},
  volume={8},
  number={2},
  pages={163--198},
  year={2002}
}

@article{herau2007short,
  title={Short and long time behavior of the {Fokker--Planck} equation in a confining potential and applications},
  author={H{\'e}rau, Fr{\'e}d{\'e}ric},
  journal={Journal of Functional Analysis},
  volume={244},
  number={1},
  pages={95--118},
  year={2007},
  doi={10.1016/j.jfa.2006.11.013},
  note={arXiv:math/0501363}
}

@article{zhang2010stochastic,
  title={Stochastic flows and {Bismut} formulas for stochastic {Hamiltonian} systems},
  author={Zhang, Xicheng},
  journal={Stochastic Processes and their Applications},
  volume={120},
  number={10},
  pages={1929--1949},
  year={2010},
  doi={10.1016/j.spa.2010.05.015}
}

@article{guillin2012degenerate,
  title={Degenerate {Fokker--Planck} equations: {Bismut} formula, gradient estimate and {Harnack} inequality},
  author={Guillin, Arnaud and Wang, Feng-Yu},
  journal={Journal of Differential Equations},
  volume={253},
  number={1},
  pages={20--40},
  year={2012},
  doi={10.1016/j.jde.2012.03.014},
  note={arXiv:1103.2817}
}

@article{wang2013derivative,
  title={Derivative formula and applications for degenerate diffusion semigroups},
  author={Wang, Feng-Yu and Zhang, Xicheng},
  journal={Journal de Math{\'e}matiques Pures et Appliqu{\'e}es},
  volume={99},
  number={6},
  pages={726--740},
  year={2013},
  doi={10.1016/j.matpur.2012.10.007},
  note={arXiv:1107.0096}
}

@article{wang2014integration,
  title={Integration by parts formula and shift {Harnack} inequality for stochastic equations},
  author={Wang, Feng-Yu},
  journal={The Annals of Probability},
  volume={42},
  number={3},
  pages={994--1019},
  year={2014},
  doi={10.1214/13-AOP875},
  note={arXiv:1203.4023}
}

\appendix

\section{Proofs of the semigroup estimates}
\label{app: proofs}

Throughout this appendix, $Z_t$ is the solution of \eqref{eq: uld} at the time $t$ from the initial state $z$, along a fixed Brownian path, and $J_t = \partial_z Z_t$. For a probability measure $\mathbb P$ with expectation $\mathbb E$ and for $p \Ge 1$, $L^p(\mathbb P)$ is the space of random variables $\xi$ with $\mathbb E|\xi|^p < \infty$, normed by $(\mathbb E|\xi|^p)^{\frac1p}$.

\subsection{Proofs of Lemma~\ref{lem: variations} and Lemma~\ref{lem: derivatives}}
\label{app: variations}

\begin{proof}[Proof of Lemma~\ref{lem: variations}]
    Write the coefficient of \eqref{eq: first variation} as
    \begin{equation}
        Q(x) := \begin{pmatrix} 0 & I_d \\ -\nabla^2U(x) & -\gamma I_d \end{pmatrix} .
        \label{eq: variation coefficient}
    \end{equation}
    By \eqref{eq: U bounds}, $\norm{Q(x)} \Le 1 + M + \gamma$ for every $x \in \mathbb R^d$, so Gr\"onwall's inequality applied to \eqref{eq: first variation} gives
    \begin{equation}
        \norm{J_t} \Le e^{(1+M+\gamma)t} , \qquad t \Ge 0 ,
        \label{eq: jacobian growth}
    \end{equation}
    hence $\norm{J_t} \Le C$ for $t \Le 1$. By \eqref{eq: U bounds} and \eqref{eq: tensor norm}, the term of \eqref{eq: second variation} with the factor $\nabla^3U(X_t)$ satisfies
    \begin{equation*}
        \biggl\lVert \begin{pmatrix} 0 & 0 \\ \nabla^3U(X_t)\bigl[ (J_t w)^x \bigr] & 0 \end{pmatrix} J_t \biggr\rVert
        \Le \norm{\nabla^3U(X_t)}\, \bigl| (J_tw)^x \bigr| \, \norm{J_t}
        \Le M \norm{J_t}^2 |w| \Le C|w| , \qquad t \Le 1 ,
    \end{equation*}
    so Gr\"onwall's inequality applied to \eqref{eq: second variation} gives $\norm{K_t[w]} \Le C t |w|\, e^{(1+M+\gamma)t} \Le C t |w|$, which proves \eqref{eq: variation bounds}.
\end{proof}

\begin{proof}[Proof of Lemma~\ref{lem: derivatives}]
    To prove \eqref{eq: gradient bound}, let $t \Ge 0$ and $w \in \mathbb R^{2d}$. The bound \eqref{eq: jacobian growth} holds at every initial state $z \in \mathbb R^{2d}$ and along every Brownian path, so the finite differences of $f(Z_t)$ in the direction $w$ are bounded by $e^{(1+M+\gamma)t} \norm{\nabla f}_\infty |w|$ and converge pathwise to $\nabla f(Z_t)\bigl[ J_tw \bigr]$. Dominated convergence therefore gives
    \begin{equation}
        \nabla (e^{t\mathcal L} f)(z) = \mathbb E^z\bigl[ J_t^\top \nabla f(Z_t) \bigr] ,
        \label{eq: gradient formula}
    \end{equation}
    which, with \eqref{eq: variation bounds}, gives \eqref{eq: gradient bound} for $0 \Le t \Le 1$. The integrand of \eqref{eq: gradient formula} is continuous in $z$ pathwise and dominated by $e^{(1+M+\gamma)t}\norm{\nabla f}_\infty$, so $e^{t\mathcal L}f \in C^1(\mathbb R^{2d})$.

    To prove \eqref{eq: second derivative pass}, let $0 \Le t \Le 1$ and $w,w' \in \mathbb R^{2d}$. Since $z \mapsto Z_t$ is twice continuously differentiable pathwise, the derivative of $\nabla f(Z_t)\bigl[ J_tw \bigr]$ in the initial state along $w'$ is $\nabla^2 f(Z_t)\bigl[ J_tw, J_tw' \bigr] + \nabla f(Z_t)\bigl[ K_t[w']\, w \bigr]$. By \eqref{eq: variation bounds}, at every initial state $z$,
    \begin{equation*}
        \bigl| \nabla^2 f(Z_t)\bigl[ J_tw, J_tw' \bigr] \bigr| \Le C \norm{\nabla^2f}_\infty |w| |w'| , \qquad
        \bigl| \nabla f(Z_t)\bigl[ K_t[w']\, w \bigr] \bigr| \Le C t \norm{\nabla f}_\infty |w| |w'| ,
    \end{equation*}
    so the finite differences of $\nabla f(Z_t)\bigl[ J_tw \bigr]$ in the direction $w'$ are bounded by $C\bigl( \norm{\nabla^2f}_\infty + t \norm{\nabla f}_\infty \bigr) |w| |w'|$. Dominated convergence therefore gives
    \begin{equation}
        \nabla^2 (e^{t\mathcal L} f)(z)[w,w'] = \mathbb E^z\Bigl[ \nabla^2 f(Z_t)\bigl[ J_tw, J_tw' \bigr] + \nabla f(Z_t)\bigl[ K_t[w']\, w \bigr] \Bigr]
        \label{eq: hessian formula}
    \end{equation}
    and thus \eqref{eq: second derivative pass}. The integrand of \eqref{eq: hessian formula} is continuous in $z$ and bounded by $C\bigl( \norm{\nabla^2f}_\infty + t \norm{\nabla f}_\infty \bigr) |w| |w'|$, so dominated convergence gives $e^{t\mathcal L}f \in C^2(\mathbb R^{2d})$.
\end{proof}

\subsection{Proof of Lemma~\ref{lem: bismut}}
\label{app: bismut}

\begin{proof}[Proof of Lemma~\ref{lem: bismut}]
    By \eqref{eq: U bounds}, \eqref{eq: eta} and \eqref{eq: hermite bounds}, and since $t \Le 1$,
    \begin{equation}
        |(\mathcal A_Z w)(s)| \Le C \Bigl( t^{-2} |w^x| + t^{-1} |w^v| \Bigr) , \qquad s \in [0,t] ,
        \label{eq: eta pointwise}
    \end{equation}
    and integrating the square of \eqref{eq: eta pointwise} over $[0,t]$ gives \eqref{eq: eta size}.

    To prove \eqref{eq: bismut}, fix $\varepsilon \neq 0$ and define the \emph{$\varepsilon$-shift} of the solution $Z_s = (X_s,V_s)$ by
    \begin{equation*}
        X^\varepsilon_s := X_s + \varepsilon w^\circleddash(s) , \qquad
        V^\varepsilon_s := V_s + \varepsilon \dot w^\circleddash(s) , \qquad s \in [0,t] ,
    \end{equation*}
    so that $\D X^\varepsilon_s = V^\varepsilon_s \,\D s$ and $(X^\varepsilon_0, V^\varepsilon_0) = z + \varepsilon w$. Replacing $\nabla^2U(X_s)\, w^\circleddash(s)$ in \eqref{eq: eta} by a finite difference, define
    \begin{equation}
        (\mathcal A_Z^{\varepsilon} w)(s) := \frac{1}{\sqrt{2\gamma}} \biggl( \frac{\nabla U(X^\varepsilon_s) - \nabla U(X_s)}{\varepsilon} + \gamma\, \dot w^\circleddash(s) + \ddot w^\circleddash(s) \biggr) .
        \label{eq: eta eps}
    \end{equation}
    By \eqref{eq: U bounds} and \eqref{eq: hermite bounds}, $|(\mathcal A_Z^{\varepsilon} w)(s)|$ satisfies \eqref{eq: eta pointwise} uniformly in $\varepsilon$. Since $\nabla^2U$ is continuous and bounded on $\mathbb R^d$, $(\mathcal A_Z^{\varepsilon} w)(s) \to (\mathcal A_Z w)(s)$ pointwise as $\varepsilon \to 0$, and dominated convergence gives
    \begin{equation}
    	\lim_{\varepsilon \to 0} \mathbb E \biggl[ \int_0^t \bigl| (\mathcal A_Z^{\varepsilon} w)(s) - (\mathcal A_Z w)(s) \bigr|^2 \D s \biggr] = 0 .
        \label{eq: eta eps limit}
    \end{equation}
    By the definition \eqref{eq: eta eps},
    \begin{equation*}
        \D V^\varepsilon_s = \D V_s + \varepsilon \ddot w^\circleddash(s)\,\D s
        = -\bigl( \nabla U(X^\varepsilon_s) + \gamma V^\varepsilon_s \bigr) \D s + \sqrt{2\gamma}\, \bigl( \D B_s + \varepsilon\, (\mathcal A_Z^{\varepsilon} w)(s) \,\D s \bigr) .
    \end{equation*}
    Thus $(X_s^\varepsilon,V_s^\varepsilon)$ solves \eqref{eq: uld} from the initial state $z + \varepsilon w$, driven by $(B^\varepsilon_s)_{s \in [0,t]}$, where
    \begin{equation*}
        B^\varepsilon_s := B_s + \varepsilon \int_0^s (\mathcal A_Z^{\varepsilon} w)(r) \,\D r , \qquad s \in [0,t] .
    \end{equation*}
    In parallel, let $\bar Z^\varepsilon_s = (\bar X^\varepsilon_s, \bar V^\varepsilon_s)$, $s \in [0,t]$, be the solution of \eqref{eq: uld} from the initial state $z + \varepsilon w$, driven by the unmodified Brownian motion $(B_s)_{s \in [0,t]}$, so that $(e^{t\mathcal L} f)(z + \varepsilon w) = \mathbb E \bigl[f(\bar Z^\varepsilon_t)\bigr]$.

    Let $\mathbb P$ denote the distribution of the paths of the Brownian motion $(B_s)_{s \in [0,t]}$. Since $\mathcal A_Z^{\varepsilon} w$ is adapted and bounded, Novikov's condition holds and
    \begin{equation}
        \mathcal D_\varepsilon := \exp\biggl( -\varepsilon \int_0^t (\mathcal A_Z^{\varepsilon} w)(s)\cdot\D B_s - \frac{\varepsilon^2}{2} \int_0^t \bigl| (\mathcal A_Z^{\varepsilon} w)(s) \bigr|^2 \D s \biggr)
        \label{eq: girsanov density}
    \end{equation}
    defines a distribution $\D \mathbb Q_\varepsilon = \mathcal D_\varepsilon \,\D\mathbb P$ of the paths, under which, by Girsanov's theorem, $(B^\varepsilon_s)_{s \in [0,t]}$ is a Brownian motion. The dynamics \eqref{eq: uld} has a unique strong solution, which can be viewed as a map of the initial state and the driving path. In particular, $(X^\varepsilon_s, V^\varepsilon_s)$ is this map evaluated at $z + \varepsilon w$ and $(B^\varepsilon_s)_{s \in [0,t]}$, and $\bar Z^\varepsilon_s$ is this map evaluated at $z + \varepsilon w$ and $(B_s)_{s \in [0,t]}$. Hence the law of $(X^\varepsilon_t, V^\varepsilon_t)$ under $\mathbb Q_\varepsilon$ is the law of $\bar Z^\varepsilon_t$ under $\mathbb P$. By \eqref{eq: hermite}, $Z_t^{\varepsilon} = (X^\varepsilon_t,V^\varepsilon_t) = Z_t$, so that
    \begin{equation}
    	(e^{t\mathcal L} f)(z + \varepsilon w) =
        \mathbb E\bigl[ f(\bar Z^\varepsilon_t) \bigr] = \mathbb E\bigl[ f(Z_t^{\varepsilon})\, \mathcal D_\varepsilon \bigr] = \mathbb E\bigl[ f(Z_t)\, \mathcal D_\varepsilon \bigr] .
        \label{eq: girsanov identity}
    \end{equation}
    The Girsanov identity \eqref{eq: girsanov identity} fully characterizes the perturbation of the semigroup in the direction $w$.

    To differentiate \eqref{eq: girsanov identity} at $\varepsilon = 0$, we show that $(\mathcal D_\varepsilon - 1)/\varepsilon$ converges to $-\int_0^t (\mathcal A_Z w)(s)\cdot\D B_s$ in $L^2(\mathbb P)$. Define the process $\mathcal D_\varepsilon(s)$, $s \in [0,t]$, by
    \begin{equation*}
        \mathcal D_\varepsilon(s) := \exp\biggl( -\varepsilon \int_0^s (\mathcal A_Z^{\varepsilon} w)(r)\cdot\D B_r - \frac{\varepsilon^2}{2} \int_0^s \bigl| (\mathcal A_Z^{\varepsilon} w)(r) \bigr|^2 \D r \biggr) ,
    \end{equation*}
    so that $\mathcal D_\varepsilon(t) = \mathcal D_\varepsilon$. It\^o's formula gives
    \begin{equation}
        \D \mathcal D_\varepsilon(s) = -\varepsilon\, \mathcal D_\varepsilon(s)\, (\mathcal A_Z^{\varepsilon} w)(s)\cdot\D B_s , \qquad \mathcal D_\varepsilon(0) = 1 ,
        \label{eq: density sde}
    \end{equation}
    so $(\mathcal D_\varepsilon - 1)/\varepsilon$ can be written as
    \begin{equation}
        \frac{\mathcal D_\varepsilon - 1}{\varepsilon} = -\int_0^t \mathcal D_\varepsilon(s)\, (\mathcal A_Z^{\varepsilon} w)(s)\cdot\D B_s .
        \label{eq: difference split}
    \end{equation}
    By Novikov's condition applied to $-2\varepsilon\mathcal A_Z^{\varepsilon} w$,
    \begin{equation*}
        \mathbb E \biggl[ \mathcal D_\varepsilon(s)^2 \exp\biggl( -\varepsilon^2 \int_0^s \bigl| (\mathcal A_Z^{\varepsilon} w)(r) \bigr|^2 \D r \biggr) \biggr] = 1 , \qquad s \in [0,t] ,
    \end{equation*}
    so \eqref{eq: eta pointwise} for $\mathcal A_Z^{\varepsilon} w$ and $\mathbb E[\mathcal D_\varepsilon(s)] = 1$ by \eqref{eq: density sde} give
    \begin{equation}
        \mathbb E \bigl| \mathcal D_\varepsilon(s) - 1 \bigr|^2 = \mathbb E \bigl[ \mathcal D_\varepsilon(s)^2 \bigr] - 1 \Le \exp\Bigl( C\varepsilon^2 \bigl( t^{-\frac32}|w^x| + t^{-\frac12}|w^v| \bigr)^2 \Bigr) - 1 , \qquad s \in [0,t] .
        \label{eq: density deviation}
    \end{equation}
    By \eqref{eq: difference split} and the It\^o isometry,
    \begin{equation*}
        \mathbb E \biggl| \frac{\mathcal D_\varepsilon - 1}{\varepsilon} + \int_0^t (\mathcal A_Z w)(s)\cdot\D B_s \biggr|^2
        \Le 2 \int_0^t \mathbb E \Bigl[ \bigl| \mathcal D_\varepsilon(s) - 1 \bigr|^2 \bigl| (\mathcal A_Z^{\varepsilon} w)(s) \bigr|^2 \Bigr] \D s
        + 2\, \mathbb E \biggl[ \int_0^t \bigl| (\mathcal A_Z^{\varepsilon} w)(s) - (\mathcal A_Z w)(s) \bigr|^2 \D s \biggr] ,
    \end{equation*}
    where the right side vanishes as $\varepsilon \to 0$, by \eqref{eq: eta eps limit} and \eqref{eq: density deviation}. Hence
    \begin{equation}
        \lim_{\varepsilon \to 0} \mathbb E \biggl| \frac{\mathcal D_\varepsilon - 1}{\varepsilon} + \int_0^t (\mathcal A_Z w)(s)\cdot\D B_s \biggr|^2 = 0 .
        \label{eq: difference integral}
    \end{equation}
    By \eqref{eq: girsanov identity} and the Cauchy--Schwarz inequality,
    \begin{equation*}
        \lim_{\varepsilon \to 0} \frac{(e^{t\mathcal L} f)(z + \varepsilon w) - (e^{t\mathcal L} f)(z)}{\varepsilon}
        = \lim_{\varepsilon \to 0} \mathbb E \biggl[ f(Z_t)\, \frac{\mathcal D_\varepsilon - 1}{\varepsilon} \biggr]
        = -\, \mathbb E^z \biggl[ f(Z_t) \int_0^t (\mathcal A_Z w)(s) \cdot \D B_s \biggr] .
    \end{equation*}
    The left side is $\nabla (e^{t\mathcal L} f)(z)[w]$, which proves \eqref{eq: bismut}.
\end{proof}

\subsection{Proof of Lemma~\ref{lem: fluctuation}}
\label{app: fluctuation}

\begin{proof}[Proof of Lemma~\ref{lem: fluctuation}]
    By mollification, it suffices to take the function $\psi \in C^2(\mathbb R^{2d})$ with $\norm{\nabla\psi}_\infty + \norm{\nabla^2\psi}_\infty < \infty$. Let $n \Ge 1$ be an integer and define the functions
    \begin{equation*}
        g_k := e^{\frac{(n-k-1)t}{n}\mathcal L}\psi , \qquad k = 0, \dots, n-1 .
    \end{equation*}
    Given the initial state $z\in \mathbb R^{2d}$, define $N_k := \mathbb E^z\bigl[ \psi(Z_t) \,\big|\, \mathcal F_{kt/n} \bigr]$ for $k = 0, 1, \dots, n$, where $(\mathcal F_s)_{s \in [0,t]}$ is the filtration generated by $(B_s)_{s \in [0,t]}$. Then $N_0 = \mathbb E^z[\psi(Z_t)]$, $N_n = \psi(Z_t)$, and $(N_k)_{k=0}^n$ forms a discrete martingale because $\mathbb E^z\bigl[ N_{k+1} \,\big|\, \mathcal F_{kt/n} \bigr] = N_k$. Hence $(N_k)_{k=0}^n$ has orthogonal increments in $L^2(\mathbb P)$, and
    \begin{equation}
        \mathbb E^z \bigl| \psi(Z_t) - \mathbb E^z \psi(Z_t) \bigr|^2 = \sum_{k=0}^{n-1} \mathbb E^z |N_{k+1} - N_k|^2 .
        \label{eq: martingale increments}
    \end{equation}
    By \eqref{eq: gradient bound} and \eqref{eq: second derivative pass},
    \begin{equation}
    	\norm{\nabla g_k}_\infty \Le C \norm{\nabla\psi}_\infty, \qquad
    	\norm{\nabla^2 g_k}_\infty \Le C \bigl( \norm{\nabla\psi}_\infty + \norm{\nabla^2\psi}_\infty \bigr).
    	\label{eq: gk bounds}
    \end{equation}
    By \eqref{eq: generator}, \eqref{eq: gradient growth} and \eqref{eq: gk bounds}, for $z' = (x',v') \in \mathbb R^{2d}$,
    \begin{align}
        \bigl| (\mathcal L g_k)(z') \bigr| & \Le |v'|\, |\nabla_x g_k(z')| + \bigl( |\nabla U(x')| + \gamma |v'| \bigr) |\nabla_v g_k(z')| + \gamma\, |\Delta_v g_k(z')| \notag \\
        & \Le \bigl( (1+\gamma)|v'| + M(\sqrt d + |x'|) \bigr) \norm{\nabla g_k}_\infty + \gamma\, d\, \norm{\nabla^2 g_k}_\infty \notag \\
        & \Le C \bigl( d + |z'| \bigr) \bigl( \norm{\nabla\psi}_\infty + \norm{\nabla^2\psi}_\infty \bigr) .
        \label{eq: generator growth}
    \end{align}
    
    By the Markov property of \eqref{eq: uld}, $N_{k+1} = g_k(Z_{(k+1)t/n})$. Applying It\^o's formula to $g_k(Z_r)$ gives
    \begin{equation*}
        N_{k+1} = g_k(Z_{kt/n}) + \sqrt{2\gamma} \int_{kt/n}^{(k+1)t/n} \nabla_v g_k(Z_r) \cdot \D B_r
        + \int_{kt/n}^{(k+1)t/n} (\mathcal L g_k)(Z_r) \,\D r .
    \end{equation*}
    Taking the conditional expectation given $\mathcal F_{kt/n}$ gives
    \begin{equation*}
        N_k = g_k(Z_{kt/n}) + \mathbb E^z\biggl[ \int_{kt/n}^{(k+1)t/n} (\mathcal L g_k)(Z_r) \,\D r \,\bigg|\, \mathcal F_{kt/n} \biggr] .
    \end{equation*}
    Subtracting the two identities gives
    \begin{equation*}
        N_{k+1} - N_k
        = \sqrt{2\gamma} \int_{kt/n}^{(k+1)t/n} \nabla_v g_k(Z_r) \cdot \D B_r
        + \int_{kt/n}^{(k+1)t/n} (\mathcal L g_k)(Z_r) \,\D r
        - \mathbb E^z\biggl[ \int_{kt/n}^{(k+1)t/n} (\mathcal L g_k)(Z_r) \,\D r \,\bigg|\, \mathcal F_{kt/n} \biggr] .
    \end{equation*}
    The It\^o isometry, \eqref{eq: gk bounds} and \eqref{eq: generator growth} give
    \begin{equation*}
        \mathbb E^z |N_{k+1} - N_k|^2
        \Le \frac{Ct}{n} \norm{\nabla\psi}_\infty^2
        + \frac{Ct^2}{n^2} \bigl( \norm{\nabla\psi}_\infty + \norm{\nabla^2\psi}_\infty \bigr)^2 \sup_{0 \Le r \Le t} \mathbb E^z \bigl( d + |Z_r| \bigr)^2 ,
    \end{equation*}
    where the supremum is finite. Summing this bound over $k$ in \eqref{eq: martingale increments} and letting $n \to \infty$ gives \eqref{eq: fluctuation}.
\end{proof}

\subsection{Proof of Theorem~\ref{thm: hessian}}
\label{app: hessian}

\begin{proof}[Proof of Theorem~\ref{thm: hessian}]
    For convenience, we prove that, for every $z \in \mathbb R^{2d}$ and $w = (w^x,w^v) \in \mathbb R^{2d}$,
    \begin{equation}
        \bigl| \nabla^2 (e^{2t\mathcal L} f)(z)[w] \bigr|
        \Le C \norm{\nabla f}_\infty \Bigl( t^{-\frac32} |w^x| + t^{-\frac12} |w^v| \Bigr) , \qquad 0<t\Le\frac12.
        \label{eq: hessian doubled}
    \end{equation}
	Put $\psi := e^{t\mathcal L} f$, so that $e^{2t\mathcal L} f = e^{t\mathcal L} \psi$. By Lemma~\ref{lem: derivatives}, $\psi \in C^1(\mathbb R^{2d})$ and $\norm{\nabla \psi}_\infty \Le C \norm{\nabla f}_\infty$.
	
	Given the Brownian motion $(B_s)_{0\Le s\Le t}$,
    let $\mathcal A_Z$ and $w^\circleddash$ be those of Definition~\ref{def: A}. At $z \in \mathbb R^{2d}$, \eqref{eq: bismut centred} applied to $\psi$ with the direction $w$ and with $c := \mathbb E^z\psi(Z_t) = (e^{2t\mathcal L}f)(z)$ gives
    \begin{equation}
        \nabla (e^{2t\mathcal L} f)(z)[w] = -\, \mathbb E^z \biggl[ \bigl( \psi(Z_t) - c \bigr) \int_0^t (\mathcal A_Z w)(s) \cdot \D B_s \biggr] .
        \label{eq: centered gradient}
    \end{equation}
    The initial state $z$ enters the right side of \eqref{eq: centered gradient} through $Z_t$, through the factor $\nabla^2U(X_s)$ of $(\mathcal A_Z w)(s)$ in \eqref{eq: eta}, and through $c = (e^{2t\mathcal L}f)(z)$, whose derivative multiplies a stochastic integral of mean zero. Differentiating \eqref{eq: centered gradient} with respect to $z$ in the direction $w'$ formally gives
    \begin{equation}
        \nabla^2 (e^{2t\mathcal L} f)(z)[w,w']
        = -\, \underbrace{\mathbb E^z \biggl[ \nabla\psi(Z_t)\bigl[ J_t w' \bigr] \int_0^t (\mathcal A_Z w)(s)\cdot\D B_s \biggr]}_{=: \mathcal T_1}
        \; -\; \underbrace{\mathbb E^z \biggl[ \bigl( \psi(Z_t) - c \bigr) \int_0^t \bigl( \partial_{w'} (\mathcal A_Z w)(s) \bigr) \cdot \D B_s \biggr]}_{=: \mathcal T_2}.
        \label{eq: two terms}
    \end{equation}

    To prove \eqref{eq: two terms}, we pass to the limit in the finite differences of both factors of \eqref{eq: centered gradient} in the direction $w'$. By \eqref{eq: gradient bound}, \eqref{eq: variation bounds} and \eqref{eq: hermite bounds}, the finite differences of $\psi(Z_t) - c$ and of $(\mathcal A_Z w)(s)$ are bounded pathwise, uniformly in $s \in [0,t]$. Since $\nabla\psi$ and $\nabla^3U$ are continuous, these finite differences converge pathwise, those of $(\mathcal A_Z w)(s)$ to its derivative in the direction $w'$,
    \begin{equation}
        \partial_{w'}(\mathcal A_Z w)(s) = \frac{1}{\sqrt{2\gamma}} \nabla^3U(X_s)\bigl[ (J_sw')^x \bigr]\, w^\circleddash(s) .
        \label{eq: eta derivative formula}
    \end{equation}
    By dominated convergence and the It\^o isometry, the finite differences of $\psi(Z_t) - c$ and of $\int_0^t (\mathcal A_Z w)(s)\cdot\D B_s$ converge in $L^2(\mathbb P)$, which proves \eqref{eq: two terms}. For a given Brownian motion, $Z_t$, $J_t$, $(\mathcal A_Z w)(s)$ and $\partial_{w'}(\mathcal A_Z w)(s)$ are continuous in $z$, so, by dominated convergence and the It\^o isometry, $\mathcal T_1$ and $\mathcal T_2$ are continuous in $z$, and $e^{2t\mathcal L}f \in C^2(\mathbb R^{2d})$.

    For $\mathcal T_1$, by \eqref{eq: variation bounds}, $\bigl| \nabla\psi(Z_t)\bigl[ J_t w' \bigr] \bigr| \Le \norm{\nabla\psi}_\infty \norm{J_t}|w'| \Le C\norm{\nabla f}_\infty |w'|$ pathwise, and \eqref{eq: eta size} gives
    \begin{equation*}
        \Biggl( \mathbb E^z \biggl| \int_0^t (\mathcal A_Z w)(s)\cdot\D B_s \biggr|^2 \Biggr)^{\frac12}
        = \biggl( \mathbb E^z \biggl[ \int_0^t \bigl| (\mathcal A_Z w)(s) \bigr|^2 \D s \biggr] \biggr)^{\frac12}
        \Le C \Bigl( t^{-\frac32}|w^x| + t^{-\frac12}|w^v| \Bigr) ,
    \end{equation*}
    so the Cauchy--Schwarz inequality gives
    \begin{equation}
        |\mathcal T_1| \Le C \norm{\nabla f}_\infty\, |w'| \Bigl( t^{-\frac32}|w^x| + t^{-\frac12}|w^v| \Bigr) .
        \label{eq: T1}
    \end{equation}

    For $\mathcal T_2$, by \eqref{eq: eta derivative formula}, \eqref{eq: variation bounds} and \eqref{eq: hermite bounds}, $|\partial_{w'}(\mathcal A_Z w)(s)| \Le C|w'|\bigl( |w^x| + t|w^v| \bigr)$ for $s \in [0,t]$, and integrating the square of this bound over $[0,t]$ gives
    \begin{equation}
        \biggl( \int_0^t \bigl| \partial_{w'}(\mathcal A_Z w)(s) \bigr|^2 \D s \biggr)^{\frac12}
        \Le C \sqrt t\, |w'| \bigl( |w^x| + t|w^v| \bigr) .
        \label{eq: eta derivative}
    \end{equation}
    Since $c = \mathbb E^z\psi(Z_t)$, \eqref{eq: fluctuation} gives
    \begin{equation*}
        \bigl( \mathbb E^z \bigl| \psi(Z_t) - c \bigr|^2 \bigr)^{\frac12} \Le C \norm{\nabla\psi}_\infty \sqrt t \Le C \norm{\nabla f}_\infty \sqrt t ,
    \end{equation*}
    so that the Cauchy--Schwarz inequality and \eqref{eq: eta derivative} give
    \begin{equation}
        |\mathcal T_2| \Le C \norm{\nabla f}_\infty\, t\, |w'| \bigl( |w^x| + t |w^v| \bigr) .
        \label{eq: T2}
    \end{equation}
    Since $t \Le \frac12$, inserting \eqref{eq: T1} and \eqref{eq: T2} into \eqref{eq: two terms} and taking the supremum over $|w'| \Le 1$ gives \eqref{eq: hessian doubled}.
\end{proof}

\subsection{Discussion of Lemma~\ref{lem: contraction}}
\label{app: contraction}

We apply Theorem~5 of \citet{schuh2024global} to \eqref{eq: uld} with the gradient split as $\nabla U(x) = \frac m2 x - \beta(x)$ into a linear part and a remainder $\beta$ of Lipschitz constant $L_\beta$. In the parameters of that paper, the inverse mass is $1$, the force is $-\nabla U$, the linear coefficient is $\frac m2 I_d$, whose smallest and largest eigenvalues are $\frac m2$, and the remainder is $\beta$. Assumptions~\ref{ass: U} and~\ref{ass: convexity} supply the remaining hypotheses of the theorem, that the force is Lipschitz and the linear coefficient positive definite. The theorem requires two conditions,
\begin{equation}
    \bigl\langle \beta(x) - \beta(x'), x - x' \bigr\rangle \Le 0 \quad \text{for } |x - x'| \Ge R_b , \qquad
    \gamma^2 > \frac{4L_\beta^2}{m} ,
    \label{eq: quoted conditions}
\end{equation}
at a radius $R_b \Ge 0$. By \eqref{eq: U bounds} and \eqref{eq: U convexity}, the one-sided condition of \eqref{eq: quoted conditions} holds at the radius $R_b$ chosen in
\begin{equation}
    L_\beta \Le M + \frac m2 , \qquad R_b = \frac{4(M+m)R}{m} .
    \label{eq: quoted constants}
\end{equation}
The one-sided condition on $\beta$ is equivalent to
\begin{equation}
    \bigl\langle \nabla U(x) - \nabla U(x'), x - x' \bigr\rangle \Ge \frac m2 |x - x'|^2 , \qquad |x - x'| \Ge R_b .
    \label{eq: one-sided gradient}
\end{equation}
To prove \eqref{eq: one-sided gradient}, write
\begin{equation*}
    \bigl\langle \nabla U(x) - \nabla U(x'), x - x' \bigr\rangle
    = \int_0^1 \bigl\langle \nabla^2U\bigl( x' + s(x - x') \bigr)(x - x'), x - x' \bigr\rangle \D s .
\end{equation*}
Since $\nabla^2U \succeq mI_d$ outside the ball of radius $R$ by \eqref{eq: U convexity}, $\nabla^2U \succeq -MI_d$ by \eqref{eq: U bounds}, and
\begin{equation*}
	\int_0^1 \mathbf 1_{\{ |x' + s(x - x')| < R \}} \D s \Le \frac{2R}{|x - x'|},
\end{equation*}
we have the dissipation condition
\begin{equation*}
    \bigl\langle \nabla U(x) - \nabla U(x'), x - x' \bigr\rangle
    \Ge m |x - x'|^2 - 2R(m+M)\, |x - x'|
    \Ge \frac m2 |x - x'|^2 , \qquad |x - x'| \Ge R_b .
\end{equation*}
By \eqref{eq: quoted constants}, $4L_\beta^2/m \Le (2M+m)^2/m$, so \eqref{eq: quoted conditions} holds as long as $\gamma^2 > (2M+m)^2/m$, which is \eqref{eq: friction condition}.

At Dirac initial laws, with $Z_t'$ the solution of \eqref{eq: uld} from the initial state $z'$, the $\mathcal W_1$ bound of \citep[Theorem~5]{schuh2024global} gives
\begin{equation*}
    \mathcal W_1\bigl( \mathrm{Law}(Z_t), \mathrm{Law}(Z_t') \bigr) \Le C e^{-\lambda t}\, |z - z'| , \qquad t \Ge 0 , \quad z, z' \in \mathbb R^{2d} .
\end{equation*}
Kantorovich--Rubinstein duality gives, for every $f \in C^1(\mathbb R^{2d})$ with $\norm{\nabla f}_\infty <\infty$,
\begin{equation*}
    \bigl| (e^{t\mathcal L} f)(z) - (e^{t\mathcal L} f)(z') \bigr| \Le C e^{-\lambda t} \norm{\nabla f}_\infty |z - z'| , \qquad t \Ge 0 , \quad z, z' \in \mathbb R^{2d} ,
\end{equation*}
and the limit $z' \to z$ gives \eqref{eq: contraction}.

Theorem~5 of \citet{schuh2024global} bounds from both sides the radius $R_1$ of the transition region of its metric and builds its rate and its prefactor from the exponent $\Lambda := \frac{m + 2L_\beta}{8} R_1^2$. For $R > 0$ the bounds on $R_1$ give
\begin{equation}
    C_0^{-1} \bigl( 1 + L_\beta R_b^2 \bigr) \Le \Lambda \Le C_0 \bigl( 1 + L_\beta R_b^2 \bigr) ,
    \label{eq: quoted exponent}
\end{equation}
with $C_0$ depending only on $m$, $M$ and $\gamma$. At $R = 0$ the theorem gives $R_1 = 0$ and the rate $0$; since Assumption~\ref{ass: convexity} with $R = 0$ implies Assumption~\ref{ass: convexity} with $R = 1$, the theorem is applied with the radius $R_b$ of \eqref{eq: quoted constants} at $R = 1$ in that case. The rate $\lambda$ and the prefactor $C$ of \eqref{eq: contraction} that the theorem gives satisfy $C_0^{-1} e^{-\Lambda} \Le \lambda \Le C_0 e^{-\Lambda}$ and $C \Le C_0 e^{\Lambda}$, after enlarging $C_0$. Since $L_\beta \Le M + \frac m2$ by \eqref{eq: quoted constants}, \eqref{eq: quoted exponent} gives
\begin{equation*}
    \lambda \Ge C_0^{-1} \exp\biggl( -C_0 \Bigl( 1 + \Bigl( M + \frac m2 \Bigr) R_b^2 \Bigr) \biggr) , \qquad
    C \Le C_0 \exp\biggl( C_0 \Bigl( 1 + \Bigl( M + \frac m2 \Bigr) R_b^2 \Bigr) \biggr) .
\end{equation*}
Hence \eqref{eq: contraction} holds with $\lambda$ and $C$ replaced by these two bounds. Since $R_b$ is proportional to $R$, the two bounds are determined by $m$, $M$, $R$ and $\gamma$, exponential in $R^2$ and free of the dimension. Remark~6 of \citet{schuh2024global} states that the rate does not depend on the dimension, and the formula for the prefactor in \citep[Theorem~5]{schuh2024global} carries no dimension either.

\subsection{Proof of Corollary~\ref{cor: uniform}}
\label{app: uniform}

\begin{proof}[Proof of Corollary~\ref{cor: uniform}]
    For $t \Le 1$, Theorem~\ref{thm: hessian} gives $e^{t\mathcal L}f \in C^2(\mathbb R^{2d})$, and \eqref{eq: hessian anisotropic} at $w = (0,w^v)$ gives \eqref{eq: uniform}, and at $w = (w^x,0)$ gives \eqref{eq: uniform xx}.
    For $t \Ge 1$, write $e^{t\mathcal L}f = e^{\mathcal L}(e^{(t-1)\mathcal L}f)$, where $e^{(t-1)\mathcal L}f \in C^1(\mathbb R^{2d})$ by Lemma~\ref{lem: derivatives}. Theorem~\ref{thm: hessian} at the time $1$ applied to $e^{(t-1)\mathcal L}f$ gives $e^{t\mathcal L}f \in C^2(\mathbb R^{2d})$ and
    \begin{equation}
        \norm{\nabla^2 (e^{t\mathcal L}f)}_\infty \Le C \norm{\nabla (e^{(t-1)\mathcal L}f)}_\infty .
        \label{eq: hessian time one}
    \end{equation}
    By \eqref{eq: contraction} at the time $t-1$,
    \begin{equation}
        \norm{\nabla (e^{(t-1)\mathcal L}f)}_\infty \Le C e^{-\lambda(t-1)} \norm{\nabla f}_\infty .
        \label{eq: contraction shifted}
    \end{equation}
    Combining \eqref{eq: hessian time one} and \eqref{eq: contraction shifted}, absorbing $e^{\lambda}$ into $C$ and using $\min(1,t) = 1$ proves \eqref{eq: uniform} and \eqref{eq: uniform xx} for $t \Ge 1$.
\end{proof}

\subsection{Proof of Lemma~\ref{lem: sharp}}
\label{app: sharp}

\begin{proof}[Proof of Lemma~\ref{lem: sharp}]
    For $U(x) = \frac12 x^2$, the coefficient $Q(x)$ of \eqref{eq: variation coefficient} is the constant matrix
    \begin{equation*}
        Q = \begin{pmatrix} 0 & 1 \\ -1 & -\gamma \end{pmatrix} ,
    \end{equation*}
    so \eqref{eq: uld} is solved by
    \begin{equation*}
        Z_t = e^{tQ} z + G_t , \qquad G_t := \sqrt{2\gamma} \int_0^t e^{(t-s)Q} \begin{pmatrix} 0 \\ 1 \end{pmatrix} \D B_s .
    \end{equation*}
    For $r \Ge 0$, write the entries of $e^{rQ}$ as
    \begin{equation*}
        e^{rQ} = \begin{pmatrix} (e^{rQ})_{xx} & (e^{rQ})_{vx} \\ (e^{rQ})_{xv} & (e^{rQ})_{vv} \end{pmatrix}  \in \mathbb R^{2\times 2}.
    \end{equation*}
    The vector $G_t = (G^x_t, G^v_t)$ is centered Gaussian, and by the It\^o isometry its blocks have the variances
    \begin{equation*}
        \sigma_x(t)^2 := \mathbb E|G^x_t|^2 = 2\gamma \int_0^t (e^{rQ})_{vx}^2 \,\D r , \qquad
        \sigma_v(t)^2 := \mathbb E|G^v_t|^2 = 2\gamma \int_0^t (e^{rQ})_{vv}^2 \,\D r .
    \end{equation*}
    The expansion $e^{rQ} = I_2 + rQ + \mathcal O(r^2)$ for small $r$ then gives
    \begin{equation}
        \sigma_x(t)^2 = \frac{2\gamma}{3}\, t^3 \bigl( 1 + \mathcal O(t) \bigr) , \qquad
        \sigma_v(t)^2 = 2\gamma\, t \bigl( 1 + \mathcal O(t) \bigr) .
        \label{eq: variance expansions}
    \end{equation}

    Since $X_t = (e^{tQ})_{xx}\,x + (e^{tQ})_{vx}\,v + G^x_t$, the function $(e^{t\mathcal L}f_x)(z) = \mathbb E^z|X_t|$ is smooth in $z$, and the chain rule gives
    \begin{equation}
        \nabla^2_{xx}(e^{t\mathcal L}f_x)(0) = \frac{2\,(e^{tQ})_{xx}^2}{\sqrt{2\pi}\,\sigma_x(t)} , \qquad
        \nabla^2_{xv}(e^{t\mathcal L}f_x)(0) = \frac{2\,(e^{tQ})_{xx}(e^{tQ})_{vx}}{\sqrt{2\pi}\,\sigma_x(t)} .
        \label{eq: sharp hessian x}
    \end{equation}
    In the same way, $V_t = (e^{tQ})_{xv}\,x + (e^{tQ})_{vv}\,v + G^v_t$ gives
    \begin{equation}
        \nabla^2_{vv}(e^{t\mathcal L}f_v)(0) = \frac{2\,(e^{tQ})_{vv}^2}{\sqrt{2\pi}\,\sigma_v(t)} .
        \label{eq: sharp hessian v}
    \end{equation}
    By \eqref{eq: variance expansions},
    \begin{equation*}
        \frac{2}{\sqrt{2\pi}\,\sigma_x(t)} = \sqrt{\frac{3}{\pi\gamma}}\, t^{-\frac32}\bigl(1 + \mathcal O(t)\bigr) , \qquad
        \frac{2}{\sqrt{2\pi}\,\sigma_v(t)} = \sqrt{\frac{1}{\pi\gamma}}\, t^{-\frac12}\bigl(1 + \mathcal O(t)\bigr) ,
    \end{equation*}
    and inserting the expansions of $(e^{tQ})_{xx}$, $(e^{tQ})_{vx}$ and $(e^{tQ})_{vv}$ into \eqref{eq: sharp hessian x} and \eqref{eq: sharp hessian v} proves \eqref{eq: sharp}.
\end{proof}

\subsection{Proof of Theorem~\ref{thm: third}}
\label{app: third}

For $0 < t \Le 1$, an initial state $z \in \mathbb R^{2d}$, a direction $w \in \mathbb R^{2d}$ and $\varepsilon \neq 0$, recall that $w^\circleddash(s)$ is the cubic interpolant \eqref{eq: hermite}, and
\begin{equation*}
    (X^\varepsilon_s, V^\varepsilon_s) = \bigl( X_s + \varepsilon w^\circleddash(s), V_s + \varepsilon \dot w^\circleddash(s) \bigr)
\end{equation*}
is the $\varepsilon$-shifted solution with the initial state $z+\varepsilon w$ in the proof of Lemma~\ref{lem: bismut}. Let $B^\varepsilon_s$, $\mathbb P$, $\mathbb Q_\varepsilon$ and $\mathcal D_\varepsilon$ be those of the same proof. Let $J^\varepsilon_s \in \mathbb R^{2d\times2d}$ solve the ordinary differential equation
\begin{equation}
    \frac{\D}{\D s} J^\varepsilon_s = Q(X^\varepsilon_s)\, J^\varepsilon_s , \qquad J^\varepsilon_0 = I_{2d} , \qquad s\in[0,t],
    \label{eq: shifted variation}
\end{equation}
with the coefficient matrix $Q(x)\in \mathbb R^{2d\times 2d}$ in \eqref{eq: variation coefficient}. We first compute the limit of $(J^\varepsilon_s - J_s)/\varepsilon$ as $\varepsilon \to 0$.

\begin{lemma}
\label{lem: shifted variation}
    Under Assumption~\ref{ass: U}, let $0 < t \Le 1$, $z \in \mathbb R^{2d}$ and $w \in \mathbb R^{2d}$, and let $w^\circleddash(s)$ be the cubic interpolant \eqref{eq: hermite}. Let $K^\circleddash_s[w] \in \mathbb R^{2d\times2d}$ solve
    \begin{equation}
        \frac{\D}{\D s} K^\circleddash_s[w] = Q(X_s)\, K^\circleddash_s[w] - \begin{pmatrix} 0 & 0 \\ \nabla^3U(X_s)\bigl[ w^\circleddash(s) \bigr] & 0 \end{pmatrix} J_s , \qquad K^\circleddash_0[w] = 0 , \qquad s\in[0,t].
        \label{eq: shift variation}
    \end{equation}
    Then, for every $\varepsilon \neq 0$, the solution $J^\varepsilon_s$ of \eqref{eq: shifted variation} satisfies
    \begin{equation}
        \norm{J^\varepsilon_s} \Le C , \qquad
        \biggl\lVert \frac{J^\varepsilon_s - J_s}{\varepsilon} \biggr\rVert \Le C s \bigl( |w^x| + t|w^v| \bigr) , \qquad \forall s\in [0,t] .
        \label{eq: shifted variation bounds}
    \end{equation}
    Moreover,
    \begin{equation}
        \lim_{\varepsilon \to 0} \frac{J^\varepsilon_s - J_s}{\varepsilon} = K^\circleddash_s[w] , \qquad \forall s \in [0,t] .
        \label{eq: shifted variation limit}
    \end{equation}
    In particular,
    \begin{equation}
        \norm{K^\circleddash_s[w]} \Le C s \bigl( |w^x| + t |w^v| \bigr) , \qquad \forall s \in [0,t] .
        \label{eq: shift variation bound}
    \end{equation}
\end{lemma}

The shift variation $K^\circleddash_s[w]$ differs from the usual second variation $K_s[w]$ in the underlying dynamics: $K^\circleddash_s[w]$ is the derivative in $\varepsilon$ at $\varepsilon = 0$ of the first variation $J^\varepsilon_s$ of the $\varepsilon$-shifted solution $(X^\varepsilon_s, V^\varepsilon_s)$ driven by $(B^\varepsilon_s)_{s \in [0,t]}$, while $K_s[w]$ is the same derivative for the solution with the initial state $z + \varepsilon w$ driven by the unmodified Brownian motion $(B_s)_{s \in [0,t]}$.

\begin{proof}[Proof of Lemma~\ref{lem: shifted variation}]
    By \eqref{eq: U bounds}, $\norm{Q(X^\varepsilon_s) - Q(X_s)} \Le M |\varepsilon|\, |w^\circleddash(s)|$. Gr\"onwall's inequality with the bounds \eqref{eq: variation bounds} and \eqref{eq: hermite bounds} gives \eqref{eq: shifted variation bounds}. Since $U \in C^3(\mathbb R^d)$, $J^\varepsilon_s$ is differentiable in $\varepsilon$, and its derivative at $\varepsilon = 0$ solves \eqref{eq: shifted variation} differentiated in $\varepsilon$ at $\varepsilon = 0$, which is \eqref{eq: shift variation}. This proves \eqref{eq: shifted variation limit}, and \eqref{eq: shifted variation bounds} gives \eqref{eq: shift variation bound}.
\end{proof}

For a bounded function $G \in C(\mathbb R^{2d};\mathbb R^{2d\times 2d})$ and directions $w', w'' \in \mathbb R^{2d}$, we compute the derivative of $\mathbb E^z\bigl[G(Z_t)[J_tw', J_tw'']\bigr]$ in the initial state $z$ along the direction $w$ without differentiating $G$, by the $\varepsilon$-shift approach in the proof of Lemma~\ref{lem: bismut}.

\begin{lemma}
\label{lem: shift}
    Under Assumption~\ref{ass: U}, let $0 < t \Le 1$, $z \in \mathbb R^{2d}$ and $w, w', w'' \in \mathbb R^{2d}$, and let $\mathcal A_Z$ be that of Definition~\ref{def: A} and $K^\circleddash_t[w]$ the solution of \eqref{eq: shift variation} at the time $t$. Let $G \in C(\mathbb R^{2d};\mathbb R^{2d\times 2d})$ be bounded. Then the derivative of $\mathbb E^z\bigl[G(Z_t)[J_tw', J_tw'']\bigr]$ along the direction $w$ is
    \begin{equation}
    \begin{aligned}
        \partial_w\, \mathbb E^z\bigl[ G(Z_t)[J_tw', J_tw''] \bigr]
        & = \mathbb E^z\Bigl[ G(Z_t)\bigl[ K^\circleddash_t[w]\, w', J_tw'' \bigr] + G(Z_t)\bigl[ J_tw', K^\circleddash_t[w]\, w'' \bigr] \Bigr] \\
        & - \mathbb E^z\biggl[ G(Z_t)[J_tw', J_tw''] \int_0^t (\mathcal A_Z w)(s)\cdot\D B_s \biggr] .
    \end{aligned}
        \label{eq: shift}
    \end{equation}
\end{lemma}

\begin{proof}[Proof of Lemma~\ref{lem: shift}]
    As in the proof of Lemma~\ref{lem: bismut}, $(X^\varepsilon_s, V^\varepsilon_s)$ is the solution of \eqref{eq: uld} from $z + \varepsilon w$ driven by $(B^\varepsilon_s)_{s \in [0,t]}$, which is a Brownian motion under $\mathbb Q_\varepsilon$, and $(X^\varepsilon_t, V^\varepsilon_t) = Z_t$. By \eqref{eq: first variation} and \eqref{eq: shifted variation}, the first variation of $(X^\varepsilon_s, V^\varepsilon_s)$ is $J^\varepsilon_s$. Hence, by Girsanov's theorem as in \eqref{eq: girsanov identity},
    \begin{equation*}
        \mathbb E^{z+\varepsilon w}\bigl[ G(Z_t)[J_tw', J_tw''] \bigr] = \mathbb E^z\bigl[ G(Z_t)[J^\varepsilon_tw', J^\varepsilon_tw'']\, \mathcal D_\varepsilon \bigr] ,
    \end{equation*}
    or, in the finite difference form,
    \begin{equation}
    	\frac1{\varepsilon} \Bigl( \mathbb E^{z+\varepsilon w}\bigl[ G(Z_t)[J_tw', J_tw''] \bigr] - \mathbb E^z\bigl[ G(Z_t)[J_tw', J_tw''] \bigr] \Bigr) = \mathcal T_1(\varepsilon) + \mathcal T_2(\varepsilon) ,
        \label{eq: shift difference}
    \end{equation}
    where
    \begin{equation*}
        \mathcal T_1(\varepsilon) := \mathbb E^z\biggl[ \frac{G(Z_t)[J^\varepsilon_tw', J^\varepsilon_tw''] - G(Z_t)[J_tw', J_tw'']}{\varepsilon}\, \mathcal D_\varepsilon \biggr] , \quad
        \mathcal T_2(\varepsilon) := \mathbb E^z\biggl[ G(Z_t)[J_tw', J_tw'']\, \frac{\mathcal D_\varepsilon-1}{\varepsilon} \biggr] .
    \end{equation*}
    Using \eqref{eq: shifted variation bounds}, the limit \eqref{eq: shifted variation limit} and $\mathcal D_\varepsilon \to 1$ in $L^2(\mathbb P)$ by \eqref{eq: density deviation}, the product rule gives
    \begin{equation*}
    	\lim_{\varepsilon\to0}\mathcal T_1(\varepsilon) = \mathbb E^z\Bigl[
    	G(Z_t)\bigl[ K^\circleddash_t[w]\, w', J_tw'' \bigr] + G(Z_t)\bigl[ J_tw', K^\circleddash_t[w]\, w'' \bigr] \Bigr] .
    \end{equation*}
    By \eqref{eq: difference integral},
    \begin{equation*}
    	\lim_{\varepsilon\to0}\mathcal T_2(\varepsilon) = 
    	- \mathbb E^z\biggl[ G(Z_t)[J_tw', J_tw''] \int_0^t (\mathcal A_Z w)(s)\cdot\D B_s \biggr] .
    \end{equation*}
    Letting $\varepsilon \to 0$ in \eqref{eq: shift difference} gives \eqref{eq: shift}.
\end{proof}

Lemma~\ref{lem: shift} with a single factor $J_tw'$ and $G = \nabla f$ gives an alternative proof of Theorem~\ref{thm: hessian}.

\begin{proof}[Proof of Theorem~\ref{thm: third}]
    To prove \eqref{eq: third}, let $w, w', w'' \in \mathbb R^{2d}$ with $w^x = w''^x = 0$. We show that $\nabla^2 (e^{2t\mathcal L} f)[w',w'']$ has a derivative in the direction $w$ at every $z \in \mathbb R^{2d}$, and that this derivative is continuous in $z$ and satisfies
    \begin{equation}
        \Bigl| \partial_w \bigl( \nabla^2 (e^{2t\mathcal L} f)[w',w''] \bigr)(z) \Bigr|
        \Le C\, t^{-1} \norm{\nabla f}_\infty\, |w|\, |w'|\, |w''| , \qquad 0 < t \Le \frac12 .
        \label{eq: third doubled}
    \end{equation}
    Put $\psi := e^{t\mathcal L} f$, so that $e^{2t\mathcal L} f = e^{t\mathcal L}\psi$. By Theorem~\ref{thm: hessian}, $\psi \in C^2(\mathbb R^{2d})$, and \eqref{eq: gradient bound} and \eqref{eq: hessian anisotropic} give
    \begin{equation}
        \norm{\nabla\psi}_\infty \Le C \norm{\nabla f}_\infty , \qquad
        \norm{\nabla^2\psi}_\infty \Le C\, t^{-\frac32} \norm{\nabla f}_\infty .
        \label{eq: third psi}
    \end{equation}
    Hence \eqref{eq: hessian formula} applies to $\psi$ and gives
    \begin{equation}
        \nabla^2 (e^{2t\mathcal L} f)(z)[w',w'']
        = \underbrace{\mathbb E^z\bigl[ \nabla^2\psi(Z_t)[J_tw', J_tw''] \bigr]}_{=: \mathcal T_1}
        + \underbrace{\mathbb E^z\bigl[ \nabla\psi(Z_t)\bigl[ K_t[w'']\, w' \bigr] \bigr]}_{=: \mathcal T_2} .
        \label{eq: third split}
    \end{equation}

    We differentiate $\mathcal T_1$ and $\mathcal T_2$ in the direction $w$ by Lemma~\ref{lem: shift}. For $s \in (0,t]$, let $\mathcal A_Z$ and $K^\circleddash_s[w]$ be those of Lemma~\ref{lem: shift} with $s$ in place of $t$, so that $K^\circleddash_s[w]$ is the solution of \eqref{eq: shift variation} with the terminal time $s$, at the time $s$. Since $w^x = w''^x = 0$, \eqref{eq: first variation}, \eqref{eq: variation bounds}, \eqref{eq: shift variation bound} and \eqref{eq: eta size} give
    \begin{equation}
        \bigl| (J_sw'')^x \bigr| \Le C s |w''| , \qquad
        \norm{K^\circleddash_s[w]} \Le C s^2 |w| , \qquad
        \biggl( \mathbb E^z \biggl| \int_0^s (\mathcal A_Z w)(r)\cdot\D B_r \biggr|^2 \biggr)^{\frac12} \Le C s^{-\frac12} |w| .
        \label{eq: third sizes}
    \end{equation}
    For $\mathcal T_1$, \eqref{eq: hessian anisotropic} with the direction $J_tw''$, \eqref{eq: variation bounds} and \eqref{eq: third sizes} give $\bigl| \nabla^2\psi(Z_t)[J_tw', J_tw''] \bigr| \Le C t^{-\frac12} \norm{\nabla f}_\infty |w'|\, |w''|$. Hence \eqref{eq: third psi}, \eqref{eq: third sizes} at $s = t$ and \eqref{eq: shift} with $G = \nabla^2\psi$ give
    \begin{equation}
        \bigl| \partial_w \mathcal T_1 \bigr| \Le C\, t^{-1} \norm{\nabla f}_\infty\, |w|\, |w'|\, |w''| .
        \label{eq: third first}
    \end{equation}
    For $\mathcal T_2$, \eqref{eq: second variation} gives
    \begin{equation*}
        K_t[w'']\, w' = - \int_0^t J_tJ_s^{-1} \begin{pmatrix} 0 \\ \nabla^3U(X_s)\bigl[ (J_sw'')^x \bigr] (J_sw')^x \end{pmatrix} \D s ,
    \end{equation*}
    where $J_tJ_s^{-1}$ is the first variation at the time $t$ of the solution of \eqref{eq: uld} started from the state $Z_s$ at the time $s$. By the Markov property and \eqref{eq: gradient formula} at the time $t - s$ applied to $\psi$,
    \begin{equation}
        \mathcal T_2
        = - \int_0^t \mathbb E^z\Bigl[ \nabla^3U(X_s)\bigl[ (J_sw')^x, (J_sw'')^x, \nabla_v (e^{(2t-s)\mathcal L} f)(Z_s) \bigr] \Bigr] \D s .
        \label{eq: third integral}
    \end{equation}
    To differentiate $\mathcal T_2$, we write the integrand of \eqref{eq: third integral} at $s \in (0,t]$ as $\mathbb E^z\bigl[ G(Z_s)[J_sw', J_sw''] \bigr]$ with
    \begin{equation*}
        G(z') := \begin{pmatrix} \nabla^3U(x')\bigl[ \nabla_v (e^{(2t-s)\mathcal L} f)(z') \bigr] & 0 \\ 0 & 0 \end{pmatrix} , \qquad z' = (x',v') \in \mathbb R^{2d} .
    \end{equation*}
    By \eqref{eq: U bounds} and Lemma~\ref{lem: derivatives}, the function $G\in C(\mathbb R^{2d};\mathbb R^{2d\times 2d})$ satisfies $\norm{G(z')} \Le C \norm{\nabla f}_\infty$, and \eqref{eq: variation bounds} and \eqref{eq: third sizes} give $\bigl| G(Z_s)[J_sw', J_sw''] \bigr| \Le C s \norm{\nabla f}_\infty |w'|\, |w''|$. Hence \eqref{eq: shift} at the time $s$ and \eqref{eq: third sizes} give
    \begin{equation}
        \biggl| \partial_w\, \mathbb E^z\Bigl[ \nabla^3U(X_s)\bigl[ (J_sw')^x, (J_sw'')^x, \nabla_v (e^{(2t-s)\mathcal L} f)(Z_s) \bigr] \Bigr] \biggr|
        \Le C\, s^{\frac12} \norm{\nabla f}_\infty\, |w|\, |w'|\, |w''| .
        \label{eq: third integrand}
    \end{equation}
    Since the bound \eqref{eq: third integrand} holds at every initial state and is integrable over $s \in (0,t]$, we differentiate \eqref{eq: third integral} in the direction $w$ under the integral over $s$ and obtain
    \begin{equation}
        \bigl| \partial_w \mathcal T_2 \bigr| \Le C\, t^{\frac32} \norm{\nabla f}_\infty\, |w|\, |w'|\, |w''| .
        \label{eq: third second}
    \end{equation}

    By \eqref{eq: third split}, \eqref{eq: third first} and \eqref{eq: third second}, the function $\nabla^2 (e^{2t\mathcal L} f)[w',w'']$ has a derivative in the direction $w$ at every state $z$. This derivative satisfies \eqref{eq: third doubled} since $t \Le \frac12$. For a given Brownian motion, $Z_s$ and $J_s$ are continuous in $z$, and so is $K^\circleddash_s[w]$, since \eqref{eq: shift variation} is a linear equation whose coefficients are bounded and continuous in $z$. Hence, by dominated convergence and the It\^o isometry as in the proof of Theorem~\ref{thm: hessian}, the right side of \eqref{eq: shift} is continuous in $z$ at every time $s \in (0,t]$. With the bound \eqref{eq: third integrand}, which is uniform in $z$ and integrable over $s \in (0,t]$, the derivative of $\nabla^2 (e^{2t\mathcal L} f)[w',w'']$ in the direction $w$ is continuous in $z$. Since $e^{2t\mathcal L}f \in C^2(\mathbb R^{2d})$ by Theorem~\ref{thm: hessian}, the derivatives of $\nabla (e^{2t\mathcal L} f)[w'']$ in the directions $w$ and $w'$ then commute, so $\nabla^2 (e^{2t\mathcal L} f)[w,w'']$ has the same derivative in the direction $w'$. Since $w$ and $w''$ are arbitrary in the velocity and $w'$ is arbitrary in $\mathbb R^{2d}$, $\nabla^2_{vv}(e^{2t\mathcal L}f)$ is continuously differentiable on $\mathbb R^{2d}$. Taking the supremum of \eqref{eq: third doubled} over $z \in \mathbb R^{2d}$ and over $|w^v| \Le 1$, $|w'| \Le 1$, $|w''^v| \Le 1$ gives \eqref{eq: third} at the time $2t \in (0,1]$.
\end{proof}

\subsection{Proof of Lemma~\ref{lem: third sharp}}
\label{app: third sharp}

\begin{proof}[Proof of Lemma~\ref{lem: third sharp}]
    Let $e^{tQ}$, $G^x_t$, $G^v_t$, $\sigma_x(t)$ and $\sigma_v(t)$ be those in the proof of Lemma~\ref{lem: sharp}. Since $X_t = (e^{tQ})_{xx}\,x + (e^{tQ})_{vx}\,v + G^x_t$, the random variable $X_t$ is Gaussian with the mean $\mathbb E^z[X_t] = (e^{tQ})_{xx}\,x + (e^{tQ})_{vx}\,v$ and the variance $\sigma_x(t)^2$. Hence the function $(e^{t\mathcal L}f_x)(z) = \mathbb E^z|X_t|$ is smooth in $z$, and the chain rule gives
    \begin{equation}
        \nabla^3_{xvv}(e^{t\mathcal L}f_x)(z) = -\frac{2\,(e^{tQ})_{xx}(e^{tQ})_{vx}^2\, \mathbb E^z[X_t]}{\sqrt{2\pi}\,\sigma_x(t)^3} \exp\biggl( -\frac{(\mathbb E^z[X_t])^2}{2\sigma_x(t)^2} \biggr) .
        \label{eq: third sharp formula x}
    \end{equation}
    In the same way, $V_t = (e^{tQ})_{xv}\,x + (e^{tQ})_{vv}\,v + G^v_t$ is Gaussian with the mean $\mathbb E^z[V_t] = (e^{tQ})_{xv}\,x + (e^{tQ})_{vv}\,v$ and the variance $\sigma_v(t)^2$, and $(e^{t\mathcal L}f_v)(z) = \mathbb E^z|V_t|$ gives
    \begin{equation}
        \nabla^3_{vvv}(e^{t\mathcal L}f_v)(z) = -\frac{2\,(e^{tQ})_{vv}^3\, \mathbb E^z[V_t]}{\sqrt{2\pi}\,\sigma_v(t)^3} \exp\biggl( -\frac{(\mathbb E^z[V_t])^2}{2\sigma_v(t)^2} \biggr) .
        \label{eq: third sharp formula}
    \end{equation}
    Since $e^{tQ}$ is invertible, $\mathbb E^z[X_t]$ and $\mathbb E^z[V_t]$ take every real value as $z$ ranges over $\mathbb R^2$. By \eqref{eq: third sharp formula x} and \eqref{eq: third sharp formula}, $|\nabla^3_{xvv}(e^{t\mathcal L}f_x)(z)|$ is maximal at $|\mathbb E^z[X_t]| = \sigma_x(t)$ and $|\nabla^3_{vvv}(e^{t\mathcal L}f_v)(z)|$ is maximal at $|\mathbb E^z[V_t]| = \sigma_v(t)$, so
    \begin{equation}
        \norm{\nabla^3_{xvv}(e^{t\mathcal L}f_x)}_\infty = \frac{2\,|(e^{tQ})_{xx}|\,(e^{tQ})_{vx}^2}{\sqrt{2\pi e}\,\sigma_x(t)^2} , \qquad
        \norm{\nabla^3_{vvv}(e^{t\mathcal L}f_v)}_\infty = \frac{2\,|(e^{tQ})_{vv}|^3}{\sqrt{2\pi e}\,\sigma_v(t)^2} .
        \label{eq: third sharp norm}
    \end{equation}
    Since $(e^{tQ})_{xx} = 1 + \mathcal O(t)$, $(e^{tQ})_{vx} = t + \mathcal O(t^2)$ and $(e^{tQ})_{vv} = 1 + \mathcal O(t)$ by the expansion $e^{tQ} = I_2 + tQ + \mathcal O(t^2)$, \eqref{eq: third sharp norm} and \eqref{eq: variance expansions} give \eqref{eq: third sharp}.
\end{proof}

\subsection{Proof of Lemma~\ref{lem: poisson}}
\label{app: poisson}

\begin{proof}[Proof of Lemma~\ref{lem: poisson}]
    Since $\pi$ is invariant under \eqref{eq: uld} and has a finite first moment, \eqref{eq: contraction} gives
    \begin{equation*}
        \bigl| (e^{t\mathcal L}f)(z) - \pi(f) \bigr| \Le C e^{-\lambda t} \norm{\nabla f}_\infty \int_{\mathbb R^{2d}} |z - z'| \,\pi(\D z') , \qquad t \Ge 0 , \quad z \in \mathbb R^{2d} .
    \end{equation*}
    Hence the integral in \eqref{eq: poisson solution} converges at every $z \in \mathbb R^{2d}$. By \eqref{eq: contraction} and \eqref{eq: uniform}, $\phi \in C^1(\mathbb R^{2d})$, $\nabla_v\phi \in C^1(\mathbb R^{2d})$ and \eqref{eq: poisson bounds} holds.
    Note that It\^o's formula \eqref{eq: ito formula} applies to $\phi$ \citep[Chapter~IV, Remarks of Theorem~3.3]{revuz1999continuous}. Taking $\mathbb E^z$ in \eqref{eq: ito formula} and using the Markov property give, for $t > 0$,
    \begin{equation}
        \int_0^t \mathbb E^z (\mathcal L\phi)(Z_s) \,\D s = (e^{t\mathcal L}\phi)(z) - \phi(z) = -\int_0^t \bigl( (e^{s\mathcal L}f)(z) - \pi(f) \bigr) \D s .
        \label{eq: phi increment}
    \end{equation}
    Dividing \eqref{eq: phi increment} by $t$ and letting $t \to 0$ gives \eqref{eq: poisson equation}.
\end{proof}

\subsection{Proof of Poincaré inequalities}
\label{app: poincare}
\label{app: l2}

Let $\bar\pi(x) \propto \exp(-U(x))$ be the distribution of $x$ under $\pi$, so that \eqref{eq: poincare} is the Poincaré inequality of $\bar\pi$.

\begin{proof}[Proof of Lemma~\ref{lem: poincare}]
    The proof utilizes the bounded perturbation trick described in \citep[Proposition~4.2.7]{bakry2014analysis}. To convexify $U$, let $E$ be the radial function on $\mathbb R^d$,
    \begin{equation}
        E(x) := \int_0^{|x|} \max\Bigl\{ 0 ,\ \Bigl( M + \frac m2 \Bigr) \min(r, R) - \frac m2\, \max(r - R, 0) \Bigr\} \D r .
        \label{eq: convexifying function}
    \end{equation}
    The integrand of \eqref{eq: convexifying function} increases with slope $M + \frac m2$ on $[0, R]$, then decreases with slope $-\frac m2$ until it vanishes. The function $E$ is nonnegative and bounded, and $\norm{E}_\infty$ depends only on $m$, $M$ and $R$. The gradient $\nabla E$ is Lipschitz, and $\nabla^2 E = (M + \frac m2) I_d$ in the ball of radius $R$. At the points outside this ball where the Hessian of $E$ exists, its eigenvalue in the direction of $x$ is the slope of the integrand and its other eigenvalues are nonnegative, so $\nabla^2 E \succeq -\frac m2 I_d$. Since $\nabla^2 U \succeq -M I_d$ by \eqref{eq: U bounds} and $\nabla^2 U \succeq m I_d$ outside the ball by \eqref{eq: U convexity}, $\nabla^2(U + E) \succeq \frac m2 I_d$ almost everywhere.

    Let $\pi_E$ be the distribution proportional to $\exp\bigl( -U(x) - E(x) \bigr) \D x$. Then \citep[Corollary~4.8.2]{bakry2014analysis} gives the Poincaré inequality of $\pi_E$, for every $\psi \in C^1(\mathbb R^d)$ with $\norm{\nabla \psi}_\infty <\infty$,
    \begin{equation}
        \pi_E\bigl( (\psi - \pi_E(\psi))^2 \bigr) \Le \frac2m\, \pi_E\bigl( |\nabla \psi|^2 \bigr) .
        \label{eq: convexified poincare}
    \end{equation}
    The distribution $\bar\pi$ has the density $e^{E}/\pi_E(e^{E})$ with respect to $\pi_E$, and $1 \Le e^{E} \Le e^{\norm{E}_\infty}$ by \eqref{eq: convexifying function}, so that
    \begin{equation*}
        e^{-\norm{E}_\infty} \Le \frac{e^{E(x)}}{\pi_E(e^{E})} \Le e^{\norm{E}_\infty} , \qquad x \in \mathbb R^d .
    \end{equation*}
    Hence \citep[Proposition~4.2.7]{bakry2014analysis} applied to \eqref{eq: convexified poincare} gives \eqref{eq: poincare} with $C_P = \frac2m e^{3\norm{E}_\infty}$.
\end{proof}

\begin{proof}[Proof of Lemma~\ref{lem: l2 decay}]
    \citet{cao2023explicit} study the dynamics \eqref{eq: uld} with the generator \eqref{eq: generator}. The hypotheses of \citep[Theorem~1]{cao2023explicit} follow from three facts. First, the Poincaré inequality \eqref{eq: poincare} extends by density to $H^1(\bar\pi)$. Second, by \eqref{eq: U bounds}, the second derivatives of $U$ satisfy
    \begin{equation*}
        \sum_{i,j=1}^d \bigl| \partial_{x_i}\partial_{x_j} U(x) \bigr|^2 \Le \max(M,1)^2 \bigl( d + |\nabla U(x)|^2 \bigr) , \qquad \Delta U(x) \Le \max(M,1)\, d , \qquad x \in \mathbb R^d .
    \end{equation*}
    Third, $U$ grows at least quadratically, so the embedding of $H^1(\bar\pi)$ into $L^2(\bar\pi)$ is compact \citep[Section~1.1]{cao2023explicit}.

    To prove \eqref{eq: l2 decay}, let first $\psi \in C^1(\mathbb R^{2d})$ with $\norm{\psi}_\infty <\infty$, $\norm{\nabla \psi}_\infty <\infty$ and $\pi(\psi) = 0$, so that $\psi, \nabla_v \psi \in L^2(\pi)$. By \eqref{eq: gradient formula} and \eqref{eq: jacobian growth}, $e^{t\mathcal L}\psi \in C^1(\mathbb R^{2d})$ with $\norm{\nabla (e^{t\mathcal L}\psi)}_\infty \Le e^{(1+M+\gamma)t} \norm{\nabla\psi}_\infty$ for every $t \Ge 0$. For $0 < t \Le 1$, Theorem~\ref{thm: hessian} gives $e^{t\mathcal L}\psi \in C^2(\mathbb R^{2d})$, and \eqref{eq: hessian anisotropic} bounds its second derivatives. For $t > 1$, Theorem~\ref{thm: hessian} with $t=1$ gives $e^{t\mathcal L}\psi = e^{\mathcal L}\bigl( e^{(t-1)\mathcal L}\psi \bigr) \in C^2(\mathbb R^{2d})$ with bounded second derivatives. It\^o's formula applied to $e^{t\mathcal L}\psi$ and the semigroup property then give
    \begin{equation*}
        \lim_{s \to 0^+} \frac{e^{(t+s)\mathcal L}\psi - e^{t\mathcal L}\psi}{s} = \mathcal L e^{t\mathcal L}\psi , \qquad t > 0 .
    \end{equation*}
    The invariance of $\pi$ gives $\pi(e^{t\mathcal L}\psi) = 0$, and integration by parts against $\pi$ gives
    \begin{equation*}
        \lim_{s \to 0^+} \frac{\pi\bigl( (e^{(t+s)\mathcal L}\psi)^2 \bigr) - \pi\bigl( (e^{t\mathcal L}\psi)^2 \bigr)}{s} = -2\gamma\, \pi\bigl( |\nabla_v (e^{t\mathcal L}\psi)|^2 \bigr) , \qquad t > 0 .
    \end{equation*}
    The identity $\pi(e^{t\mathcal L}\psi) = 0$ and these two limits are the properties of the solution $e^{t\mathcal L}\psi$ of the backward Kolmogorov equation that the proof of \citep[Theorem~1]{cao2023explicit} uses. Since $\nabla^2 U \succeq -M I_d$ by \eqref{eq: U bounds}, that theorem gives, with a universal constant $c_0$,
    \begin{equation*}
        \pi\bigl( (e^{t\mathcal L}\psi)^2 \bigr)^{\frac12} \Le c_0 \exp\biggl( -\frac{\gamma\, t}{c_0\, C_P \bigl( C_P^{-\frac12} + \sqrt M + \gamma \bigr)^2} \biggr) \pi(\psi^2)^{\frac12} , \qquad t \Ge 0 ,
    \end{equation*}
    which gives \eqref{eq: l2 decay} for such $\psi$.

    For a general $\psi \in L^2(\pi)$ with $\pi(\psi) = 0$, the invariance of $\pi$ under \eqref{eq: uld} shows that $\mathbb E^z\bigl| \psi(Z_t) \bigr|$ is finite for $\pi$-almost every $z$, and Jensen's inequality gives $\pi\bigl( (e^{t\mathcal L}\psi)^2 \bigr) \Le \pi(\psi^2)$. The set $\bigl\{ g \in C^1(\mathbb R^{2d}) : \norm{g}_\infty <\infty ,\ \norm{\nabla g}_\infty <\infty ,\ \pi(g) = 0 \bigr\}$ is dense in the functions of $L^2(\pi)$ with mean zero, so \eqref{eq: l2 decay} extends to $\psi$.
\end{proof}

\section{Proofs of the uniform-in-time error}
\label{app: sgubu}

Throughout this appendix, the letter $C$ denotes a positive constant depending only on $m$, $M$, $R$ and $\gamma$, and it changes from line to line. In particular, $C$ does not depend on the dimension $d$ or on $\Sigma$.

\subsection{Lyapunov function}
\label{app: lyapunov}

Following \citep{schuh2024convergence, lu2025mse}, the proof of Lemma~\ref{lem: sgubu moments} applies a Lyapunov function to the SG-BU integrator. One step of SG-BU with step size $h$ is $\mathscr U_h \circ \mathscr B_h(\xi)$. By \eqref{eq: sgubu}, the states $\mathscr U_{h/2} \mathbf Z_n$ form an SG-BU solution,
\begin{equation}
    \mathscr U_{h/2} \mathbf Z_{n+1} = \bigl( \mathscr U_h \circ \mathscr B_h(\xi_n) \bigr) \mathscr U_{h/2} \mathbf Z_n , \qquad n = 0,1,2,\dots ,
    \label{eq: half step chain}
\end{equation}
where $\mathscr U_h$ is driven by the increments of $(B_t)_{t \Ge 0}$ on $[(n+\frac12)h, (n+\frac32)h]$. Define the Lyapunov function
\begin{equation}
    \mathcal V(x,v) := |\gamma x + v|^2 + |v|^2 , \qquad (x,v) \in \mathbb R^{2d} ,
    \label{eq: lyapunov}
\end{equation}
which satisfies
\begin{equation}
    |x|^2 \Le \frac{2}{\gamma^2}\, \mathcal V(x,v) , \qquad
    |v|^2 \Le \mathcal V(x,v) , \qquad
    \mathcal V(x,v) \Le 2\gamma^2 |x|^2 + 3 |v|^2 .
    \label{eq: lyapunov equivalence}
\end{equation}

Next we show that
\begin{equation}
    x \cdot \nabla U(x) \Ge \frac{3m}8 |x|^2 - C d , \qquad \forall x \in \mathbb R^d .
    \label{eq: dissipativity}
\end{equation}
Let $R_b$ be the radius in \eqref{eq: quoted constants}. For $|x| \Ge R_b$, \eqref{eq: one-sided gradient} at $x' = 0$ and \eqref{eq: U bounds} give
\begin{equation*}
    x \cdot \nabla U(x) \Ge \frac m2 |x|^2 - M \sqrt d\, |x| \Ge \frac{3m}8 |x|^2 - C d .
\end{equation*}
For $|x| < R_b$, \eqref{eq: gradient growth} and $d \Ge 1$ give
\begin{equation*}
    x \cdot \nabla U(x) \Ge -M \bigl( \sqrt d + R_b \bigr) R_b \Ge \frac{3m}8 |x|^2 - C d .
\end{equation*}

\subsection{Proof of Lemma~\ref{lem: sgubu moments}}
\label{app: sgubu moments}

\begin{proof}[Proof of Lemma~\ref{lem: sgubu moments}]
    Fix a state $(x,v) \in \mathbb R^{2d}$ and let $(x',v') := \bigl( \mathscr U_h \circ \mathscr B_h(\xi) \bigr)(x,v)$ be the one-step update of SG-BU, where $\mathscr U_h$ is driven by a Brownian motion $(B_t)_{t \Ge 0}$ independent of $\xi$. We first prove that, for $h \Le h_0$,
    \begin{equation}
        \mathbb E\bigl[ \mathcal V(x',v')^2 \bigr] \Le \Bigl( 1 - \frac{mh}{32\gamma} \Bigr) \mathcal V(x,v)^2 + C d^2 h .
        \label{eq: lyapunov contraction}
    \end{equation}

    To prove \eqref{eq: lyapunov contraction}, define the deterministic vector $\mathfrak a$, the vector $\mathfrak b$ of the stochastic gradient noise \eqref{eq: sg noise} and the Gaussian vector $\zeta$, all in $\mathbb R^{2d}$, by
    \begin{equation*}
    \begin{aligned}
        & \mathfrak a := \Bigl( \gamma x + v - h \nabla U(x) ,\; e^{-\gamma h} \bigl( v - h \nabla U(x) \bigr) \Bigr) , \\
        & \mathfrak b := \bigl( \hat\eta(x;\xi) ,\; e^{-\gamma h}\, \hat\eta(x;\xi) \bigr) , \qquad
        \zeta := \sqrt{2\gamma}\, \biggl( B_h ,\; \int_0^h e^{-\gamma(h-s)} \,\D B_s \biggr) .
    \end{aligned}
    \end{equation*}
    By \eqref{eq: B flow} and \eqref{eq: U flow}, $(\gamma x' + v', v') = \mathfrak a - h \mathfrak b + \zeta$, so \eqref{eq: lyapunov} gives
    \begin{equation}
        \mathcal V(x',v') = |\mathfrak a - h \mathfrak b + \zeta|^2 .
        \label{eq: lyapunov one step}
    \end{equation}

    To bound $|\mathfrak a|^2$, note that
    \begin{equation*}
        |\mathfrak a|^2 = \mathcal V(x,v) - \bigl( 1 - e^{-2\gamma h} \bigr) |v|^2 - 2\gamma h\, x \cdot \nabla U(x) - 2h \bigl( 1 + e^{-2\gamma h} \bigr)\, v \cdot \nabla U(x) + h^2 \bigl( 1 + e^{-2\gamma h} \bigr) |\nabla U(x)|^2 .
    \end{equation*}
    By \eqref{eq: step size}, $\gamma h \Le \frac12 \bigl( 1 - \frac{4M^2}{\gamma^2 m} \bigr)$ and $4hM^2 \Le \frac{\gamma m}8$. With \eqref{eq: gradient growth}, these two bounds give
    \begin{equation*}
    \begin{aligned}
        & 1 - e^{-2\gamma h} \Ge 2\gamma h ( 1 - \gamma h ) \Ge \Bigl( \gamma + \frac{4M^2}{\gamma m} \Bigr) h , \qquad
        h^2 \bigl( 1 + e^{-2\gamma h} \bigr) |\nabla U(x)|^2 \Le \frac{\gamma m}8\, h \bigl( d + |x|^2 \bigr) , \\
        & 2h \bigl( 1 + e^{-2\gamma h} \bigr)\, |v \cdot \nabla U(x)| \Le h \biggl( \frac{\gamma m}2 |x|^2 + \Bigl( \frac{8M^2}{\gamma m} + \frac{m}{2\gamma} \Bigr) |v|^2 + C d \biggr) .
    \end{aligned}
    \end{equation*}
    By \eqref{eq: friction condition}, $\gamma \bigl( 1 - \frac{4M^2}{\gamma^2 m} \bigr) \Ge \frac m\gamma$. Hence  \eqref{eq: dissipativity} and \eqref{eq: lyapunov equivalence} give
    \begin{equation}
        |\mathfrak a|^2 \Le \mathcal V(x,v) - h \Bigl( \frac{\gamma m}8 |x|^2 + \frac{m}{2\gamma} |v|^2 \Bigr) + C d h \Le \Bigl( 1 - \frac{mh}{16\gamma} \Bigr) \mathcal V(x,v) + C d h .
        \label{eq: a bound}
    \end{equation}

    To average over the Brownian motion, note that $\zeta$ is centered and independent of $\xi$. The It\^o isometry gives $\mathbb E |\zeta|^2 = d \bigl( 2\gamma h + 1 - e^{-2\gamma h} \bigr) \Le 4\gamma h d$, and $\mathbb E |\zeta|^4 \Le 3 ( \mathbb E |\zeta|^2 )^2$ since $\zeta$ is Gaussian. Hence
    \begin{equation}
        \mathbb E \bigl[ |\mathfrak a - h\mathfrak b + \zeta|^4 \,\big|\, \xi \bigr]
        \Le |\mathfrak a - h\mathfrak b|^4 + 6\, \mathbb E |\zeta|^2\, |\mathfrak a - h\mathfrak b|^2 + \mathbb E |\zeta|^4
        \Le |\mathfrak a - h\mathfrak b|^4 + C d h\, |\mathfrak a - h\mathfrak b|^2 + C d^2 h^2 .
        \label{eq: gaussian fourth}
    \end{equation}
    To average over $\xi$, note that $\mathbb E_\xi[\mathfrak b] = 0$ by Assumption~\ref{ass: sg} and $|\mathfrak b|^2 \Le 2 |\hat\eta(x;\xi)|^2$. By \eqref{eq: sg variance} and H\"older's inequality, $\mathbb E_\xi |\mathfrak b|^k \Le ( 2\Sigma d )^{\frac k2}$ for $k = 2, 3, 4$. By \eqref{eq: step size}, $h \Sigma \Le \frac{\gamma m}{576}$. Hence
    \begin{equation}
        \mathbb E_\xi |\mathfrak a - h\mathfrak b|^4 \Le |\mathfrak a|^4 + 6h^2 |\mathfrak a|^2\, \mathbb E_\xi |\mathfrak b|^2 + 4h^3 |\mathfrak a|\, \mathbb E_\xi |\mathfrak b|^3 + h^4\, \mathbb E_\xi |\mathfrak b|^4
        \Le |\mathfrak a|^4 + C d h\, |\mathfrak a|^2 + C d^2 h^2 .
        \label{eq: xi moments}
    \end{equation}

    Taking the expectation of \eqref{eq: gaussian fourth} with respect to $\xi$ and inserting \eqref{eq: lyapunov one step}, \eqref{eq: xi moments}, \eqref{eq: a bound} and $h \Le \frac1{2\gamma}$ gives
    \begin{equation*}
        \mathbb E\bigl[ \mathcal V(x',v')^2 \bigr] \Le |\mathfrak a|^4 + C d h\, |\mathfrak a|^2 + C d^2 h^2
        \Le \Bigl( 1 - \frac{mh}{16\gamma} \Bigr)^2 \mathcal V(x,v)^2 + C d h\, \mathcal V(x,v) + C d^2 h .
    \end{equation*}
    By $h \Le \frac1{2\gamma}$ and \eqref{eq: friction condition}, $\frac{mh}{16\gamma} \Le \frac{m}{32\gamma^2} \Le 1$, so with $C d h\, \mathcal V(x,v) \Le \frac{mh}{32\gamma} \mathcal V(x,v)^2 + C d^2 h$ we prove \eqref{eq: lyapunov contraction}.

    To iterate \eqref{eq: lyapunov contraction} along the SG-BU solution \eqref{eq: half step chain}, note that $\mathscr U_{h/2} \mathbf Z_n$ is independent of $\xi_n$ and of the increments of $(B_t)_{t \Ge 0}$ after the time $(n+\frac12)h$. Hence \eqref{eq: lyapunov contraction} at $(x,v) = \mathscr U_{h/2} \mathbf Z_n$ gives
    \begin{equation*}
        \mathbb E\bigl[ \mathcal V( \mathscr U_{h/2} \mathbf Z_{n+1} )^2 \bigr] \Le \Bigl( 1 - \frac{mh}{32\gamma} \Bigr)\, \mathbb E\bigl[ \mathcal V( \mathscr U_{h/2} \mathbf Z_n )^2 \bigr] + C d^2 h , \qquad n \Ge 0 .
    \end{equation*}
    Since $\mathbf X_0 = 0$, $\mathbf V_0 \sim \mathcal N(0, I_d)$ and $h \Le \frac1{2\gamma}$, \eqref{eq: U flow} and \eqref{eq: lyapunov equivalence} give $\mathbb E\bigl[ \mathcal V( \mathscr U_{h/2} \mathbf Z_0 )^2 \bigr] \Le C\, \mathbb E |\mathscr U_{h/2} \mathbf Z_0|^4 \Le C d^2$. Induction on $n$ gives
    \begin{equation}
        \sup_{n \Ge 0} \mathbb E\bigl[ \mathcal V( \mathscr U_{h/2} \mathbf Z_n )^2 \bigr] \Le C d^2 .
        \label{eq: half step moments}
    \end{equation}

    Note that \eqref{eq: gradient growth}, \eqref{eq: sg variance}, $h \Le \frac1{2\gamma}$ and $h \Sigma \Le \frac{\gamma m}{576}$ give $h^4\, \mathbb E_\xi |\hat g(x;\xi)|^4 \Le C \bigl( d + |x|^2 \bigr)^2$ for every $x \in \mathbb R^d$. For $(x,v) \in \mathbb R^{2d}$ and a map $\mathscr U_{h/2}$ driven by a Brownian motion independent of $\xi$, \eqref{eq: B flow}, \eqref{eq: U flow} and \eqref{eq: lyapunov equivalence} then give
    \begin{equation*}
        \mathbb E \bigl| \bigl( \mathscr U_{h/2} \circ \mathscr B_h(\xi) \bigr)(x,v) \bigr|^4 \Le C \bigl( |x|^4 + |v|^4 + h^4\, \mathbb E_\xi |\hat g(x;\xi)|^4 + d^2 \bigr) \Le C \bigl( \mathcal V(x,v)^2 + d^2 \bigr) .
    \end{equation*}
    By \eqref{eq: sgubu}, this bound at $(x,v) = \mathscr U_{h/2} \mathbf Z_n$ and \eqref{eq: half step moments} give $\mathbb E |\mathbf Z_{n+1}|^4 \Le C d^2$ for every $n \Ge 0$. Since $\mathbb E |\mathbf Z_0|^4 = \mathbb E |\mathbf V_0|^4 \Le 3d^2$, this proves \eqref{eq: sgubu moments}.
\end{proof}

\subsection{Local error of the full gradient step}
\label{app: local error}

For $n \Ge 0$ and $s \in [0,h]$, let $\mathbf Z_n(s) = (\mathbf X_n(s), \mathbf V_n(s))$ be the solution to \eqref{eq: uld} at the time $nh + s$ started from $\mathbf Z_n$ at the time $nh$ and driven by $(B_t)_{t \Ge 0}$. The local error of the step \eqref{eq: full gradient step} with the gradient $\nabla U$ is
\begin{equation}
    \Delta\mathbf Z_n = (\Delta\mathbf X_n, \Delta\mathbf V_n) := \bar{\mathbf Z}_{n+1} - \mathbf Z_n(h) \in \mathbb R^{2d} .
    \label{eq: full gradient error}
\end{equation}
As in \citep{sanz2021wasserstein, lu2025mse}, Lemma~\ref{lem: local error} splits $\Delta\mathbf Z_n$ into a \emph{mean-zero} term $\mathcal M(\Delta\mathbf Z_n)$, of zero conditional mean given $\mathbf Z_n$, and a \emph{remainder} $\mathcal R(\Delta\mathbf Z_n)$.

\begin{lemma}
\label{lem: local error}
    Under Assumptions~\ref{ass: U}--\ref{ass: sg}, let $h \Le h_0$ with $h_0$ in \eqref{eq: step size}. For every $n \Ge 0$,
    \begin{equation}
        \bigl( \mathbb E |\mathbf X_n(h) - \mathbf X_n|^2 \bigr)^{\frac12} \Le C d^{\frac12} h , \qquad
        \bigl( \mathbb E |\mathbf V_n(h) - \mathbf V_n|^2 \bigr)^{\frac12} \Le C d^{\frac12} h^{\frac12} ,
        \label{eq: one step}
    \end{equation}
    and $\Delta\mathbf Z_n$ can be decomposed as $\Delta\mathbf Z_n = \mathcal M(\Delta\mathbf Z_n) + \mathcal R(\Delta\mathbf Z_n)$ with
    \begin{equation}
        \mathbb E\bigl[ \mathcal M(\Delta\mathbf Z_n) \,\big|\, \mathbf Z_n \bigr] = 0 , \qquad
        \bigl( \mathbb E |\mathcal M(\Delta\mathbf Z_n)|^2 \bigr)^{\frac12} \Le C d^{\frac12} h^{\frac52} , \qquad
        \bigl( \mathbb E |\mathcal R(\Delta\mathbf Z_n)|^2 \bigr)^{\frac12} \Le C d h^3 .
        \label{eq: local error}
    \end{equation}
\end{lemma}

\begin{proof}[Proof of Lemma~\ref{lem: local error}]
    Note that \eqref{eq: uld}, \eqref{eq: gradient growth} and Gr\"onwall's inequality on $[0,h]$, with $h \Le \frac1{2\gamma}$, give
    \begin{equation}
        \sup_{0 \Le s \Le h} |\mathbf Z_n(s) - \mathbf Z_n| \Le C \Bigl( \bigl( \sqrt d + |\mathbf Z_n| \bigr) h + \sup_{0 \Le s \Le h} |B_{nh+s} - B_{nh}| \Bigr) .
        \label{eq: displacement sup}
    \end{equation}
    Since $\mathbb E |\mathbf Z_n|^4 \Le C d^2$ by \eqref{eq: sgubu moments} and $\mathbb E \bigl[ \sup_{0 \Le s \Le h} |B_{nh+s} - B_{nh}|^4 \bigr] \Le C d^2 h^2$, \eqref{eq: displacement sup} gives
    \begin{equation}
        \biggl( \mathbb E\Bigl[ \sup_{0 \Le s \Le h} |\mathbf Z_n(s) - \mathbf Z_n|^4 \Bigr] \biggr)^{\frac14} \Le C d^{\frac12} h^{\frac12} , \qquad
        \biggl( \mathbb E \Bigl( \sqrt d + \sup_{0 \Le s \Le h} |\mathbf Z_n(s)| \Bigr)^4 \biggr)^{\frac14} \Le C d^{\frac12} .
        \label{eq: displacement moments}
    \end{equation}
    With $\mathbf X_n(h) - \mathbf X_n = \int_0^h \mathbf V_n(s) \,\D s$, \eqref{eq: displacement moments} gives \eqref{eq: one step}.
    For $s \in [0,h]$, define the weights
    \begin{equation}
        \beta_x(s) := \frac{1-e^{-\gamma(h-s)}}{\gamma} , \qquad
        \beta_v(s) := e^{-\gamma(h-s)} , \qquad
        \beta(s) := (\beta_x(s), \beta_v(s)) ,
        \label{eq: weights}
    \end{equation}
    which satisfy $0 \Le \beta_x(s) \Le h - s$ and $|\beta_x| + |\beta_x'| + |\beta_x''| + |\beta_v| + |\beta_v'| + |\beta_v''| \Le C$ on $[0,h]$. The product $\beta(s) \nabla U(x)$ denotes $\bigl( \beta_x(s) \nabla U(x), \beta_v(s) \nabla U(x) \bigr) \in \mathbb R^{2d}$, and similarly for other vectors in $\mathbb R^d$.

    To split $\Delta\mathbf Z_n$, note that \eqref{eq: uld} is the dynamics $\mathscr U$ with the additional force $-\nabla U(X_t)$ in the velocity. Hence \eqref{eq: U flow}, \eqref{eq: weights} and \eqref{eq: full gradient step} give
    \begin{equation}
        \mathbf Z_n(h) = \mathscr U_h \mathbf Z_n - \int_0^h \beta(s) \nabla U(\mathbf X_n(s)) \,\D s , \qquad
        \bar{\mathbf Z}_{n+1} = \mathscr U_h \mathbf Z_n - h\, \beta\bigl( \tfrac h2 \bigr) \nabla U(\mathbf Y_n) ,
        \label{eq: step states}
    \end{equation}
    where $\mathscr U_h$ is driven by the increments of $(B_t)_{t \Ge 0}$ on $[nh, (n+1)h]$. Inserting \eqref{eq: step states} into \eqref{eq: full gradient error} gives
    \begin{equation}
    	\Delta\mathbf X_n = 
        \underbrace{\int_0^h \beta_x(s)\, \nabla U(\mathbf X_n(s)) \,\D s - h\, \beta_x\bigl( \tfrac h2 \bigr) \nabla U\bigl( \mathbf X_n(\tfrac h2) \bigr)}_{\mathcal T_x} + h\, \beta_x\bigl( \tfrac h2 \bigr) \Bigl( \nabla U\bigl( \mathbf X_n(\tfrac h2) \bigr) - \nabla U(\mathbf Y_n) \Bigr) ,
        \label{eq: midpoint form}
    \end{equation}
    and \eqref{eq: midpoint form} holds with $\Delta\mathbf V_n$, $\beta_v$ and $\mathcal T_v$ in place of $\Delta\mathbf X_n$, $\beta_x$ and $\mathcal T_x$.
	
	In \eqref{eq: midpoint form}, $\mathcal T_x$ is the error of the midpoint rule for the function $s \mapsto \beta_x(s) \nabla U(\mathbf X_n(s))$, whose derivative is $\beta_x'(s) \nabla U(\mathbf X_n(s)) + \beta_x(s) \nabla^2 U(\mathbf X_n(s))\, \mathbf V_n(s)$. Hence It\^o's formula with \eqref{eq: uld} gives
    \begin{align}
        \mathcal T_x & = \int_{\frac h2}^h \int_{\frac h2}^s \int_{h-s'}^{s'} \Bigl( \beta_x''(r) \nabla U(\mathbf X_n(r)) + \bigl( 2\beta_x'(r) - \gamma \beta_x(r) \bigr) \nabla^2 U(\mathbf X_n(r))\, \mathbf V_n(r) \Bigr) \D r \,\D s' \,\D s \notag \\
        & + \int_{\frac h2}^h \int_{\frac h2}^s \int_{h-s'}^{s'} \beta_x(r) \Bigl( \nabla^3 U(\mathbf X_n(r)) [\mathbf V_n(r), \mathbf V_n(r)] - \nabla^2 U(\mathbf X_n(r)) \nabla U(\mathbf X_n(r)) \Bigr) \D r \,\D s' \,\D s \notag \\
        & + \sqrt{2\gamma} \int_{\frac h2}^h \int_{\frac h2}^s \int_{h-s'}^{s'} \beta_x(r) \nabla^2 U(\mathbf X_n(r)) \,\D B_{nh+r} \,\D s' \,\D s ,
        \label{eq: midpoint expansion}
    \end{align}
    and the same identity holds with $\beta_v$ and $\mathcal T_v$ in place of $\beta_x$ and $\mathcal T_x$. Define the term $\mathcal M(\Delta\mathbf Z_n)$ by
    \begin{equation}
        \mathcal M(\Delta\mathbf Z_n) := \sqrt{2\gamma} \int_{\frac h2}^h \int_{\frac h2}^s \int_{h-s'}^{s'} \beta(r) \nabla^2 U(\mathbf X_n(r)) \,\D B_{nh+r} \,\D s' \,\D s ,
        \label{eq: mean-zero term}
    \end{equation}
    and the remainder by $\mathcal R(\Delta\mathbf Z_n) := \Delta\mathbf Z_n - \mathcal M(\Delta\mathbf Z_n)$.

    The state $\mathbf Z_n$ is independent of the increments of $(B_t)_{t \Ge 0}$ after the time $nh$, and the integrand of \eqref{eq: mean-zero term} is adapted and bounded, so $\mathbb E[ \mathcal M(\Delta\mathbf Z_n) \,|\, \mathbf Z_n ] = 0$. By \eqref{eq: U bounds}, $\sum_{i,j=1}^d | \partial_{x_i}\partial_{x_j} U(x) |^2 \Le d M^2$, so the It\^o isometry gives $\mathbb E |\mathcal M(\Delta\mathbf Z_n)|^2 \Le C d h^5$. Since $\mathbf Y_n$ is the position of $\mathscr U_{h/2} \mathbf Z_n$, \eqref{eq: U flow}, $0 \Le \beta_x(s) \Le h - s$ and \eqref{eq: gradient growth} give
    \begin{equation*}
        \Bigl| \mathbf X_n\bigl( \tfrac h2 \bigr) - \mathbf Y_n \Bigr| = \biggl| \int_0^{\frac h2} \beta_x\bigl( s + \tfrac h2 \bigr) \nabla U(\mathbf X_n(s)) \,\D s \biggr| \Le C h^2 \Bigl( \sqrt d + \sup_{0 \Le s \Le h} |\mathbf Z_n(s)| \Bigr) .
    \end{equation*}
    Hence \eqref{eq: midpoint form}, \eqref{eq: midpoint expansion}, \eqref{eq: U bounds} and \eqref{eq: gradient growth} give
    \begin{equation*}
        |\mathcal R(\Delta\mathbf Z_n)| \Le C h^3 \Bigl( \sqrt d + \sup_{0 \Le s \Le h} |\mathbf Z_n(s)| \Bigr)^2 .
    \end{equation*}
    With \eqref{eq: displacement moments}, this bound gives $\bigl( \mathbb E |\mathcal R(\Delta\mathbf Z_n)|^2 \bigr)^{\frac12} \Le C d h^3$, which completes the proof of \eqref{eq: local error}.
\end{proof}

\subsection{Proof of Theorem~\ref{thm: sgubu}}
\label{app: sgubu error}

Lemma~\ref{lem: half step} below bounds the derivatives of the mean of a test function over the map $\mathscr U_{h/2}$.

\begin{lemma}
\label{lem: half step}
    Let $h \Le 1$. Then for every $f \in C^1(\mathbb R^{2d})$ with $\norm{\nabla f}_\infty < \infty$, the function
    \begin{equation}
        f_{h}(z) := \mathbb E\bigl[ f( \mathscr U_{h/2} z ) \bigr] , \qquad z \in \mathbb R^{2d} ,
        \label{eq: half step}
    \end{equation}
    belongs to $C^2(\mathbb R^{2d})$, $\nabla^2_{vv}f_{h}$ is continuously differentiable on $\mathbb R^{2d}$, and
    \begin{equation}
        \norm{\nabla f_{h}}_\infty \Le C \norm{\nabla f}_\infty , \quad
        \norm{\nabla^2_{vv}f_{h}}_\infty \Le C h^{-\frac12} \norm{\nabla f}_\infty , \quad
        \norm{\nabla^3_{xvv}f_{h}}_\infty + \norm{\nabla^3_{vvv}f_{h}}_\infty \Le C h^{-1} \norm{\nabla f}_\infty .
        \label{eq: half step bounds}
    \end{equation}
\end{lemma}

\begin{proof}[Proof of Lemma~\ref{lem: half step}]
    The dynamics $\mathscr U$ is \eqref{eq: uld} with $U = 0$, which satisfies Assumption~\ref{ass: U}, and $f_{h}$ is $e^{\frac h2\mathcal L}f$ for this $U=0$. Hence Theorems~\ref{thm: hessian} and~\ref{thm: third} at $t = \frac h2$ give $f_{h} \in C^2(\mathbb R^{2d})$ with $\nabla^2_{vv}f_{h}$ continuously differentiable on $\mathbb R^{2d}$, and \eqref{eq: gradient bound}, \eqref{eq: hessian anisotropic} and \eqref{eq: third} give \eqref{eq: half step bounds}.
\end{proof}

\begin{proof}[Proof of Theorem~\ref{thm: sgubu}]
    Assume $n \Ge 1$, and put $t_k := (n-k-1)h$ for $0 \Le k \Le n-1$.
    Since $\mathbf Z_k(h)$ is the solution to \eqref{eq: uld} in time $h$,
    \begin{equation*}
        \mathbb E\bigl[ (e^{t_k\mathcal L}f)(\mathbf Z_k(h)) \bigr] = \mathbb E\bigl[ (e^{(t_k+h)\mathcal L}f)(\mathbf Z_k) \bigr] , \qquad 0 \Le k \Le n-1 .
    \end{equation*}
    With $Z_0 = \mathbf Z_0$, this identity gives the telescoping sum
    \begin{equation}
        \mathbb E [f(\mathbf Z_n)] - \mathbb E [f(Z_{nh})]
        = \sum_{k=0}^{n-1} \mathbb E\Bigl[ (e^{t_k\mathcal L}f)(\mathbf Z_{k+1}) - (e^{t_k\mathcal L}f)(\mathbf Z_k(h)) \Bigr] .
        \label{eq: sgubu telescoping}
    \end{equation}
    In each term of \eqref{eq: sgubu telescoping}, we insert the state $\bar{\mathbf Z}_{k+1}$ of \eqref{eq: full gradient step} and aim to prove, for the stochastic gradient error $\mathbf Z_{k+1} - \bar{\mathbf Z}_{k+1}$ and for the local error $\Delta\mathbf Z_k = \bar{\mathbf Z}_{k+1} - \mathbf Z_k(h)$ of \eqref{eq: full gradient error},
    \begin{align}
        & \biggl| \sum_{k=0}^{n-1} \mathbb E\Bigl[ (e^{t_k\mathcal L}f)(\mathbf Z_{k+1}) - (e^{t_k\mathcal L}f)(\bar{\mathbf Z}_{k+1}) \Bigr] \biggr| \Le C \norm{\nabla f}_\infty\, \Sigma d h ,
        \label{eq: sgubu sg part} \\
        & \biggl| \sum_{k=0}^{n-1} \mathbb E\Bigl[ (e^{t_k\mathcal L}f)(\bar{\mathbf Z}_{k+1}) - (e^{t_k\mathcal L}f)(\mathbf Z_k(h)) \Bigr] \biggr| \Le C \norm{\nabla f}_\infty\, d h^2 .
        \label{eq: sgubu local part}
    \end{align}
    With \eqref{eq: sgubu telescoping}, the bounds \eqref{eq: sgubu sg part} and \eqref{eq: sgubu local part} give \eqref{eq: sgubu error}.\\[6pt]
    \textbf{Proof of \eqref{eq: sgubu sg part}.} Note that $\xi_k$ is independent of $\mathbf Y_k$, so \eqref{eq: sg variance} gives
    \begin{equation}
        \mathbb E |\hat\eta(\mathbf Y_k;\xi_k)|^2 \Le \Sigma d , \qquad k \Ge 0 .
        \label{eq: sgubu inputs}
    \end{equation}
    Fix $0 \Le k \Le n-2$. Using $e^{t_k\mathcal L}f \in C^2(\mathbb R^{2d})$ by Corollary~\ref{cor: uniform}, we write
    \begin{equation}
        (e^{t_k\mathcal L}f)(\mathbf Z_{k+1}) - (e^{t_k\mathcal L}f)(\bar{\mathbf Z}_{k+1}) = \mathcal T_1 + \mathcal T_2 ,
        \label{eq: sgubu taylor}
    \end{equation}
    where the terms $\mathcal T_1$ and  $\mathcal T_2$ are
    \begin{align*}
        \mathcal T_1 & := \nabla (e^{t_k\mathcal L}f)(\bar{\mathbf Z}_{k+1}) \cdot ( \mathbf Z_{k+1} - \bar{\mathbf Z}_{k+1} ) , \\
        \mathcal T_2 & := \int_0^1 (1-\theta)\, \nabla^2 (e^{t_k\mathcal L}f)\bigl( \bar{\mathbf Z}_{k+1} + \theta ( \mathbf Z_{k+1} - \bar{\mathbf Z}_{k+1} ) \bigr)\bigl[ \mathbf Z_{k+1} - \bar{\mathbf Z}_{k+1} , \mathbf Z_{k+1} - \bar{\mathbf Z}_{k+1} \bigr] \D\theta .
    \end{align*}
    By \eqref{eq: sg error}, $\mathbf Z_{k+1} - \bar{\mathbf Z}_{k+1} = -h \bigl( \frac{1-e^{-\frac{\gamma h}2}}{\gamma} , e^{-\frac{\gamma h}2} \bigr) \hat\eta(\mathbf Y_k;\xi_k)$, which has mean zero given $\mathbf Z_k$ and $(B_t)_{t \Ge 0}$. Hence $\mathcal T_1$ of \eqref{eq: sgubu taylor} has mean zero. For $\mathcal T_2$, \eqref{eq: uniform} and \eqref{eq: uniform xx} give, for every $z \in \mathbb R^{2d}$,
    \begin{align}
        & \quad \Bigl| \nabla^2 (e^{t_k\mathcal L}f)(z)\bigl[ \mathbf Z_{k+1} - \bar{\mathbf Z}_{k+1} , \mathbf Z_{k+1} - \bar{\mathbf Z}_{k+1} \bigr] - h^2\, \nabla^2_{vv} (e^{t_k\mathcal L}f)(z)\bigl[ \hat\eta(\mathbf Y_k;\xi_k), \hat\eta(\mathbf Y_k;\xi_k) \bigr] \Bigr| \notag \\
        & \Le h^2 \Bigl( \frac{h^2}4 \norm{\nabla^2_{xx} (e^{t_k\mathcal L}f)}_\infty + h \norm{\nabla^2_{xv} (e^{t_k\mathcal L}f)}_\infty + \gamma h \norm{\nabla^2_{vv} (e^{t_k\mathcal L}f)}_\infty \Bigr) |\hat\eta(\mathbf Y_k;\xi_k)|^2 \notag \\
        & \Le C \norm{\nabla f}_\infty e^{-\lambda t_k} \Bigl( h^3 \min(1,t_k)^{-\frac12} + h^4 \min(1,t_k)^{-\frac32} \Bigr) |\hat\eta(\mathbf Y_k;\xi_k)|^2 .
        \label{eq: hessian blocks}
    \end{align}
    For $0 \Le k \Le n-2$, $t_k \Ge h$. Since $h \Le \min(1,t_k)$, \eqref{eq: hessian blocks} and \eqref{eq: uniform} give
    \begin{equation}
        |\mathcal T_2| \Le C \norm{\nabla f}_\infty e^{-\lambda t_k} \min(1,t_k)^{-\frac12} h^2\, |\hat\eta(\mathbf Y_k;\xi_k)|^2 .
        \label{eq: sgubu second order}
    \end{equation}
    With $\lambda$ as in Lemma~\ref{lem: contraction} and $h \Le \frac1{2\gamma}$ by \eqref{eq: step size},
    \begin{equation}
        \sum_{j=1}^\infty h \min(1,jh)^{-\frac12} e^{-\lambda jh} \Le C , \qquad
        \sum_{j=1}^\infty h \min(1,jh)^{-\frac32} e^{-\lambda jh} \Le C h^{-\frac12} .
        \label{eq: time sums}
    \end{equation}
    Hence \eqref{eq: sgubu inputs}--\eqref{eq: time sums} give
    \begin{equation*}
        \biggl| \sum_{k=0}^{n-2} \mathbb E\Bigl[ (e^{t_k\mathcal L}f)(\mathbf Z_{k+1}) - (e^{t_k\mathcal L}f)(\bar{\mathbf Z}_{k+1}) \Bigr] \biggr| \Le C \norm{\nabla f}_\infty\, \Sigma d h .
    \end{equation*}
    For $k = n-1$, \eqref{eq: sgubu}, \eqref{eq: full gradient step} and \eqref{eq: B flow} give
    \begin{equation*}
        \mathbf Z_n = \mathscr U_{h/2} \Bigl( \mathscr B_h \mathscr U_{h/2} \mathbf Z_{n-1} - h \bigl( 0, \hat\eta(\mathbf Y_{n-1};\xi_{n-1}) \bigr) \Bigr) , \qquad
        \bar{\mathbf Z}_n = \mathscr U_{h/2} \bigl( \mathscr B_h \mathscr U_{h/2} \mathbf Z_{n-1} \bigr) .
    \end{equation*}
    Here the map $\mathscr U_{h/2}$ applied last is driven by the increments of $(B_t)_{t \Ge 0}$ on $[(n-\frac12)h, nh]$, which are independent of $\mathscr U_{h/2} \mathbf Z_{n-1}$ and $\xi_{n-1}$. Hence, with $f_{h}$ in \eqref{eq: half step},
    \begin{equation}
        \mathbb E\bigl[ f(\mathbf Z_n) - f(\bar{\mathbf Z}_n) \bigr] = \mathbb E\Bigl[ f_{h}\Bigl( \mathscr B_h \mathscr U_{h/2} \mathbf Z_{n-1} - h \bigl( 0, \hat\eta(\mathbf Y_{n-1};\xi_{n-1}) \bigr) \Bigr) - f_{h}\bigl( \mathscr B_h \mathscr U_{h/2} \mathbf Z_{n-1} \bigr) \Bigr] .
        \label{eq: last step}
    \end{equation}
    Since $\xi_{n-1}$ is independent of $\mathscr U_{h/2} \mathbf Z_{n-1}$, the term of first order in Taylor's formula for $f_{h}$ at $\mathscr B_h \mathscr U_{h/2} \mathbf Z_{n-1}$ has mean zero, so \eqref{eq: half step bounds} and \eqref{eq: sgubu inputs} give
    \begin{equation*}
        \bigl| \mathbb E\bigl[ f(\mathbf Z_n) - f(\bar{\mathbf Z}_n) \bigr] \bigr| \Le \frac{h^2}2 \norm{\nabla^2_{vv}f_{h}}_\infty\, \mathbb E |\hat\eta(\mathbf Y_{n-1};\xi_{n-1})|^2 \Le C \norm{\nabla f}_\infty\, \Sigma d h^{\frac32} ,
    \end{equation*}
    which completes the proof of \eqref{eq: sgubu sg part}.\\[6pt]
    \textbf{Proof of \eqref{eq: sgubu local part}.} Fix $0 \Le k \Le n-2$. By \eqref{eq: uniform} and \eqref{eq: uniform xx}, for every $z = (x,v)$ and $z' = (x',v')$ in $\mathbb R^{2d}$,
    \begin{equation}
        \bigl| \nabla (e^{t_k\mathcal L}f)(z) - \nabla (e^{t_k\mathcal L}f)(z') \bigr| \Le C \norm{\nabla f}_\infty e^{-\lambda t_k} \Bigl( \min(1,t_k)^{-\frac32} |x - x'| + \min(1,t_k)^{-\frac12} |v - v'| \Bigr) .
        \label{eq: gradient lipschitz}
    \end{equation}
    Since $\bar{\mathbf Z}_{k+1} = \mathbf Z_k(h) + \Delta\mathbf Z_k$, \eqref{eq: gradient lipschitz} at $z = \mathbf Z_k(h) + \theta\, \Delta\mathbf Z_k$, $\theta \in [0,1]$, and $z' = \mathbf Z_k$ gives
    \begin{equation}
    \begin{aligned}
        & \Bigl| (e^{t_k\mathcal L}f)(\bar{\mathbf Z}_{k+1}) - (e^{t_k\mathcal L}f)(\mathbf Z_k(h)) - \nabla (e^{t_k\mathcal L}f)(\mathbf Z_k) \cdot \Delta\mathbf Z_k \Bigr| \\
        & \quad \Le C \norm{\nabla f}_\infty e^{-\lambda t_k} \Bigl( \min(1,t_k)^{-\frac32} \bigl( |\mathbf X_k(h) - \mathbf X_k| + |\Delta\mathbf Z_k| \bigr) + \min(1,t_k)^{-\frac12} |\mathbf V_k(h) - \mathbf V_k| \Bigr) |\Delta\mathbf Z_k| .
    \end{aligned}
        \label{eq: sgubu remainder}
    \end{equation}
    Since $\mathbb E[ \mathcal M(\Delta\mathbf Z_k) \,|\, \mathbf Z_k ] = 0$, \eqref{eq: contraction} and \eqref{eq: local error} give
    \begin{equation}
        \Bigl| \mathbb E\bigl[ \nabla (e^{t_k\mathcal L}f)(\mathbf Z_k) \cdot \Delta\mathbf Z_k \bigr] \Bigr| = \Bigl| \mathbb E\bigl[ \nabla (e^{t_k\mathcal L}f)(\mathbf Z_k) \cdot \mathcal R(\Delta\mathbf Z_k) \bigr] \Bigr| \Le C \norm{\nabla f}_\infty e^{-\lambda t_k}\, d h^3 .
        \label{eq: sgubu mean}
    \end{equation}
    Since $h \Le \frac1d$, \eqref{eq: local error} gives $\bigl( \mathbb E |\Delta\mathbf Z_k|^2 \bigr)^{\frac12} \Le C d^{\frac12} h^{\frac52}$, and with \eqref{eq: one step},
    \begin{equation}
        \mathbb E\bigl[ |\mathbf X_k(h) - \mathbf X_k|\, |\Delta\mathbf Z_k| \bigr] + \mathbb E |\Delta\mathbf Z_k|^2 \Le C d h^{\frac72} , \qquad
        \mathbb E\bigl[ |\mathbf V_k(h) - \mathbf V_k|\, |\Delta\mathbf Z_k| \bigr] \Le C d h^3 .
        \label{eq: sgubu products}
    \end{equation}
    By \eqref{eq: sgubu remainder}--\eqref{eq: sgubu products},
    \begin{equation*}
        \Bigl| \mathbb E\bigl[ (e^{t_k\mathcal L}f)(\bar{\mathbf Z}_{k+1}) - (e^{t_k\mathcal L}f)(\mathbf Z_k(h)) \bigr] \Bigr|
        \Le C \norm{\nabla f}_\infty e^{-\lambda t_k}\, d \Bigl( \min(1,t_k)^{-\frac32} h^{\frac72} + \min(1,t_k)^{-\frac12} h^3 \Bigr) .
    \end{equation*}
    Summing this bound over $0 \Le k \Le n-2$ with \eqref{eq: time sums} gives
    \begin{equation*}
        \biggl| \sum_{k=0}^{n-2} \mathbb E\Bigl[ (e^{t_k\mathcal L}f)(\bar{\mathbf Z}_{k+1}) - (e^{t_k\mathcal L}f)(\mathbf Z_k(h)) \Bigr] \biggr| \Le C \norm{\nabla f}_\infty\, d h^2 .
    \end{equation*}
    For $k = n-1$, \eqref{eq: local error} gives
    \begin{equation*}
        \bigl| \mathbb E\bigl[ f(\bar{\mathbf Z}_n) - f(\mathbf Z_{n-1}(h)) \bigr] \bigr| \Le \norm{\nabla f}_\infty\, \mathbb E |\Delta\mathbf Z_{n-1}| \Le C \norm{\nabla f}_\infty \bigl( d^{\frac12} h^{\frac52} + d h^3 \bigr) \Le C \norm{\nabla f}_\infty\, d h^2 ,
    \end{equation*}
    which completes the proof of \eqref{eq: sgubu local part} and of \eqref{eq: sgubu error}.
\end{proof}

\subsection{Proof of Lemma~\ref{lem: invariant}}
\label{app: invariant}

\begin{proof}[Proof of Lemma~\ref{lem: invariant}]
    For $z \in \mathbb R^{2d}$, let $\mathbb E^z$ also denote the expectation for SG-UBU with the initial state $\mathbf Z_0 = z$. We first prove (i) for every bounded measurable function $\psi$ on $\mathbb R^{2d}$,
    \begin{equation}
        z \mapsto \mathbb E^z\bigl[ \psi(\mathbf Z_1) \bigr] \text{ is continuous on } \mathbb R^{2d} ;
        \label{eq: feller}
    \end{equation}
    (ii) for the initial state $\mathbf Z_0 = z$, the law of $\mathbf Z_1$ has a density $\mathfrak p(z,\cdot)$ in the Lebesgue measure, with
    \begin{equation}
        \mathfrak p(z,z') > 0 , \qquad \forall z, z' \in \mathbb R^{2d} .
        \label{eq: positive density}
    \end{equation}

    By \eqref{eq: U flow}, $\mathscr U_{h/2} z$ is a Gaussian vector whose mean is linear in $z$ and whose covariance matrix is positive definite and does not depend on $z$. Hence the law of $\mathscr U_{h/2} z$ has a strictly positive density in the Lebesgue measure and depends continuously on $z$ in the total variation distance.

    To prove \eqref{eq: feller}, note that in the step \eqref{eq: sgubu} from $\mathbf Z_0 = z$, the map $\mathscr U_{h/2}$ applied first, the random variable $\xi_0$ and the map $\mathscr U_{h/2}$ applied last are independent. Hence $\mathbb E^z[\psi(\mathbf Z_1)]$ is the mean under the law of $\mathscr U_{h/2} z$ of the function
    \begin{equation*}
        z' \mapsto \mathbb E\bigl[ \psi\bigl( \mathscr U_{h/2} \mathscr B_h(\xi_0) z' \bigr) \bigr] , \qquad z' \in \mathbb R^{2d} ,
    \end{equation*}
    where the expectation is over $\xi_0$ and the map $\mathscr U_{h/2}$ applied last. This function is bounded and does not depend on $z$, so the continuity of the law of $\mathscr U_{h/2} z$ in the total variation distance gives \eqref{eq: feller}.

    To prove \eqref{eq: positive density}, note that the state $\mathscr B_h(\xi_0) \mathscr U_{h/2} z$ is independent of the map $\mathscr U_{h/2}$ applied last, which maps this state to $\mathbf Z_1$. Hence, given this state, the law of $\mathbf Z_1$ is Gaussian with a strictly positive density. The mean of this density over the law of the state is the density $\mathfrak p(z,\cdot)$ of $\mathbf Z_1$, which gives \eqref{eq: positive density}.

    To prove the existence of an invariant distribution, let $(\mathbf Z_n)_{n \Ge 0}$ be the SG-UBU solution with $\mathbf X_0 = 0$ and $\mathbf V_0 \sim \mathcal N(0,I_d)$. By \eqref{eq: sgubu moments}, the mean of $|z|^4$ under the probability distribution
    \begin{equation}
        \frac1n \sum_{k=0}^{n-1} \mathrm{Law}(\mathbf Z_k) , \qquad n \Ge 1 ,
        \label{eq: averaged laws}
    \end{equation}
    is $\frac1n \sum_{k=0}^{n-1} \mathbb E |\mathbf Z_k|^4 \Le C d^2$. Hence the distributions \eqref{eq: averaged laws} are tight, and by Prokhorov's theorem, a subsequence converges weakly to a probability distribution $\pi_h$. This weak convergence gives \eqref{eq: invariant moments}:
    \begin{equation*}
        \int_{\mathbb R^{2d}} |z|^4 \,\pi_h(\D z) \Le \sup_{n \Ge 1} \frac1n \sum_{k=0}^{n-1} \mathbb E |\mathbf Z_k|^4 \Le C d^2 .
    \end{equation*}

    To prove that $\pi_h$ is invariant, let $\psi\in C(\mathbb R^{2d})$ be bounded. By the Markov property of the SG-UBU solution, the function $z \mapsto \mathbb E^z[\psi(\mathbf Z_1)]$ has the mean $\mathbb E[\psi(\mathbf Z_{k+1})]$ under $\mathrm{Law}(\mathbf Z_k)$. Hence 
    \begin{equation*}
        \biggl| \frac1n \sum_{k=0}^{n-1} \mathbb E [\psi(\mathbf Z_{k+1})] - \frac1n \sum_{k=0}^{n-1} \mathbb E [\psi(\mathbf Z_k)] \biggr| = \frac1n \Bigl| \mathbb E [\psi(\mathbf Z_n)] - \mathbb E [\psi(\mathbf Z_0)] \Bigr| \Le \frac{2 \norm{\psi}_\infty}{n} .
    \end{equation*}
    Since $z \mapsto \mathbb E^z[\psi(\mathbf Z_1)]$ is bounded, and continuous by \eqref{eq: feller}, letting $n \to \infty$ along the subsequence gives
    \begin{equation*}
        \int_{\mathbb R^{2d}} \mathbb E^z\bigl[ \psi(\mathbf Z_1) \bigr] \,\pi_h(\D z) = \int_{\mathbb R^{2d}} \psi(z) \,\pi_h(\D z)
    \end{equation*}
    for every bounded $\psi\in C(\mathbb R^{2d})$. Hence $\pi_h$ is an invariant distribution of SG-UBU.

    To prove the uniqueness, let $\pi_h'$ be an invariant distribution of SG-UBU. Since the law of $\mathbf Z_1$ for the initial state $\mathbf Z_0 = z$ has the density $\mathfrak p(z,\cdot)$, the invariant distributions $\pi_h$ and $\pi_h'$ have densities in the Lebesgue measure, written $\pi_h(z)$ and $\pi_h'(z)$, and for almost every $z' \in \mathbb R^{2d}$,
    \begin{equation}
        \pi_h(z') - \pi_h'(z') = \int_{\mathbb R^{2d}} \bigl( \pi_h(z) - \pi_h'(z) \bigr)\, \mathfrak p(z,z') \,\D z .
        \label{eq: density difference}
    \end{equation}
    Since $\int_{\mathbb R^{2d}} \mathfrak p(z,z') \,\D z' = 1$, integrating the absolute value of \eqref{eq: density difference} in $z'$ gives
    \begin{equation*}
        \int_{\mathbb R^{2d}} |\pi_h(z') - \pi_h'(z')| \,\D z' \Le \int_{\mathbb R^{2d}} \int_{\mathbb R^{2d}} |\pi_h(z) - \pi_h'(z)|\, \mathfrak p(z,z') \,\D z \,\D z' = \int_{\mathbb R^{2d}} |\pi_h(z) - \pi_h'(z)| \,\D z .
    \end{equation*}
    Hence this inequality is an equality, and by \eqref{eq: positive density}, either $\pi_h(z) \Ge \pi_h'(z)$ for almost every $z$, or $\pi_h(z) \Le \pi_h'(z)$ for almost every $z$. Then $\int_{\mathbb R^{2d}} \bigl( \pi_h(z) - \pi_h'(z) \bigr) \D z = 0$ gives $\pi_h = \pi_h'$.
\end{proof}

\subsection{Proof of Theorem~\ref{thm: bias}}
\label{app: bias}

In this subsection, $(\mathbf Z_n)_{n \Ge 0}$ is the SG-UBU solution with the initial state $\mathbf Z_0 \sim \pi_h$, and $Z_t$ is the solution to \eqref{eq: uld} with $Z_0 = \mathbf Z_0$. The random variable $\xi_0$ of the step from $\mathbf Z_0$ is the random variable $\xi$ of \eqref{eq: bias coefficient}. Since $\pi_h$ is invariant, $\mathbf Z_n \sim \pi_h$ for every $n \Ge 0$, so \eqref{eq: invariant moments} gives the moment bound \eqref{eq: sgubu moments} for $(\mathbf Z_n)_{n \Ge 0}$. Using the independence of $\mathbf Z_0$, $(\xi_n)_{n \Ge 0}$ and $(B_t)_{t \Ge 0}$, the proofs of Lemma~\ref{lem: local error} and Theorem~\ref{thm: sgubu} apply to $(\mathbf Z_n)_{n \Ge 0}$.

\begin{proof}[Proof of Theorem~\ref{thm: bias}]
    Since the state $\mathbf Z_k \sim \pi_h$, the random variable $\xi_k$ and the increments of $(B_t)_{t \Ge 0}$ on $[kh, (k+1)h]$ are independent, \eqref{eq: sgubu} and \eqref{eq: full gradient step} give
    \begin{equation}
        \mathrm{Law}\bigl( \mathbf Z_{k+1}, \bar{\mathbf Z}_{k+1} \bigr) = \mathrm{Law}\bigl( \mathbf Z_1, \bar{\mathbf Z}_1 \bigr) , \qquad k \Ge 0 .
        \label{eq: stationary pair}
    \end{equation}
    Since $\mathbb E [f(\mathbf Z_n)] = \pi_h(f)$, \eqref{eq: sgubu telescoping}, \eqref{eq: sgubu local part} and \eqref{eq: stationary pair} with $j = n-k-1$ give
    \begin{equation}
        \biggl| \pi_h(f) - \mathbb E [f(Z_{nh})] - \sum_{j=0}^{n-1} \mathbb E^{\pi_h}_\xi\Bigl[ (e^{jh\mathcal L}f)(\mathbf Z_1) - (e^{jh\mathcal L}f)(\bar{\mathbf Z}_1) \Bigr] \biggr| \Le C \norm{\nabla f}_\infty\, d h^2 , \qquad n \Ge 1 .
        \label{eq: stationary error}
    \end{equation}
    For $j \Ge 1$, \eqref{eq: sgubu inputs}--\eqref{eq: sgubu second order} hold with $k = 0$ and $jh$ in place of $t_0$, and $\mathcal T_1$ has mean zero. Hence \eqref{eq: time sums} gives
    \begin{equation*}
        \sum_{j=1}^\infty \Bigl| \mathbb E^{\pi_h}_\xi\Bigl[ (e^{jh\mathcal L}f)(\mathbf Z_1) - (e^{jh\mathcal L}f)(\bar{\mathbf Z}_1) \Bigr] \Bigr| \Le C \norm{\nabla f}_\infty\, \Sigma d h .
    \end{equation*}
    Since $Z_0 \sim \pi_h$ and $\pi$ is invariant under \eqref{eq: uld}, \eqref{eq: contraction} gives, with $\lambda$ as in Lemma~\ref{lem: contraction},
    \begin{equation}
        \bigl| \mathbb E [f(Z_t)] - \pi(f) \bigr| = \bigl| \pi_h(e^{t\mathcal L}f) - \pi(e^{t\mathcal L}f) \bigr| \Le C e^{-\lambda t} \norm{\nabla f}_\infty\, \mathcal W_1(\pi_h,\pi) , \qquad t \Ge 0 .
        \label{eq: stationary limit}
    \end{equation}
    Hence letting $n \to \infty$ in \eqref{eq: stationary error} gives
    \begin{equation}
        \biggl| \pi_h(f) - \pi(f) - \sum_{j=0}^\infty \mathbb E^{\pi_h}_\xi\Bigl[ (e^{jh\mathcal L}f)(\mathbf Z_1) - (e^{jh\mathcal L}f)(\bar{\mathbf Z}_1) \Bigr] \biggr| \Le C \norm{\nabla f}_\infty\, d h^2 .
        \label{eq: bias series}
    \end{equation}

    By \eqref{eq: hessian blocks} at $z = \bar{\mathbf Z}_1 + \theta ( \mathbf Z_1 - \bar{\mathbf Z}_1 )$ and \eqref{eq: time sums}, the coefficient $\mathcal A_h(f)$ of \eqref{eq: bias coefficient} satisfies
    \begin{equation}
        \biggl| \sum_{j=1}^\infty \mathbb E^{\pi_h}_\xi\Bigl[ (e^{jh\mathcal L}f)(\mathbf Z_1) - (e^{jh\mathcal L}f)(\bar{\mathbf Z}_1) \Bigr] - h\, \mathcal A_h(f) \biggr| \Le C \norm{\nabla f}_\infty\, \Sigma d h^2 .
        \label{eq: stationary sg sum}
    \end{equation}
    By \eqref{eq: bias series} and \eqref{eq: stationary sg sum}, we only need to prove
    \begin{equation}
        \Bigl| \mathbb E^{\pi_h}_\xi\bigl[ f(\mathbf Z_1) - f(\bar{\mathbf Z}_1) \bigr] \Bigr| \Le C \norm{\nabla f}_\infty\, \Sigma (1+\Sigma)\, d^2 h^2 .
        \label{eq: stationary last step}
    \end{equation}

    Note that $A \in \mathbb R^{d\times d}$ satisfies $A : \mathcal C(x) = \mathbb E_\xi\bigl[ A[\hat\eta(x;\xi), \hat\eta(x;\xi)] \bigr]$, so \eqref{eq: sg variance} and \eqref{eq: sg lipschitz} give, for every $x, x' \in \mathbb R^d$,
    \begin{equation}
        |A : \mathcal C(x)| \Le \norm{A}\, \Sigma d , \qquad
        \bigl| A : \bigl( \mathcal C(x) - \mathcal C(x') \bigr) \bigr| \Le 2 \norm{A}\, \Sigma d^{\frac12} |x - x'| .
        \label{eq: covariance bounds}
    \end{equation}
    Since $\xi$ is independent of $\mathscr U_{h/2} \mathbf Z_0$, \eqref{eq: last step} at $n = 1$, Taylor's formula for $f_{h}$ at $\mathscr B_h \mathscr U_{h/2} \mathbf Z_0$ and \eqref{eq: half step bounds} give
    \begin{equation}
        \biggl| \mathbb E^{\pi_h}_\xi\bigl[ f(\mathbf Z_1) - f(\bar{\mathbf Z}_1) \bigr] - \frac{h^2}2\, \mathbb E^{\pi_h}\Bigl[ \nabla^2_{vv}f_{h}\bigl( \mathscr B_h \mathscr U_{h/2} \mathbf Z_0 \bigr) : \mathcal C(\mathbf Y_0) \Bigr] \biggr| \Le C \norm{\nabla f}_\infty (\Sigma d)^{\frac32} h^2 .
        \label{eq: last step covariance}
    \end{equation}
    Here $\mathbb E^{\pi_h}$ is written for $\mathbb E^{\pi_h}_\xi$ when the random variable does not involve $\xi$. By \eqref{eq: half step bounds} and \eqref{eq: covariance bounds}, for every $z = (x,v)$ and $z' = (x',v')$ in $\mathbb R^{2d}$,
    \begin{equation}
        \Bigl| \nabla^2_{vv}f_{h}(z) : \mathcal C(x) - \nabla^2_{vv}f_{h}(z') : \mathcal C(x') \Bigr| \Le C h^{-1} \norm{\nabla f}_\infty\, \Sigma d\, |z - z'| .
        \label{eq: half step lipschitz}
    \end{equation}
    By \eqref{eq: U flow}, \eqref{eq: B flow} and \eqref{eq: displacement moments}, the state $\mathscr B_h \mathscr U_{h/2} \mathbf Z_0$, whose position is $\mathbf Y_0$, and the state $\mathbf Z_0\bigl( \frac h2 \bigr) = Z_{h/2}$ satisfy
    \begin{equation}
        \mathbb E^{\pi_h} \Bigl| \mathscr B_h \mathscr U_{h/2} \mathbf Z_0 - \mathbf Z_0\bigl( \tfrac h2 \bigr) \Bigr| \Le C d^{\frac12} h \Longrightarrow
        \mathbb E^{\pi_h} \bigl| \mathbf Y_0 - X_{h/2} \bigr| \Le C d^{\frac12} h .
        \label{eq: half step state}
    \end{equation}
    Since $\mathbf Z_n \sim \pi_h$, \eqref{eq: sgubu error} for $(\mathbf Z_n)_{n \Ge 0}$ and $h \Le 1$ give
    \begin{equation*}
        \bigl| \pi_h(f) - \mathbb E^{\pi_h} [f(Z_{nh})] \bigr| \Le C \norm{\nabla f}_\infty (1+\Sigma)\, d h , \qquad n \Ge 0 .
    \end{equation*}
    By \eqref{eq: stationary limit}, $\mathbb E^{\pi_h} [f(Z_{nh})] \to \pi(f)$ as $n \to \infty$, so
    \begin{equation*}
        \bigl| \pi_h(f) - \pi(f) \bigr| \Le C \norm{\nabla f}_\infty (1+\Sigma)\, d h .
    \end{equation*}
    Since $f \in C^1(\mathbb R^{2d})$ with $\norm{\nabla f}_\infty < \infty$ is arbitrary, Kantorovich--Rubinstein duality on $\mathbb R^{2d}$ gives
    \begin{equation}
        \mathcal W_1(\pi_h,\pi) \Le C (1+\Sigma)\, d h .
        \label{eq: phase space w1}
    \end{equation}
    The bounds \eqref{eq: half step lipschitz} and \eqref{eq: half step state} give
    \begin{equation}
        \Bigl| \mathbb E^{\pi_h}\Bigl[ \nabla^2_{vv}f_{h}\bigl( \mathscr B_h \mathscr U_{h/2} \mathbf Z_0 \bigr) : \mathcal C(\mathbf Y_0) \Bigr] - \mathbb E^{\pi_h}\Bigl[ \nabla^2_{vv}f_{h}( Z_{h/2} ) : \mathcal C( X_{h/2} ) \Bigr] \Bigr| \Le C \norm{\nabla f}_\infty\, \Sigma d^{\frac32} .
        \label{eq: half step shift}
    \end{equation}
    By mollification, \eqref{eq: stationary limit} also holds for every Lipschitz function $f$, with the Lipschitz constant of $f$ in place of $\norm{\nabla f}_\infty$. At $t = \frac h2$ and for the function $(x,v) \mapsto \nabla^2_{vv}f_{h}(x,v) : \mathcal C(x)$, with the Lipschitz constant $C h^{-1} \norm{\nabla f}_\infty\, \Sigma d$ of \eqref{eq: half step lipschitz}, \eqref{eq: stationary limit} and the bound \eqref{eq: phase space w1} give
    \begin{equation}
        \Bigl| \mathbb E^{\pi_h}\Bigl[ \nabla^2_{vv}f_{h}( Z_{h/2} ) : \mathcal C( X_{h/2} ) \Bigr] - \pi\bigl( \nabla^2_{vv}f_{h} : \mathcal C \bigr) \Bigr| \Le C h^{-1} \norm{\nabla f}_\infty\, \Sigma d\, \mathcal W_1(\pi_h,\pi) \Le C \norm{\nabla f}_\infty\, \Sigma (1+\Sigma)\, d^2 .
        \label{eq: half step exact}
    \end{equation}
    Since $d \Ge 1$, \eqref{eq: half step shift} and \eqref{eq: half step exact} give
    \begin{equation}
        \Bigl| \mathbb E^{\pi_h}\Bigl[ \nabla^2_{vv}f_{h}\bigl( \mathscr B_h \mathscr U_{h/2} \mathbf Z_0 \bigr) : \mathcal C(\mathbf Y_0) \Bigr] - \pi\bigl( \nabla^2_{vv}f_{h} : \mathcal C \bigr) \Bigr| \Le C \norm{\nabla f}_\infty\, \Sigma (1+\Sigma)\, d^2 .
        \label{eq: half step law}
    \end{equation}
    
    Under $\pi$, the velocity $v \sim \mathcal N(0,I_d)$ and $( v \otimes v - I_d ) : \mathcal C(x)$ has mean zero, hence \eqref{eq: sg variance} gives
    \begin{equation}
        \pi\Bigl( \bigl( ( v \otimes v - I_d ) : \mathcal C(x) \bigr)^2 \Bigr) \Le 2 \Sigma^2 d^2 .
        \label{eq: noise variance}
    \end{equation}
    Integration by parts in $v$, the Cauchy--Schwarz inequality, \eqref{eq: noise variance}, \eqref{eq: poincare full} and \eqref{eq: half step bounds} give
    \begin{equation}
        \bigl| \pi\bigl( \nabla^2_{vv}f_{h} : \mathcal C \bigr) \bigr| = \Bigl| \pi\Bigl( \bigl( f_{h} - \pi(f_{h}) \bigr)\, \bigl( v \otimes v - I_d \bigr) : \mathcal C(x) \Bigr) \Bigr| \Le C \norm{\nabla f}_\infty\, \Sigma d .
        \label{eq: half step mean}
    \end{equation}
    Since $d \Ge 1$ and $\Sigma^{\frac12} \Le 1 + \Sigma$, \eqref{eq: last step covariance}, \eqref{eq: half step law} and \eqref{eq: half step mean} give \eqref{eq: stationary last step}, which completes the proof of \eqref{eq: bias}.

    Finally, \eqref{eq: uniform} at the times $t = jh$, \eqref{eq: sgubu inputs} at $k = 0$ and \eqref{eq: time sums} give \eqref{eq: bias coefficient bound}.
\end{proof}

\subsection{Proof of Theorem~\ref{thm: coefficient limit}}
\label{app: coefficient limit}

\begin{proof}[Proof of Theorem~\ref{thm: coefficient limit}]
    In this proof, $\mathbf Z_0 \sim \pi_h$, and $\mathbf Z_1 = (\mathbf X_1, \mathbf V_1)$ is the state after one step of \eqref{eq: sgubu} from $\mathbf Z_0$, so that $\mathbf Z_1 \sim \pi_h$. The symbols $\mathbf Y_0$ and $\bar{\mathbf Z}_1$ are those of \eqref{eq: full gradient step}. By \eqref{eq: sg error} and $h \Le \frac1d \Le 1$,
    \begin{equation}
        |\mathbf Z_1 - \bar{\mathbf Z}_1| \Le 2h\, |\hat\eta(\mathbf Y_0;\xi)| .
        \label{eq: sg error size}
    \end{equation}
    By \eqref{eq: step size}, $h\Sigma \Le C$, so \eqref{eq: U flow}, \eqref{eq: B flow}, \eqref{eq: sg variance} and \eqref{eq: invariant moments} give
    \begin{equation}
        \bigl( \mathbb E^{\pi_h}_\xi |\mathbf X_1 - \mathbf Y_0|^2 \bigr)^{\frac12} \Le C d^{\frac12} h .
        \label{eq: position difference}
    \end{equation}
    With $\lambda$ as in Lemma~\ref{lem: contraction} and $h \Le 1$,
    \begin{equation}
        \sum_{j=1}^\infty h \min(1,jh)^{-1} e^{-\lambda jh} \Le C \Bigl( 1 + \log\frac1h \Bigr) .
        \label{eq: log sum}
    \end{equation}
    Next we prove \eqref{eq: coefficient limit} in four steps.\\[6pt]
    \textbf{Step 1:} Prove
    \begin{equation}
        \biggl| \mathcal A_h(f) - \frac h2 \sum_{j=1}^\infty \mathbb E^{\pi_h}\Bigl[ \nabla^2_{vv}(e^{jh\mathcal L}f)(\bar{\mathbf Z}_1) : \mathcal C(\mathbf Y_0) \Bigr] \biggr| \Le C \norm{\nabla f}_\infty (\Sigma d)^{\frac32} h \Bigl( 1 + \log\frac1h \Bigr) .
        \label{eq: limit theta}
    \end{equation}
    Since $\xi$ is independent of the pair $(\mathbf Y_0, \bar{\mathbf Z}_1)$, \eqref{eq: third uniform}, \eqref{eq: sg error size} and \eqref{eq: sg variance} give, for every $t > 0$,
    \begin{equation*}
    \begin{aligned}
        & \biggl| \mathbb E^{\pi_h}_\xi\biggl[ \int_0^1 (1-\theta)\, \nabla^2_{vv}(e^{t\mathcal L}f)\bigl( \bar{\mathbf Z}_1 + \theta ( \mathbf Z_1 - \bar{\mathbf Z}_1 ) \bigr)\bigl[ \hat\eta(\mathbf Y_0;\xi), \hat\eta(\mathbf Y_0;\xi) \bigr] \D\theta \biggr] \\
        & \qquad - \frac12\, \mathbb E^{\pi_h}\Bigl[ \nabla^2_{vv}(e^{t\mathcal L}f)(\bar{\mathbf Z}_1) : \mathcal C(\mathbf Y_0) \Bigr] \biggr| \Le C \min(1,t)^{-1} e^{-\lambda t} \norm{\nabla f}_\infty\, (\Sigma d)^{\frac32} h .
    \end{aligned}
    \end{equation*}
    Multiplying this bound by $h$ and summing over $t = jh$, $j \Ge 1$, with \eqref{eq: log sum} gives \eqref{eq: limit theta}.\\[6pt]
    \textbf{Step 2:} Prove
    \begin{equation}
    \begin{aligned}
        & \biggl| \frac h2 \sum_{j=1}^\infty \mathbb E^{\pi_h}\Bigl[ \nabla^2_{vv}(e^{jh\mathcal L}f)(\bar{\mathbf Z}_1) : \mathcal C(\mathbf Y_0) \Bigr] - \frac h2 \sum_{j=1}^\infty \pi_h\bigl( \nabla^2_{vv}(e^{jh\mathcal L}f) : \mathcal C \bigr) \biggr| \\
        & \qquad \Le C \norm{\nabla f}_\infty \Bigl( (\Sigma d)^{\frac32} h \Bigl( 1 + \log\frac1h \Bigr) + \Sigma d h \Bigr) .
    \end{aligned}
        \label{eq: limit state}
    \end{equation}
    For every $t > 0$, \eqref{eq: third uniform}, \eqref{eq: sg error size}, \eqref{eq: covariance bounds} and \eqref{eq: sg variance} give
    \begin{equation*}
        \Bigl| \mathbb E^{\pi_h}_\xi\Bigl[ \Bigl( \nabla^2_{vv}(e^{t\mathcal L}f)(\bar{\mathbf Z}_1) - \nabla^2_{vv}(e^{t\mathcal L}f)(\mathbf Z_1) \Bigr) : \mathcal C(\mathbf Y_0) \Bigr] \Bigr| \Le C \min(1,t)^{-1} e^{-\lambda t} \norm{\nabla f}_\infty\, (\Sigma d)^{\frac32} h ,
    \end{equation*}
    and \eqref{eq: uniform}, \eqref{eq: covariance bounds} and \eqref{eq: position difference} give
    \begin{equation*}
        \Bigl| \mathbb E^{\pi_h}_\xi\Bigl[ \nabla^2_{vv}(e^{t\mathcal L}f)(\mathbf Z_1) : \bigl( \mathcal C(\mathbf Y_0) - \mathcal C(\mathbf X_1) \bigr) \Bigr] \Bigr| \Le C \min(1,t)^{-\frac12} e^{-\lambda t} \norm{\nabla f}_\infty\, \Sigma d\, h .
    \end{equation*}
    Since $\mathbf Z_1 \sim \pi_h$, $\mathbb E^{\pi_h}_\xi\bigl[ \nabla^2_{vv}(e^{t\mathcal L}f)(\mathbf Z_1) : \mathcal C(\mathbf X_1) \bigr] = \pi_h\bigl( \nabla^2_{vv}(e^{t\mathcal L}f) : \mathcal C \bigr)$. Hence adding these two bounds gives,
    \begin{equation*}
        \Bigl| \mathbb E^{\pi_h}\Bigl[ \nabla^2_{vv}(e^{t\mathcal L}f)(\bar{\mathbf Z}_1) : \mathcal C(\mathbf Y_0) \Bigr] - \pi_h\bigl( \nabla^2_{vv}(e^{t\mathcal L}f) : \mathcal C \bigr) \Bigr| \Le C e^{-\lambda t} \norm{\nabla f}_\infty \Bigl( \min(1,t)^{-1} (\Sigma d)^{\frac32} h + \min(1,t)^{-\frac12}\, \Sigma d\, h \Bigr) .
    \end{equation*}
    Multiplying this bound by $\frac h2$ and summing over $t = jh$, $j \Ge 1$, with \eqref{eq: time sums} and \eqref{eq: log sum} gives \eqref{eq: limit state}.\\[6pt]
    \textbf{Step 3:} Prove
    \begin{equation}
        \frac h2 \sum_{j=1}^\infty \Bigl| \pi_h\bigl( \nabla^2_{vv}(e^{jh\mathcal L}f) : \mathcal C \bigr) - \pi\bigl( \nabla^2_{vv}(e^{jh\mathcal L}f) : \mathcal C \bigr) \Bigr| \Le C \norm{\nabla f}_\infty\, \Sigma (1+\Sigma)\, d^2 h \Bigl( 1 + \log\frac1h \Bigr) .
        \label{eq: limit law}
    \end{equation}
    For every $t > 0$, $z = (x,v)$ and $z' = (x',v')$ in $\mathbb R^{2d}$, \eqref{eq: uniform}, \eqref{eq: third uniform} and \eqref{eq: covariance bounds} give
    \begin{equation*}
        \Bigl| \nabla^2_{vv}(e^{t\mathcal L}f)(z) : \mathcal C(x) - \nabla^2_{vv}(e^{t\mathcal L}f)(z') : \mathcal C(x') \Bigr| \Le C \min(1,t)^{-1} e^{-\lambda t} \norm{\nabla f}_\infty\, \Sigma d\, |z - z'| .
    \end{equation*}
    Hence, at the times $t = jh$, Kantorovich--Rubinstein duality, \eqref{eq: phase space w1} and \eqref{eq: log sum} give \eqref{eq: limit law}.\\[6pt]
    \textbf{Step 4:} Prove
    \begin{equation}
        \biggl| \frac h2 \sum_{j=1}^\infty \pi\bigl( \nabla^2_{vv}(e^{jh\mathcal L}f) : \mathcal C \bigr) - \frac12\, \pi\bigl( \nabla^2_{vv}\phi : \mathcal C \bigr) \biggr| \Le C \norm{\nabla f}_\infty\, \Sigma d^2 h .
        \label{eq: limit riemann}
    \end{equation}
    Under $\pi$, the velocity $v\sim\mathcal N(0,I_d)$ is independent of the position, so integration by parts in $v$ gives
    \begin{equation}
        \pi\bigl( \nabla^2_{vv}(e^{t\mathcal L}f) : \mathcal C \bigr) = \pi\Bigl( (e^{t\mathcal L}f)\, \bigl( v \otimes v - I_d \bigr) : \mathcal C(x) \Bigr) , \qquad \forall t>0.
        \label{eq: gaussian parts}
    \end{equation}
    For $0 < s < t$, since $e^{s\mathcal L}f \in C^2(\mathbb R^{2d})$ by Corollary~\ref{cor: uniform}, It\^o's formula gives
    \begin{equation}
        e^{t\mathcal L}f - e^{s\mathcal L}f = \int_0^{t-s} e^{r\mathcal L}\bigl( \mathcal L e^{s\mathcal L}f \bigr) \D r .
        \label{eq: limit ito}
    \end{equation}
    By \eqref{eq: generator}, \eqref{eq: contraction} and \eqref{eq: uniform}, for every $z = (x,v) \in \mathbb R^{2d}$,
    \begin{equation*}
        \bigl| ( \mathcal L e^{s\mathcal L}f )(z) \bigr| \Le C \norm{\nabla f}_\infty e^{-\lambda s} \Bigl( |v| + |\nabla U(x)| + d \min(1,s)^{-\frac12} \Bigr) .
    \end{equation*}
    Integration by parts in $x$ and \eqref{eq: U bounds} give $\pi(|\nabla U|^2) = \pi(\Delta U) \Le M d$, and $\pi(|v|^2) = d$, so
    \begin{equation}
        \pi\bigl( ( \mathcal L e^{s\mathcal L}f )^2 \bigr)^{\frac12} \Le C \norm{\nabla f}_\infty e^{-\lambda s} \min(1,s)^{-\frac12}\, d .
        \label{eq: limit generator}
    \end{equation}
    Since $\pi$ is invariant, $\pi\bigl( (e^{r\mathcal L}\psi)^2 \bigr) \Le \pi(\psi^2)$ for every $\psi \in L^2(\pi)$ and $r \Ge 0$. Hence \eqref{eq: gaussian parts}--\eqref{eq: limit generator}, \eqref{eq: noise variance} and the Cauchy--Schwarz inequality give
    \begin{equation}
        \Bigl| \pi\bigl( \nabla^2_{vv}(e^{t\mathcal L}f) : \mathcal C \bigr) - \pi\bigl( \nabla^2_{vv}(e^{s\mathcal L}f) : \mathcal C \bigr) \Bigr| \Le C e^{-\lambda s} \min(1,s)^{-\frac12} \norm{\nabla f}_\infty\, \Sigma d^2 (t-s) , \qquad 0 < s < t .
        \label{eq: limit lipschitz}
    \end{equation}
    By \eqref{eq: poisson solution} and \eqref{eq: uniform}, $\nabla^2_{vv}\phi = \int_0^\infty \nabla^2_{vv}(e^{t\mathcal L}f) \,\D t$. Hence integrating \eqref{eq: limit lipschitz} at $t = jh$ over $s \in [(j-1)h, jh]$ and summing over $j \Ge 1$ gives \eqref{eq: limit riemann}.

    Since $d \Ge 1$, $h \Le 1$ and $\Sigma^{\frac12} \Le 1 + \Sigma$, adding \eqref{eq: limit theta}--\eqref{eq: limit riemann} gives \eqref{eq: coefficient limit}.
\end{proof}

\subsection{Proof of Corollary~\ref{lem: w1 limit}}
\label{app: w1 limit}

\begin{proof}[Proof of Corollary~\ref{lem: w1 limit}]
    Since $\mathfrak A \in L^2(\pi)$, the Cauchy--Schwarz inequality gives
    \begin{equation}
        \int_{\mathbb R^{2d}} \bigl( 1 + |z| \bigr)\, |\mathfrak A(z)| \,\pi(\D z) \Le \pi\bigl( ( 1 + |z| )^2 \bigr)^{\frac12}\, \pi(\mathfrak A^2)^{\frac12} < \infty .
        \label{eq: response moment}
    \end{equation}
    By \eqref{eq: response} with $f = 1$, for which $\phi = 0$ in \eqref{eq: poisson solution}, $\pi(\mathfrak A) = 0$. Hence, with \eqref{eq: response moment}, the signed measure $\mathfrak A \pi$ has zero total mass and a finite first moment, so $\norm{\mathfrak A \pi}_{\mathrm{KR}}$ is defined.

    Let $h \Le \min\{ h_0, \frac1d \}$ with $h_0$ in \eqref{eq: step size}, and let $f \in C^1(\mathbb R^{2d})$ with $\norm{\nabla f}_\infty \Le 1$. By \eqref{eq: response} and \eqref{eq: bias expansion},
    \begin{equation}
        \Bigl| \pi_h(f) - \pi(f) - h\, \pi(\mathfrak A f) \Bigr| \Le C (1+\Sigma)^2\, d^2 h^2 \Bigl( 1 + \log\frac1h \Bigr) .
        \label{eq: w1 expansion}
    \end{equation}
    By mollification and \eqref{eq: response moment}, \eqref{eq: w1 expansion} holds for every function $f$ on $\mathbb R^{2d}$ with Lipschitz constant at most $1$. With the definition of $\norm{\cdot}_{\mathrm{KR}}$, taking the supremum over these functions $f$ gives
    \begin{equation*}
        \Bigl| \mathcal W_1(\pi_h,\pi) - h\, \norm{\mathfrak A \pi}_{\mathrm{KR}} \Bigr| \Le C (1+\Sigma)^2\, d^2 h^2 \Bigl( 1 + \log\frac1h \Bigr) .
    \end{equation*}
    Dividing this bound by $h$ and letting $h \to 0$ gives \eqref{eq: w1 limit}.
\end{proof}

\end{document}